\documentclass[reqno]{amsart}

\usepackage{tikz-cd}
\usetikzlibrary{knots}

\usepackage{fullpage,dsfont}

\usepackage{hyperref} % this should come before amsrefs

\usepackage{parskip}
\makeatletter % need this to avoid the conflict between amsthm and parskip
\def\thm@space@setup{%
 \thm@preskip=\parskip \thm@postskip=0pt
}
\def\thm@remark{%
  \thm@headfont{\itshape}%
  \normalfont % body font
  \thm@preskip\parskip \thm@postskip=0pt
}
\makeatother

\usepackage[nobysame,alphabetic,initials,msc-links]{amsrefs}

\usepackage[final]{pdfpages}

\usepackage{multirow}

\usepackage{array}

\usepackage{spath3}

\usepackage{float}

\DefineSimpleKey{bib}{how}
\DefineSimpleKey{bib}{mrclass}
\DefineSimpleKey{bib}{mrnumber}
\DefineSimpleKey{bib}{fjournal}
\DefineSimpleKey{bib}{mrreviewer}

\renewcommand{\PrintDOI}[1]{%
  \href{http://dx.doi.org/#1}{{\tt DOI:#1}}%
}
\renewcommand{\eprint}[1]{#1}
\BibSpec{book}{%
    +{}  {\PrintPrimary}                {transition}
    +{.} { \PrintDate}                  {date}
    +{.} { \textit}                     {title}
    +{.} { }                            {part}
    +{:} { \textit}                     {subtitle}
    +{,} { \PrintEdition}               {edition}
    +{}  { \PrintEditorsB}              {editor}
    +{,} { \PrintTranslatorsC}          {translator}
    +{,} { \PrintContributions}         {contribution}
    +{,} { }                            {series}
    +{,} { \voltext}                    {volume}
    +{,} { }                            {publisher}
    +{,} { }                            {organization}
    +{,} { }                            {address}
    +{,} { }                            {status}
    +{,} { \PrintDOI}                   {doi}
    +{,} { \PrintISBNs}                 {isbn}
    +{}  { \parenthesize}               {language}
    +{}  { \PrintTranslation}           {translation}
    +{;} { \PrintReprint}               {reprint}
    +{.} { }                            {note}
    +{.} {}                             {transition}
    +{}  {\SentenceSpace \PrintReviews} {review}
}
\BibSpec{article}{%
    +{}  {\PrintAuthors}                {author}
    +{,} { \textit}                     {title}
    +{.} { }                            {part}
    +{:} { \textit}                     {subtitle}
    +{,} { \PrintContributions}         {contribution}
    +{.} { \PrintPartials}              {partial}
    +{,} { }                            {journal}
    +{}  { \textbf}                     {volume}
    +{}  { \PrintDatePV}                {date}
    +{,} { \issuetext}                  {number}
    +{,} { \eprintpages}                {pages}
    +{,} { }                            {status}
    +{,} { \PrintDOI}                   {doi}
    +{,} { \eprint}                     {eprint}
    +{}  { \parenthesize}               {language}
    +{}  { \PrintTranslation}           {translation}
    +{;} { \PrintReprint}               {reprint}
    +{.} { }                            {note}
    +{.} {}                             {transition}
    +{}  {\SentenceSpace \PrintReviews} {review}
}
\BibSpec{collection.article}{%
    +{}  {\PrintAuthors}                {author}
    +{,} { \textit}                     {title}
    +{.} { }                            {part}
    +{:} { \textit}                     {subtitle}
    +{,} { \PrintContributions}         {contribution}
    +{,} { \PrintConference}            {conference}
    +{}  {\PrintBook}                   {book}
    +{,} { }                            {booktitle}
    +{,} { \PrintDateB}                 {date}
    +{,} { pp.~}                        {pages}
    +{,} { }                            {publisher}
    +{,} { }                            {organization}
    +{,} { }                            {address}
    +{,} { }                            {status}
    +{,} { \PrintDOI}                   {doi}
    +{,} { \eprint}                     {eprint}
    +{}  { \parenthesize}               {language}
    +{}  { \PrintTranslation}           {translation}
    +{;} { \PrintReprint}               {reprint}
    +{.} { }                            {note}
    +{.} {}                             {transition}
    +{}  {\SentenceSpace \PrintReviews} {review}
}

\usepackage{amssymb, amsfonts, amsxtra, amsmath}
\usepackage{mathrsfs}
\usepackage{mathdots}
\usepackage{wasysym}
\usepackage[all]{xy}
\usepackage{bbm}
\usepackage{calc}
\usepackage{accents}
\usepackage{enumitem}
\usepackage{fouridx}

\numberwithin{equation}{section}

\newtheorem{Theorem}{Theorem}[section]
\newtheorem*{Theorem*}{Theorem}
\newtheorem{Def}[Theorem]{Definition}
\newtheorem{Lem}[Theorem]{Lemma}
\newtheorem{Prop}[Theorem]{Proposition}
\newtheorem{Cor}[Theorem]{Corollary}

\newtheorem{Rem}[Theorem]{Remark}

\newcommand\bp{\begin{proof}}
\newcommand\ep{\end{proof}}

\mathchardef\mhyph="2D

\DeclareMathOperator{\End}{\mathrm{End}}

\DeclareMathOperator{\fin}{\mathrm{f}}

\DeclareMathOperator{\Fun}{\mathrm{Fun}}

\DeclareMathOperator{\Hom}{\mathrm{Hom}}
\DeclareMathOperator{\id}{\mathrm{id}}

\DeclareMathOperator{\Lin}{\mathrm{Lin}}

\DeclareMathOperator{\Ker}{\mathrm{Ker}}

\DeclareMathOperator{\can}{\mathrm{can}}

\newcommand{\ev}{\mathrm{ev}}
\newcommand{\cop}{\mathrm{cop}}
\newcommand{\op}{\mathrm{op}}

\newcommand{\wM}{\widetilde{M}}

\newcommand{\msB}{\mathscr{B}}

\newcommand{\msD}{\mathscr{D}}

\newcommand{\msF}{\mathscr{F}}

\newcommand{\msH}{\mathscr{H}}
\newcommand{\msI}{\mathscr{I}}

\newcommand{\msU}{\mathscr{U}}

\newcommand{\msX}{\mathscr{X}}

\newcommand{\mcD}{\mathcal{D}}

\newcommand{\Hsp}{\mathcal{H}}

\newcommand{\mbr}{\mathbf{r}}

\newcommand{\mbs}{\mathbf{s}}
\newcommand{\mbt}{\mathbf{t}}

\newcommand{\mbunit}{\mathbf{1}}

\newcommand{\C}{\mathbb{C}}

\newcommand{\Q}{\mathbb{Q}}

\newcommand{\Z}{\mathbb{Z}}

\newcommand{\opp}{\mathrm{op}}

\newcommand{\Span}{\mathrm{Span}}

\newcommand{\Grt}[3]{#1{\tiny {\begin{pmatrix} #2\\#3\end{pmatrix}}}}

\newcommand{\Unit}{\mathbf{1}}

\newcommand{\UnitC}[2]{\Grt{\mathbf{1}}{#1}{#2}} 
\newcommand{\UnitCC}[2]{\Grt{\check{\mathbf{1}}}{#1}{#2}}

\newcommand{\Gr}[5]{\fourIdx{#2}{#4}{#3}{#5}{#1}}%TODO: better typesetting
\newcommand{\emdash}{\;\textrm{---}\;}

\title{Partial $*$-algebraic quantum groups and Drinfeld doubles of partial compact quantum groups}
\author{K. De Commer and J. Konings}
\address{Vrije Universiteit Brussel}
\email{kenny.de.commer@vub.be}
\email{johan.konings@vub.be}

\begin{document}
\maketitle

\begin{abstract}
We introduce a notion of \emph{partial algebraic quantum group}. This is an important special case of a weak multiplier Hopf algebra with integrals, as introduced in the work of Van Daele and Wang. At the same time, it generalizes the notion of partial compact quantum group as introduced by De Commer and Timmermann. As an application, we show that the Drinfeld double of a partial compact quantum group can be defined as a partial $*$-algebraic quantum group. 
\end{abstract}

\section*{Introduction}

Tensor C$^*$-categories are important mathematical structures arising in many different contexts, such as for example the representation theory of compact quantum groups, subfactor theory and conformal field theory. In case they are closed under subobjects, have an irreducible unit and are rigid (i.e.\ any object has a dual), we refer to them as \emph{unitary} tensor C$^*$-categories (following the terminology of \cite{Pen20}).

Although many examples of unitary tensor C$^*$-categories arise from compact quantum groups, not all of them can be described in this way. However, continuing the work initiated by Hayashi \cites{Hay93, Hay99}, it was shown in \cite{DCT15} that any unitary tensor C$^*$-category can always be realized as the representation theory of a \emph{partial compact quantum group}. Such partial compact quantum groups can be seen as quantum groupoids with a classical set of objects. In a sense, since the theory of partial compact quantum groups is quite similar to that of compact quantum groups, this provides a more direct way to use well-known constructions from the theory of compact quantum groups in the general theory of unitary tensor C$^*$-categories. For the case of tensor C$^*$-categories with finitely many objects, see e.g.\ \cite{NiV02}.  

In this paper, we will develop a more general theory of \emph{partial algebraic quantum groups}, offering a partial version of the theory developed in \cite{VDae98}. In fact, our objects can be described through the general theory of weak multiplier bialgebras \cite{BG-TL-C15}, weak multiplier Hopf algebras (with integrals) \cites{VDW13,VDW20,VDW17} and multiplier Hopf algebroids with integrals \cites{Tim16,Tim17} (see \cite[Section 6]{TVDW22} for some of the history on these structures). These latter theories start out from a total algebra with a comultiplication, and from their axioms construct a \emph{base algebra} which can be interpreted as the algebra of functions on the `object set' of the corresponding quantum groupoid. In our setting, we \emph{start} with a given base algebra, which will simply be the (classical) function algebra on a given set. We believe that this special case is worth developing independently, as it allows to bypass some of the technicalities of the more general setting, as it provides a more direct entry point into the theory of weak multiplier Hopf algebras and quantum groupoids, and as this case is then already sufficiently general to cover many of the relevant applications (maybe up to thickening the base algebra through a Morita equivalence). However, once we have established our structures as special cases of the above theories, we will borrow results from the previous articles whenever convenient, as to expedite our getting to the main results of the paper. 

The main result of this paper will be the construction of the \emph{Drinfeld double} for a partial compact quantum group. Such a Drinfeld double has been constructed for finite-dimensional weak Hopf algebras \cite{Nen02}, but as far as we are aware this construction has not been considered yet in the setting of  weak (multiplier) Hopf algebras. Applications of this construction to the theory of tensor C$^*$-categories will appear elsewhere.

%As an application of this construction, we show that the tube $*$-algebra of a unitary tensor C$^*$-category is Morita equivalent to the Drinfeld double of some (canonical) partial compact quantum group. This Morita equivalent copy then possesses extra structure, such as the presence of a comultiplication, revealing in a more concrete way the monoidal structure of the representation category of the tube $^*$-algebra, which is simply the center of the original unitary tensor C$^*$-category \cite{PV15,NY16,GJ16}. 

Our paper is structured as follows. 

In the \emph{first section}, we treat in detail various aspects of non-unital algebras with a grading. Using this formalism, we then introduce the notion of a partial regular multiplier Hopf algebra, and show how the latter fits into the general formalism of weak multiplier bialgebras (as developed by B\"{o}hm, Gomez-Torrecillas and Lopez-Centella) and weak multiplier Hopf algebras (as developed by Van Daele and Wang). 

In the \emph{second section}, we describe the structure of left and right invariant functionals on partial algebraic quantum groups, and describe the associated duality theory. We then introduce partial $*$-algebraic quantum groups, where a compatible $*$-structure is present, and show that the resulting $*$-algebras allow for a good representation theory on Hilbert spaces. We note that most of the contents of this section are obtained from directly applying the results in \cites{VDW13,VDW17} and \cite{Tim16}, but on the one hand it is worthwhile to spell out the details concretely for future reference, and on the other hand we prove strictly stronger results in the case of partial $*$-algebraic quantum groups, proving for example complete diagonalizability of all the structure maps. 

In the \emph{third section}, we introduce the notion of a Drinfeld double for a partial Hopf algebra with invariant integrals, and we show that the resulting object is a partial algebraic quantum group. We then show the correspondence between representations of the Drinfeld double and Yetter-Drinfeld modules for the original partial Hopf algebra. 

%In the final \emph{fourth section}, we recall from \cite{DCT15} the tensor C$^*$-category associated to a partial compact quantum group. We then show that the tube $^*$-algebra associated to this tensor C$^*$-category is Morita equivalent (in the setting of $*$-algebras) to the Drinfeld double of the partial compact quantum group. 

\thanks{\emph{Acknowledgments}: The work of K.~De Commer was supported by the FWO grant G032919N.} We thank the referee for their valuable comments. 

\section{Regular partial multiplier Hopf algebras}

\subsection{$I$-partial algebras}
Let $I$ be a set. Let $A$ be a vector space (over $\C$) with an $I^2$-grading. We write this as 
\[
A = \underset{r,s \in I}{\oplus}\; \Gr{A}{}{}{r}{s}.
\]
We write $A_s = \underset{r}{\oplus}\, \Gr{A}{}{}{r}{s}$ and $\Gr{A}{}{}{r}{} = \underset{s}{\oplus}\, \Gr{A}{}{}{r}{s}$.

We next assume that $A$ has an (associative) algebra structure $m$, which we write
\[
m: A \times A \rightarrow A,\qquad (a,b)\mapsto a\cdot b = ab. 
\]
For $V,W\subseteq A$ we write
\[
V\cdot W = \left\{\sum_i v_iw_i \mid v_i\in V,w_i\in W\right\}.
\]

\begin{Def}
We call $A$ an \emph{$I$-partial algebra} if $m$ has the following compatibility with the grading: $\Gr{A}{}{}{r}{s} \cdot \Gr{A}{}{}{s'}{t} = \{0\}$ for $s\neq s'$ and 
\[
\Gr{A}{}{}{r}{s} \cdot \Gr{A}{}{}{s}{t} \subseteq \;\Gr{A}{}{}{r}{t}. 
\]
\end{Def}

Note that $\Gr{A}{}{}{}{r}$ and $\Gr{A}{}{}{r}{}$ are subalgebras. When referring to $A$ as an algebra with the grading forgotten, we will call $A$ the \emph{total algebra}. 

We may equivalently view the $\Gr{A}{}{}{r}{s}$ as morphism spaces of a small $\C$-linear category with object set $I$, \emph{except} for the fact that we are not imposing the existence of identity morphisms. The lack of the latter has to be compensated by requiring nice regularity assumptions of the $I$-partial algebra. More precisely, the natural conditions are to impose regularity conditions on the algebras $\Gr{A}{}{}{r}{r}$, \emph{as well as} on the morphism spaces $\Gr{A}{}{}{r}{s}$ as $\Gr{A}{}{}{r}{r}$-$\Gr{A}{}{}{s}{s}$-bimodules. This leads to the following conditions which we now introduce. 

Recall (see e.g.\ the appendix of \cite{VDae94}) that the \emph{multiplier algebra} $M(A)$ of a (not necessarily unital) algebra $A$ is defined as the subalgebra of couples $m = (\lambda,\rho) \in \End(A)\oplus \End(A)^{\opp}$ such that, writing $ma = \lambda(a)$ and $am = \rho(a)$ for all $a\in A$, one has
\[
a(mb) = (am)b,\qquad \forall a,b\in A.  
\]
We endow $M(A)$ with the \emph{strict topology}, whereby a net $m_{\alpha}$ converges to $m$ if and only if for all $a\in A$ one has that $m_{\alpha}a = ma$ and $am_{\alpha} = am$ eventually.

If $A$ is an $I$-partial algebra, we can for example define idempotent and orthogonal elements $\Unit(r)\in M(A)$, the \emph{base units}, uniquely determined by  
\[
a \Unit(s) = \delta_{s,s'}a,\qquad \Unit(r)a = \delta_{r,r'}a,\quad a\in \Gr{A}{}{}{r'}{s'}.
\]

Recall that an algebra $A$ is called \emph{left non-degenerate} if $a\in A$ and $aA=0$ implies $a=0$. Similarly, an algebra $A$ is called \emph{right non-degenerate} if $a\in A$ and $Aa=0$ implies $a=0$. One calls $A$ \emph{non-degenerate} if it is left and right non-degenerate.

\begin{Def}
We call an $I$-partial algebra $A$ \emph{non-degenerate} if for each $r\in I$ the algebra $\Gr{A}{}{}{}{r}$ is left non-degenerate and the algebra $\Gr{A}{}{}{r}{}$ is right non-degenerate for any $r\in I$.
\end{Def}

This is easily seen to be equivalent with the following condition: if $a\in \Gr{A}{}{}{r}{s}$ and $a\Gr{A}{}{}{s}{s}=0$ or $\Gr{A}{}{}{r}{r}a =0$, then $a=0$, i.e.\ each $\Gr{A}{}{}{r}{s}$ is non-degenerate as a left $\Gr{A}{}{}{r}{r}$-module and as a right $\Gr{A}{}{}{s}{s}$-module (following the terminology in \cite{VDae08}).

\begin{Lem}\label{LemNonDeg}
If $A$ is non-degenerate as a partial algebra, then the total algebra $A$ is also non-degenerate. 
\end{Lem} 
\begin{proof}
Assume that $A$ is non-degenerate as a partial algebra, and assume $a = \sum_r a_r$ with $a_r \in \Gr{A}{}{}{}{r}$ and $aA=0$. Then in particular for any $r \in I$ we have
\[
0 = a\cdot \Gr{A}{}{}{r}{r}=a_r \cdot \Gr{A}{}{}{r}{r} = a_r \cdot \Gr{A}{}{}{}{r},
\]
so $a_r =0$ as $\Gr{A}{}{}{}{r}$ is left non-degenerate. Hence $a=0$, and $A$ is left non-degenerate. Similarly, one shows $A$ is right non-degenerate.
\end{proof} 

\begin{Rem}\label{RemNonDeg}
The converse statement does not hold. Indeed, consider  the $\Q$-partial algebra $A$ of strictly upper triangular matrices in $M_{\Q}(\C)$ with finitely many non-zero entries, with $\Gr{A}{}{}{r}{s}$ spanned by the matrix unit $e_{rs}$ if $r<s$ and $\Gr{A}{}{}{r}{s}=0$ if $r\geq s$. Then it is easily seen that $A$ is a degenerate $\Q$-partial algebra whose total algebra is non-degenerate. 
\end{Rem}

By Lemma \ref{LemNonDeg}, it follows that we can identify $A \subseteq M(A)$ if $A$ is a non-degenerate $I$-partial algebra. 

Recall that an algebra $A$ is called \emph{idempotent} if $A\cdot A =A$.

\begin{Def}
We call an $I$-partial algebra $A$ \emph{idempotent} if each $\Gr{A}{}{}{}{r}$ and $\Gr{A}{}{}{r}{}$ is idempotent. 
\end{Def}

This is easily seen to be equivalent with the following condition: $\Gr{A}{}{}{r}{s}\Gr{A}{}{}{s}{s} = \Gr{A}{}{}{r}{s}$ and $\Gr{A}{}{}{r}{r}\Gr{A}{}{}{r}{s} = \Gr{A}{}{}{r}{s}$ for all $r,s$, which means respectively that $\Gr{A}{}{}{r}{s}$ is unital as a right $\Gr{A}{}{}{s}{s}$-module and as a left $\Gr{A}{}{}{r}{r}$-module (following the terminology of \cite{DVDZ99}).

\begin{Lem} 
If $A$ is idempotent as a partial algebra, then the total algebra $A$ (forgetting the grading) is also idempotent. 
\end{Lem}
\begin{proof}
This is immediate. 
\end{proof}

\begin{Rem}
The converse statement does not hold, by the same example as in Remark \ref{RemNonDeg}. To obtain an example of a \emph{non-degenerate} non-idempotent partial algebra, yet whose total algebra \emph{is} idempotent, we modify the above example by putting $\Gr{A}{}{}{r}{s} = \C[x]$ for $r<s$ and $\Gr{A}{}{}{r}{s} = x\C[x]$ for $r \leq s$. On the other hand, note that if $A$ is an idempotent partial algebra whose total algebra is non-degenerate, then $A$ is a non-degenerate partial algebra: if $a \in \Gr{A}{}{}{}{r}$ and $a \Gr{A}{}{}{}{r} = 0$, then by partial idempotency 
\[
0 = a  \cdot \Gr{A}{}{}{}{r}\cdot \Gr{A}{}{}{r}{} =  a\cdot\Gr{A}{}{}{r}{}\cdot \Gr{A}{}{}{r}{} = a \cdot \Gr{A}{}{}{r}{} = aA,
\]
hence $a=0$ by total non-degeneracy. Similarly for the right non-degeneracy of the $\Gr{A}{}{}{r}{}$.
\end{Rem}

Note that eventually, the algebras we are interested in will have local units (in their total algebra), which will immediately imply their non-degeneracy and idempotency as partial algebras. Nevertheless, we believe it is important conceptually to note the difference between the notions of partial non-degeneracy/idempotency and total non-degeneracy/idempotency in the context of partial algebras (i.e.\ small linear categories without identity morphisms).

In what follows, we will work also with $I^2$-partial algebras, where we view elements of $I^2 = I\times I$ as columns $\begin{pmatrix} r \\ t\end{pmatrix}$. In this case, we write the base units as $\UnitC{r}{t}$ and 
\[
A =\underset{r,s,t,u}{\oplus}  \UnitC{r}{t}A\UnitC{s}{u}=  \underset{r,s,t,u}{\oplus}  \Gr{A}{r}{s}{t}{u}. 
\]
We can further group the base units $\UnitC{r}{t}$ together as
\[
\Unit^r = \sum_t \UnitC{r}{t} \in M(A),\qquad \Unit_t = \sum_r \UnitC{r}{t} \in M(A),
\]
where the sums converge for the strict topology on $M(A)$. We then also write 
\[
A^s = \underset{r,t,u}{\oplus}\, \Gr{A}{r}{s}{t}{u}
\]
etc.

\subsection{Tensor product of $I^2$-partial algebras}

Let $A$ and $B$ be $I^2$-partial algebras. We define 
\[
A \otimes^I B := \underset{r,s,t,u,v,w}{\oplus}\; \Gr{A}{r}{s}{t}{u}\otimes \Gr{B}{t}{u}{v}{w} \subseteq A \otimes B,
\]
which inherits the tensor product algebra structure. It becomes an $I^2$-partial algebra when putting 
\[
\Gr{(A\otimes^I B)}{r}{s}{v}{w} = \underset{t,u}{\oplus} \;\Gr{A}{r}{s}{t}{u} \otimes \Gr{B}{t}{u}{v}{w}. 
\]
Note that we can write 
\[
A \otimes^I B = E(A\otimes B)E
\]
where $E \in M(A\otimes B)$ is the idempotent
\[
E =\sum_{t} \Unit_{t} \otimes \Unit^t. 
\]
This sum and each of its summands have an obvious meaning as elements of $M(A\otimes B)$, and the sum converges for the strict topology on $M(A\otimes B)$.

We also have a natural homomorphism
\begin{equation}\label{EqQuotMap}
EM(A\otimes B)E \rightarrow M(A\otimes^I B)
\end{equation}
by restricting the action of multipliers. 

\begin{Lem}\label{LemTenProdND}
Let $A$ and $B$ be non-degenerate and idempotent $I^2$-partial algebras. Then $A\otimes^I B$ is a non-degenerate and idempotent $I^2$-partial algebra.
\end{Lem}
When considering a tensor $\sum_i a_i \otimes b_i \in \Gr{(A\otimes^I B)}{r}{s}{v}{w}$, we will tacitly assume that the simple tensors in the expression lie in some component $\Gr{A}{r}{s}{t}{u} \otimes \Gr{B}{t}{u}{v}{w}$.
\begin{proof}
Assume that $\sum_i a_i\otimes b_i \in \Gr{(A\otimes^I B)}{}{s}{}{w}$ and
\[
\sum_i a_i c\otimes b_id =0,\qquad \forall t,u\in I, \forall c \in \Gr{A}{}{s}{t}{u}, \forall d\in \Gr{B}{t}{u}{}{w}.
\]
Then also
\[
\sum_i a_i c\otimes b_id =0,\qquad \forall u\in I,\forall c\in A_u^s,\forall d\in B^u_w.
\]
Hence for $\omega$ a linear functional on $B$, we have for any fixed $u\in I$ and $d \in B^u_w$ that
\[
\sum_i \omega(b_id)a_i c=0,\qquad \forall   c\in A_u^s,
\]
so $\sum_i \omega(b_id)a_i \Unit_u=0$ by left non-degeneracy of $\Gr{A}{}{s}{}{u}$. Hence $\sum_i a_i\Unit_u \otimes b_id = 0$ for all $u\in I$ and $d\in B^u_w$. By symmetry, also $\sum_i a_i\Unit_u \otimes b_i\Unit^u = 0$ for all $u\in I$. As $u$ was arbitrary, we conclude that $\sum_i a_i\otimes b_i =0$. This shows that $\Gr{(A\otimes^I B)}{}{s}{}{w}$ is left non-degenerate. Right non-degeneracy of $\Gr{(A\otimes^I B)}{s}{}{w}{}$ is checked similarly.

To see that $A\otimes^I B$ is idempotent, choose $a \in \Gr{A}{}{s}{t}{u}$ and $b\in \Gr{B}{t}{u}{}{w}$. By idempotency of $A$ and $B$, we can then choose $c_i,d_i \in  \Gr{A}{}{s}{}{u}$ and $e_j,f_j \in \Gr{B}{}{u}{}{w}$ such that 
\[
\sum_i c_id_i = a,\qquad \sum_j e_jf_j = b. 
\]
With 
\[
x_{ij} = \sum_{ij} \Unit_tc_i \otimes \Unit^te_j,\qquad y_{ij} = \sum_{ij} \Unit_u d_i\otimes \Unit^uf_j,
\]
we have $x_{ij},y_{ij} \in \Gr{(A\otimes^I B)}{}{s}{}{w}$ with
\[
a\otimes b = \sum_{ij} x_{ij}y_{ij}. 
\]
The other idempotency condition is checked similarly.
\end{proof} 

\begin{Lem}\label{LemEqQuotMap}
Assume that $A,B$ are non-degenerate and idempotent $I^2$-partial algebras. Then the map \eqref{EqQuotMap} is an isomorphism of algebras.
\end{Lem}
\begin{proof}
Assume that $x\in A\otimes B$. Then by the same proof as for Lemma \ref{LemTenProdND}, we see that the non-degeneracy conditions on the multiplication of $A,B$ imply that
\begin{equation}\label{EqIdNonDeg}
\left(\forall y\in A\otimes^I B: xy = 0\right) \Rightarrow xE =0, \qquad
\left(\forall y\in A\otimes^I B: yx = 0\right) \Rightarrow Ex =0.
\end{equation}
To see then that \eqref{EqQuotMap} is injective, pick $m \in EM(A\otimes B)E$ and assume $m(A\otimes^I B) =0$. With $z \in A\otimes B$ and $x = zm \in A\otimes B$, it then follows from \eqref{EqIdNonDeg} that $zm =xE =0$. As $z \in A\otimes B$ was arbitrary, we conclude $m=0$ (by non-degeneracy of the total algebra $A\otimes B$). This proves injectivity.

To show the surjectivity, consider the subalgebra 
\[
A\,{}^I\!\otimes^IB = \underset{r,t,v}{\oplus} (\Gr{A}{r}{r}{t}{t}\otimes \Gr{B}{t}{t}{v}{v}) \subseteq A\otimes^I B. 
\]
Then similarly again as in the proof for  Lemma \ref{LemTenProdND}, the idempotency condition on $A,B$ implies that 
\begin{equation}\label{EqIdempotentnondeg}
(A\otimes B)(A\,{}^I\!\otimes^IB) = (A\otimes B)E, \qquad 
(A\,{}^I\!\otimes^IB)(A\otimes B) = E(A\otimes B).
\end{equation}
Hence if $m$ is a multiplier of $A\otimes^I B$, we can define, unambigously by \eqref{EqIdNonDeg}, on $A \otimes B$ the left multiplier
\[
\widetilde{m}x := \sum_i (my_i)z_i,\qquad x\in A\otimes B,\quad Ex = \sum_i y_iz_i\textrm{ for }y_i \in A\,{}^I\!\otimes^IB, z_i \in A\otimes B,
\]
and right multiplier
\[
x\widetilde{m} := \sum_i z_i(y_im),\qquad x\in A\otimes B,\quad xE = \sum_i z_iy_i\textrm{ for }y_i \in A\,{}^I\!\otimes^IB, z_i \in A\otimes B.
\]
It is then easy to check that $\widetilde{m}$ defines a multiplier in $M(A\otimes B)$, that $\widetilde{m}\in EM(A\otimes B)E$ and that $\widetilde{m}x = mx$ and $x\widetilde{m} = xm$ for $x\in A\otimes^IB \subseteq A\otimes B$, showing that the map \eqref{EqQuotMap} is surjective. 
\end{proof}

We will also need a particular subalgebra of the multiplier algebra of a tensor product. 

\begin{Def}\label{DefRestrMult}
Let $A,B$ be two non-degenerate, idempotent $I^2$-partial algebras. The \emph{restricted} multiplier algebra of the tensor product is defined as the subalgebra
\[
\widetilde{M}(A\otimes^I B) \subseteq M(A\otimes^I B)
\]
consisting of all $x\in  M(A\otimes^I B)$ such that, for any $a\in A,b\in B$, all of the following expressions
\[
x(a\otimes 1),\quad x(1\otimes b),\quad (a\otimes 1)x,\quad (1\otimes b)x
\]
lie in $A\otimes B$. 
\end{Def}
Here the $a\otimes 1$ etc.\ can easily be interpreted as elements of $M(A\otimes B)$, and we are viewing 
$x\in M(A\otimes^IB)\subseteq M(A\otimes B)$ by Lemma \ref{LemEqQuotMap}. In other words, we have 
\begin{equation}\label{EqTildeTens}
\widetilde{M}(A\otimes^I B) = \widetilde{M}(A\otimes B) \cap M(A\otimes^I B),
\end{equation}
where the restricted multiplier algebra $\widetilde{M}(A\otimes B)$ is defined exactly as in Definition \ref{DefRestrMult}. The following lemma is then simply borrowed from the usual theory of multiplier algebras.

\begin{Lem}\label{LemSlice}
Let $\omega \in B^*$, the linear dual of $B$, and let $x \in \widetilde{M}(A\otimes^I B)$. Then there exists a unique element $(\id\otimes \omega)(x)\in M(A)$ such that 
\[
(\id\otimes \omega)(x)a = (\id\otimes \omega)(x(a\otimes 1)),\qquad a(\id \otimes \omega)(x) = (\id\otimes \omega)((a\otimes 1)x),\qquad \forall a\in A.
\]
\end{Lem} 
Similarly, we can slice with functionals on the left. 

\subsection{Morphisms between $I$-partial algebras}

Let $A,B$ be non-degenerate, idempotent $I$-partial algebras. 

\begin{Def}
By a \emph{morphism} from $A$ to $B$ we mean an algebra homomorphism 
\[
f: A \rightarrow M(B)
\]
such that for all indices $r$ it holds that
\begin{equation}\label{EqMor2}
\Gr{B}{}{}{}{r} \cdot f(\Gr{A}{}{}{}{r}) = \Gr{B}{}{}{}{r},\qquad   f(\Gr{A}{}{}{r}{})\cdot  \Gr{B}{}{}{r}{} = \Gr{B}{}{}{r}{}.
\end{equation}
\end{Def}
If $f$ is a morphism, it then follows that $f$ is a non-degenerate morphism (in the sense of \cite[Appendix]{VDae94}) between the total algebras, and that its extension to $M(A)$ satisfies 
\begin{equation}\label{EqMor3}
f(\Unit_A(r)) = \Unit_B(r).
\end{equation}
Conversely, any non-degenerate homomorphism $f: A \rightarrow M(B)$ satisfying \eqref{EqMor3} defines a morphism of $I$-partial algebras, since then for example $B\cdot f(A_r) = B_r$ and hence 
\[
B_r\cdot f(A_r) = B\cdot f(A_r^2) = B\cdot f(A_r) = B_r.
\] 
It is then easy to see that morphisms can be composed.

Note that the idempotency condition on our $I$-partial algebras implies that the identity map is a morphism. 

The following lemma is immediate, using the map \eqref{EqQuotMap} and a similar proof as at the end of the proof of Lemma \ref{LemTenProdND}.
\begin{Lem}\label{LemTensProdMorph}
Let $A,B,C,D$ be non-degenerate, idempotent $I^2$-partial algebras, and let $f: A \rightarrow M(B)$ and $g : C\rightarrow M(D)$ be morphisms. Then, viewing $M(B)\otimes M(D) \subseteq M(B\otimes D)$, the map
\[
f\otimes g: A \otimes C \rightarrow M(B\otimes D)
\]
(co-)restricts to a morphism of $I^2$-partial algebras
\[
f\otimes g: A \otimes^I C \rightarrow M(B\otimes^I D)
\]
\end{Lem}

\subsection{Regular partial multiplier bialgebras}
\begin{Def}\label{DefNondegcom}
Let $A$ be a non-degenerate, idempotent $I^2$-partial algebra. A (coassociative) \emph{regular comultiplication} is a morphism $\Delta$ from $A$ to $A\otimes^I A$ such that
\[
\Delta(A) \subseteq \widetilde{M}(A\otimes^I A)
\]
and such that
\begin{align}
(\Delta\otimes \id)\Delta = (\id\otimes \Delta)\Delta.\label{EqCoAssoc}    
\end{align}
\end{Def}
To interpret the coassociativity condition, observe that Lemma \ref{LemTensProdMorph} gives us that 
\[
\Delta \otimes \id: A \otimes^I A \rightarrow M((A\otimes^I A)\otimes^I A)
\]
is non-degenerate as a homomorphism between non-unital algebras, and hence extends to a homomorphism
\[
\Delta \otimes \id: M(A \otimes^I A) \rightarrow M((A\otimes^I A)\otimes^I A).
\]
A similar result holds for $\id\otimes \Delta$, and the coassociativity condition then becomes meaningful by the natural identifications
\[
A\otimes^I(A\otimes^I A) = (A\otimes^I A)\otimes^I A = A\otimes^I A \otimes^I A,
\]
with the last entry defined in the obvious way.

By \eqref{EqMor3}, we then have
\[
\Delta\left(\UnitC{r}{s}\right) = \sum_t \UnitC{r}{t}\otimes \UnitC{t}{s},
\]
and similarly
\[
\Delta(\Unit_s) = \sum_t \Unit_t \otimes \UnitC{t}{s},\qquad \Delta(\Unit^r) = \sum_t \UnitC{r}{t}\otimes \Unit^t,\qquad \Delta(1) = E. 
\]

\begin{Rem}
In the setting of multiplier Hopf algebras, the `regularity' refers to the fact that we are assuming that all spaces 
\[
\Delta(A)(1\otimes A),\qquad (A\otimes 1)\Delta(A),\qquad \Delta(A)(A\otimes 1), \qquad (1\otimes A)\Delta(A)
\]
lie in $A\otimes A$, whereas a non-regular comultiplication would only require this to be the case for the first two spaces. For the theory we are eventually interested in, there is no need to consider the more technical case of non-regular comultiplications. 
\end{Rem}

\begin{Lem}\label{LemProdWher}
Let $a\in \Gr{A}{}{}{}{u}$ and $b \in \Gr{A}{}{}{p}{}$. Then 
\[
\Delta(a)(b\otimes 1) \in A\otimes \Gr{A}{}{p}{}{u}.
\]
\end{Lem} 
We of course have corresponding identities for multiplications at the other places, from left or right.
\begin{proof}
This follows from the fact that 
\[
\Delta(\Unit_u)(\Unit_p\otimes 1) = \sum_{k}\Unit_k\Unit_p\otimes \UnitC{k}{u} = \Unit_p\otimes \UnitC{p}{u}.
\]
\end{proof}

\begin{Def}\label{DefPMB}
A \emph{regular $I$-partial multiplier bialgebra} consists of a triple $(A,\Delta,\varepsilon)$ where $A$ is a non-degenerate, idempotent $I^2$-partial algebra with regular comultiplication $\Delta$, and where 
\[
\varepsilon: A \rightarrow \C
\]
is a linear functional satisfying the following conditions: 
\begin{enumerate}
\item\label{EqCounit1} We have $\varepsilon\left(\Gr{A}{r}{s}{t}{u}\right)=0$ when $r\neq t$ or $s\neq u$,
\item\label{EqCounit2} For each $r\in I$ we have $\varepsilon(\Gr{A}{}{}{}{r}) \neq 0 \neq \varepsilon(\Gr{A}{}{}{r}{})$,
\item\label{EqCounit3} The following homomorphism property holds: for all $s,t\in I$, one has
\[
\varepsilon(ab) = \varepsilon(a)\varepsilon(b)\qquad \textrm{whenever}\quad a\in A_t^s,b\in {}^s_tA,
\]
\item\label{EqCounit4} The following identity holds (using Lemma \ref{LemSlice} to make sense of it):
\begin{equation}\label{EqCounitCond}
(\id\otimes \varepsilon)\Delta(a) = (\varepsilon\otimes \id)\Delta(a) = a,\qquad \forall a \in A.
\end{equation}
\end{enumerate}
\end{Def}
Note that by the third condition, $\varepsilon$ defines a homomorphism
\[
\widetilde{\varepsilon}: A \rightarrow M_{00}(I),\qquad \widetilde{\varepsilon}(a)_{r,s} = \varepsilon(\UnitC{r}{r}a \UnitC{s}{s}), 
\]
where $M_{00}(I)$ is the algebra of $I$-indexed matrices with only finitely many non-zero entries.

We call a map $\varepsilon$ as above a \emph{counit} of $(A,\Delta)$.

The first two conditions in Definition 1.17 might seem artificial a priori. In the next two lemmas, we show that they are in fact equivalent with natural conditions on the pair $(A,\Delta)$. 

\begin{Lem}\label{LemCounitZero}
Let $A$ be a non-degenerate, idempotent $I^2$-partial algebra with regular comultiplication $\Delta$, and let $\varepsilon: A \rightarrow\C$ be a map satisfying \eqref{EqCounit1},  \eqref{EqCounit3} and \eqref{EqCounit4} in Definition \ref{DefPMB}. Then $\varepsilon(\Gr{A}{}{}{}{r}) =0$ or $\varepsilon(\Gr{A}{}{}{r}{}) =0$ if and only if $\Unit_r = 0= \Unit^r$.
\end{Lem}
\begin{proof}
Assume $\varepsilon(\Gr{A}{}{}{}{r})=0$. 

Take $a \in \Gr{A}{}{}{}{r}$ and $b\in \Gr{A}{}{}{r}{}$. By Lemma \ref{LemProdWher} 
\[
\Delta(a)(b\otimes 1) \in A\otimes \Gr{A}{}{r}{}{r}.
\]
Applying $\varepsilon$ to the second leg, we  find 
\[
ab = (\id\otimes \varepsilon)(\Delta(a)(b\otimes 1)) =0.
\]
The non-degeneracy condition on $A$ then implies $\Gr{A}{}{}{}{r} = \Gr{A}{}{}{r}{} =0$, hence $\Unit_r=0$.

Similarly, if $a\in \Gr{A}{}{r}{}{}$ and $b\in \Gr{A}{r}{}{}{}$, then $\Delta(a)(1\otimes b)  \in \Gr{A}{}{r}{}{r}\otimes A$, hence 
\[
ab = (\varepsilon \otimes \id)(\Delta(a)(1\otimes b)) = 0,
\]
so $\Gr{A}{}{r}{}{} = \Gr{A}{r}{}{}{}=0$, hence $\Unit^r=0$. 

Similar arguments apply under the condition $\varepsilon(\Gr{A}{}{}{r}{})=0$. 
\end{proof}

Condition \eqref{EqCounit2} for the counit in Definition \ref{DefPMB} can hence be replaced by the assumption on $A$ that none of the $\Gr{A}{r}{}{}{},\Gr{A}{}{r}{}{},\Gr{A}{}{}{r}{},\Gr{A}{}{}{}{r}$ are zero. 

Note that if $A= (A,\Delta)$ is a regular $I$-partial multiplier bialgebra, then also $A^{\cop} = (A,\Delta^{\opp})$ is a regular $I$-partial multiplier bialgebra for the grading 
\[
\Gr{(A^{\cop})}{r}{s}{t}{u} = \Gr{A}{t}{u}{r}{s}
\]
and the \emph{opposite comultiplication} $\Delta^{\opp}: A^{\cop} \rightarrow \widetilde{M}(A^{\cop} \otimes^I A^{\cop})$ uniquely defined by
\[
\Delta^{\opp}(a)(b\otimes c) = \sum_{i} x_i\otimes y_i\qquad \textrm{if}\qquad \Delta(a)(c\otimes b) = \sum_i y_i\otimes x_i.
\]
If $\varepsilon$ is a counit for $(A,\Delta)$, it remains a counit for $(A,\Delta^{\opp})$. We further have 
\[
\Delta^{\opp}(1) = E^{\opp} = \sum_t \Unit^t\otimes \Unit_t.
\]

Similarly, denote $A^{\opp}$ for the opposite algebra of $A$. We  identify $A$ and $A^{\opp}$ as vector spaces by a map $a \mapsto a^{\opp}$. We endow $A^{\opp}$ with the $I^2$-grading $\Gr{(A^{\opp})}{r}{s}{t}{u} = \left(\Gr{A}{s}{r}{u}{t}\right)^{\opp}$. Note that one can naturally identify $(A\otimes^I A)^{\opp} \cong A^{\opp}\otimes^I A^{\opp}$ and $M(A\otimes^IA)^{\opp} \cong M((A\otimes^I A)^{\opp})$. An inspection of the definition of a regular $I$-partial multiplier bialgebra then reveals that also $(A^{\opp},\Delta)$ is a regular $I$-partial multiplier bialgebra.

Recall now that a comultiplication $\Delta$ is called \emph{full} if $A$ is the smallest vector space $V$ for which $\Delta(A)(1\otimes A) \subseteq V\otimes A$ and the smallest vector space $W$ for which $(A\otimes 1)\Delta(A) \subseteq A\otimes W$, cf.~\cite[Definition 1.4]{VDW13}. Equivalently, this states that any element $a\in A$ can be written as a finite sum 
\begin{equation}\label{EqLeftFull}
a = \sum_i (\id\otimes \omega_i)(\Delta(b_i)(1\otimes c_i)),\qquad b_i,c_i \in A,\omega_i \in \Lin(A,\C),
\end{equation}
\begin{equation}\label{EqRightFull}
a = \sum_i (\omega_i\otimes \id)((b_i\otimes 1)\Delta(c_i),\qquad b_i,c_i \in A,\omega_i \in \Lin(A,\C).
\end{equation}

\begin{Lem}\label{LemFullCopr}
Let $A$ be a non-degenerate idempotent $I^2$-partial algebra with regular comultiplication $\Delta$, and let $\varepsilon: A \rightarrow\C$ be a map satisfying \eqref{EqCounit2},  \eqref{EqCounit3} and \eqref{EqCounit4} in Definition \ref{DefPMB}. Then the following are equivalent: 
\begin{itemize}
    \item The map $\varepsilon$ satisfies \eqref{EqCounit1} in Definition \ref{DefPMB} (i.e.\ $\varepsilon$ is a counit map).
    \item The comultiplication $\Delta$ is full.
    \item The comultiplication $\Delta^{\opp}$ is full. 
\end{itemize}
\end{Lem} 
\begin{proof}
Assume first that $\varepsilon$ is a counit map. By the partial multiplicativity of $\varepsilon$ and the counit property, we have 
\[
(\id\otimes \varepsilon)(\Delta(\Gr{A}{}{r}{}{s})(1\otimes \Gr{A}{s}{}{s}{})) = \varepsilon(\Gr{A}{s}{}{s}{}) (\id\otimes \varepsilon)(\Delta(\Gr{A}{}{r}{}{s})(1\otimes \Unit^s))=  (\id\otimes \varepsilon)(\Delta(\Gr{A}{}{r}{}{s})) = \Gr{A}{}{r}{}{s},
\]
where in the last but one step we used properties \eqref{EqCounit1} and \eqref{EqCounit2} of Definition \ref{DefPMB}. A second analogous computation for the right leg shows that $\Delta$ is full. As $\varepsilon$ is also a counit for $(A,\Delta^{\opp})$ it follows that also $\Delta^{\opp}$ is full. 

Assume now conversely that $\Delta$ is full. Then necessarily $\Gr{A}{}{k}{}{l}$ is the smallest subspace $V$ such that 
\[
\sum_m \Delta(\Gr{A}{}{k}{}{m})(1\otimes \Gr{A}{l}{}{m}{}) \subseteq V \otimes A. 
\]
But if $k \neq l$, then applying the counit to the first leg of the left hand expression gives $0$ by \eqref{EqCounit4}. Hence $\varepsilon$ vanishes on $\Gr{A}{}{k}{}{l}$. Similarly, the other condition for the fullness of $\Delta$ gives that $\varepsilon$ vanishes on $\Gr{A}{k}{}{l}{}$ if $k\neq l$. Hence \eqref{EqCounit1} is satisfied.

Again, by symmetry we can also conclude \eqref{EqCounit1} from the fullness of $\Delta^{\opp}$. 
\end{proof}

From the fullness of $\Delta$ and condition \eqref{EqCounit4} in Definition \ref{DefPMB}, we immediately obtain:

\begin{Cor}\label{CorUniqueCounit}
Let $A$ be a non-degenerate, idempotent $I^2$-partial algebra, and let $\Delta: A \rightarrow \widetilde{M}(A\otimes^I A)$ be a regular comultiplication. Then there exists at most one counit map.
\end{Cor} 

We hence call $\varepsilon$ \emph{the} counit map if it exists. In the following, we will simply refer to $(A,\Delta)$ as a regular $I$-partial multiplier bialgebra, and write $\varepsilon$ by default for the counit.

\begin{Def}
We define 
\[
\msB_A^{\mbs}= \Span \{\Unit_s\mid s\in I\} \subseteq M(A),\qquad \msB_A^{\mbt}= \Span\{\Unit^s\mid s\in I\} \subseteq M(A),
\]
and call them respectively the \emph{source} and \emph{target} algebras of $A$. We refer to $\msB_A^{\mbs}$ and $\msB_A^{\mbt}$ together as the \emph{base algebras}. 
\end{Def}

By Lemma \ref{LemCounitZero}, both $\msB_A^{\mbs}$ and $\msB_A^{\mbt}$ are faithful copies of the algebra
\[
\msB_A = \Fun_{\fin}(I) = \{f: I \rightarrow \C\mid f\textrm{ has finite support}\},
\]
which we will at times also refer to simply as \emph{the base algebra}. We then write 
\[
\mbs: \msB_A \rightarrow \msB_A^{\mbs},\quad f\mapsto \sum_s f(s)\Unit_s,\qquad \mbt: \msB_A \rightarrow \msB_A^{\mbt},\quad f\mapsto \sum_sf(s)\Unit^s. 
\]

To end, we point out the connection with the theory of weak multiplier bialgebras. To do this, we first recall the notion of a weak multiplier bialgebra, as defined in \cite[Definition 2.1 and Definition 2.3]{BG-TL-C15}.

\begin{Def}\label{DefWeakMultBiAlg}
    Let $A$ be an idempotent, non-degenerate algebra with multiplication $\mu: A\otimes A\to A$. Let $\Delta: A \to M(A\otimes A)$ be a multiplicative linear map (called a \emph{comultiplication} on $A$). 
    The pair $(A,\Delta)$ is called a \emph{regular weak multiplier bialgebra} if there exist both an idempotent element $E$ in $M(A\otimes A)$ and a linear map $\varepsilon: A\to \C$ (called a \emph{counit}), such that the following six conditions are satisfied:
    \begin{enumerate}[label=(\roman*)]
        \item The right hand sides of the expressions below,
            \begin{align*}
                T_1(a\otimes b) = \Delta(a)(1\otimes b),\\
                T_2(a\otimes b) = (a\otimes 1)\Delta(b),\\
                T_3(a\otimes b) = (1\otimes b)\Delta(a),\\
                T_4(a\otimes b) = \Delta(b)(a\otimes 1),
            \end{align*}
            have image in the two-sided ideal $A\otimes A$ inside $M(A\otimes A)$, and hence define maps $T_i: A \otimes A \rightarrow A \otimes A$,
        \item The maps $T_1$ and $T_2$ satisfy the following commutation relation:
            \[(T_2\otimes \id)(\id \otimes T_1) = (\id\otimes T_1)(T_2\otimes\id).\]
            or, equivalently, the maps $T_3$ and $T_4: A\otimes A\to A\otimes A$ satisfy:
            \[(T_4\otimes \id)(\id \otimes T_3) = (\id\otimes T_3)(T_4\otimes \id).\]
            These equivalent conditions are called \emph{coassociativity}.
        \item The counit $\varepsilon$ satisfies
            \[(\varepsilon\otimes \id)T_1 = \mu = (\id\otimes \varepsilon)T_2,\]
            or, equivalently,
            \[(\varepsilon\otimes \id)T_3 = \mu^\op = (\id\otimes \varepsilon)T_4.\]
        \item We have the following equalities
        \begin{align*}
        &\Span\{\Delta(a)(b\otimes c)\mid a,b,c\in A\} = \Span\{E(b\otimes c)\mid b,c\in A \},\\
        &\Span\{(a\otimes b)\Delta(c)\mid a,b,c\in A\} = \Span\{(a\otimes b)E\mid a,b\in A \}
        \end{align*}
        of subspaces of $A\otimes A$.
        \item The elements
        \begin{align}\label{EqCondEBI}
            (E\otimes 1)(1\otimes E) = (1\otimes E)(E\otimes 1) = (\Delta\otimes \id)(E) = (\id\otimes \Delta)(E)
        \end{align}
        of $M(A\otimes A\otimes A)$ are all equal
        \item For any $a,b,c\in A$, we have these two equalities
        \begin{align*}
            &(\varepsilon\otimes\id)((1\otimes a)E(b\otimes c))
            = (\varepsilon \otimes \id)(\Delta(a)(b\otimes c))
            \textrm{ and }\\
            &(\varepsilon\otimes\id)((a\otimes b)E(1\otimes c)) = (\varepsilon\otimes\id)((a\otimes b)\Delta(c)).
        \end{align*} 
    \end{enumerate}
\end{Def}

\begin{Prop}\label{PropWeakMultBiAlg}
Let $(A,\Delta)$ be a regular $I$-partial multiplier bialgebra. Then $(A,\Delta)$ is a regular weak multiplier bialgebra as defined in Definition \ref{DefWeakMultBiAlg}. 
\end{Prop}
\begin{proof}
We have a natural candidate for the idempotent multiplier that Definition \ref{DefWeakMultBiAlg} requires, namely $E = \sum_s \Unit_s \otimes \Unit^s\in M(A\otimes A)$. Similarly, we set the linear map $\varepsilon$ to be the `partial' counit. We then check whether the six conditions in Definition \ref{DefWeakMultBiAlg} hold, for these choices of $E$ and $\varepsilon$.

By the definition of a `partial' comultiplication, $\Delta$ has image in the \emph{restricted} multiplier algebra $\wM(A\otimes^I A)$. This implies that the maps $T_1, \dots , T_4$ all have image inside $A\otimes^I A \subseteq A\otimes A$. Thus, condition $(i)$ is seen to be satisfied. The condition $(ii)$ follows immediately from the coassociativity condition \eqref{EqCoAssoc}, and condition $(iii)$ is immediate from Definition \ref{DefPMB}\eqref{EqCounit4}.

For condition $(iv)$, we will only check the first equality of subspaces, the second check is very similar.
Expressions of the form $\Delta(a)(b\otimes c)$, which span the subspace on the left, can be written as follows:
\[\Delta(a)(b\otimes c) = \Delta(1\cdot a)(b\otimes c) = \Delta(1)\Delta(a)(1\otimes c)(b\otimes 1) = ET_1(a\otimes c)(b\otimes 1).\]
These are elements of the right hand subspace, since the middle factor is, by $(i)$, an element of $A\otimes A$.
On the other hand, since $\Delta$ is a morphism, it is in particular non-degenerate as a homomorphism $A \rightarrow M(A\otimes^I A)$, hence by \eqref{EqIdempotentnondeg} we get
\[
\Delta(A)(A\otimes A) = \Delta(A)E(A\otimes A)= \Delta(A)(A\otimes^IA)(A\otimes A) = (A\otimes^IA)(A\otimes A) =  E(A\otimes A).
\]

The final condition $(vi)$ requires to check that for all $a,b,c\in A$ we have
\[
(\varepsilon \otimes \id)((1\otimes a)E(b\otimes c)) = (\varepsilon\otimes \id)(\Delta(a)(b\otimes c)),
\]
\[
(\varepsilon \otimes \id)((a\otimes b)E(1\otimes c)) = (\varepsilon \otimes \id)((a\otimes b)\Delta(c)). 
\]
Let us verify the first identity. We have, referring to the numbering in Definition \ref{DefPMB},  
\begin{eqnarray*}
(\varepsilon\otimes \id)(\Delta(a)(b\otimes c)) &=& \sum_{r,s,t} (\varepsilon\otimes \id)(\Delta(a\UnitC{r}{t})(\UnitC{r}{s}b\otimes \UnitC{s}{t}c))\\
&\overset{\eqref{EqCounit3}}{=}& \sum_{r,s,t} \varepsilon(\UnitC{r}{s}b)  (\varepsilon\otimes \id)(\Delta(a\UnitC{r}{t})(\UnitC{r}{s}\otimes \UnitC{s}{t}c))\\
&=& \sum_{r,s,t} \varepsilon(\UnitC{r}{s}b)  (\varepsilon\otimes \id)(\Delta(a\UnitC{r}{t})(1\otimes \UnitC{s}{t}c))\\
&\overset{\eqref{EqCounit4}}{=}& \sum_{r,s,t} \varepsilon(\UnitC{r}{s}b)  a\UnitC{r}{t}\UnitC{s}{t}c \\
&=& \sum_{r} \varepsilon(\UnitC{r}{r}b)  a\Unit^rc\\
&\overset{\eqref{EqCounit1}}{=}& \sum_{r,t} \varepsilon(\Unit_rb)  a\Unit^rc\\
&=& (\varepsilon \otimes \id)((1\otimes a)E(b\otimes c)).
\end{eqnarray*}
The other identity follows similarly.
\end{proof}

We will not consider here the converse question of when a (regular) weak multiplier bialgebra (with full comultiplication) is a regular $I$-partial multiplier bialgebra, but see Proposition \ref{PropWeakToPartial}.

\subsection{Partial multiplier Hopf algebras}

\begin{Def}\label{DefCanMaps}
A regular $I$-partial multiplier bialgebra $(A,\Delta)$ is called a \emph{regular $I$-partial multiplier Hopf algebra} if and only if the following four maps are bijective:

\begin{align}
\label{canr}&\can_r: \underset{r}{\oplus} (\Gr{A}{}{}{}{r} \otimes \Gr{A}{}{}{r}{})  \rightarrow \underset{r}{\oplus}(\Gr{A}{}{}{r}{} \otimes \Gr{A}{r}{}{}{}), \qquad a\otimes b \mapsto \Delta(a)(1\otimes b),\\
\label{canl}&\can_l: \underset{r}{\oplus}(\Gr{A}{}{r}{}{} \otimes \Gr{A}{r}{}{}{})  \rightarrow \underset{r}{\oplus} (\Gr{A}{}{}{}{r} \otimes \Gr{A}{}{r}{}{}),\qquad a\otimes b \mapsto (a\otimes 1)\Delta(b),\\
\label{canrc}&\can_r^c: \underset{r}{\oplus} (\Gr{A}{}{r}{}{} \otimes \Gr{A}{r}{}{}{})  \rightarrow \underset{r}{\oplus}(\Gr{A}{r}{}{}{} \otimes \Gr{A}{}{}{r}{}), \qquad a\otimes b \mapsto \Delta^{\opp}(a)(1\otimes b),\\
\label{canlc}&\can_l^c: \underset{r}{\oplus}(\Gr{A}{}{}{}{r} \otimes \Gr{A}{}{}{r}{})  \rightarrow \underset{r}{\oplus} (\Gr{A}{}{r}{}{} \otimes \Gr{A}{}{}{}{r}),\qquad a\otimes b \mapsto (a\otimes 1)\Delta^{\opp}(b).
\end{align}
\end{Def}

As is easily seen, $(A,\Delta)$ is a regular $I$-partial multiplier Hopf algebra if and only if $(A,\Delta^{\cop})$ or $(A^{\opp},\Delta)$ are. 

In the following, we will also use the notation
\begin{equation}\label{canrgen} 
\widetilde{\can}_r: A\otimes A\rightarrow A\otimes A, \qquad a\otimes b \mapsto \Delta(a)(1\otimes b)
\end{equation}
for the natural extension of $\can_r$ to $A\otimes A$, and similarly for the other maps. Note that these maps already appeared in Definition \ref{DefWeakMultBiAlg}. In general, we will use the notation $T_i$ to refer to general weak multiplier bialgebras, and use the notation of Definition \ref{DefCanMaps} for the specific case of $I$-partial multiplier bialgebras.

We relate now the notion of a regular $I$-partial multiplier Hopf algebra to that of a regular weak multiplier Hopf algebra as introduced in \cite[Definition 1.14]{VDW13}. By \cite[Theorem 2.9]{BG-TL-C15}, any such regular weak multiplier Hopf algebra will be a regular weak multiplier bialgebra, so we can skip any condition in \cite[Definition 1.14]{VDW13} referring only to the underlying regular weak multiplier bialgebra structure. 

\begin{Def}\label{DefCorrWHOPF}
    Let $(A,\Delta)$ be a regular weak multiplier bialgebra such that $\Delta$ and $\Delta^{\opp}$ are full (cf.\ \eqref{EqLeftFull} and \eqref{EqRightFull}). Then $(A,\Delta)$ is called a \emph{regular weak multiplier Hopf algebra} if the following conditions are satisfied:
\begin{enumerate}
\item The maps $T_1,T_4$ in Definition \ref{DefWeakMultBiAlg} have image $E(A\otimes A)$, while $T_2,T_3$ have image $(A\otimes A)E$.
\item The kernels of the canonical maps are of the form
        \[ 
\Ker(T_i) = (\id-G_i)(A\otimes A)
\]
        where the linear maps $G_1, G_2,G_3,G_4: A\otimes A\to A\otimes A$ are characterized by the equalities
        \begin{align}
            (G_1 \otimes\id)(\Delta_{13}(a)(1\otimes b\otimes c)) = \Delta_{13}(a)(1\otimes E)(1\otimes b\otimes c),\label{EqG1}\\
            (\id\otimes G_2)((a\otimes b\otimes 1)\Delta_{13}(c)) = (a\otimes b\otimes 1)(E\otimes 1)\Delta_{13}(c),\label{EqG2}\\
          (G_3\otimes\id)((1\otimes b\otimes c)\Delta_{13}(a)) = (1\otimes b\otimes c)(1\otimes E)\Delta_{13}(a),\label{EqG3}\\
            (\id\otimes G_4)(\Delta_{13}(c)(a\otimes b\otimes 1)) = \Delta_{13}(c)(E\otimes 1)(a\otimes b\otimes 1),\label{EqG4}
        \end{align}
        for all $a,b,c\in A$.
\end{enumerate}

\end{Def}

%The next proposition shows that a regular $I$-partial multiplier Hopf algebra $(A,\Delta)$ is a regular weak multiplier Hopf algebra in the sense of \cite{VDW13}. 

\begin{Prop}\label{PropPartialToWeak}
A regular $I$-partial multiplier Hopf algebra $(A,\Delta)$ is a regular weak multiplier Hopf algebra.
\end{Prop} 
\begin{proof}
We know from Lemma \ref{LemFullCopr} that $\Delta$ is full, and it is also easily verified directly that 
\[
E = \sum_p \Unit_p \otimes \Unit^p
\]
satisfies $T_1(A\otimes A) = E(A\otimes A)$ and $T_2(A\otimes A) = (A\otimes A)E$ by bijectivity of $\can_r$ and $\can_l$. By considering the regular $I$-partial multiplier Hopf algebra $(A^{\opp},\Delta)$, the conditions on the images of $T_3,T_4$ follow as well. 

To see that also condition (2) is satisfied above, let us first show that the map $G_1$ is given by the idempotent
\[
\widetilde{G}_1: A \otimes A \rightarrow A\otimes A,\qquad a\otimes b \mapsto \sum_p a\Unit_r \otimes \Unit_rb.  
\]
Indeed, for $a\in A,b\in \Gr{A}{}{}{r}{}$ and $c\in \Gr{A}{s}{}{}{}$ we have 
\[
\Delta_{13}(a)(1\otimes b\otimes c) \in \Gr{A}{}{}{}{s}\otimes \Gr{A}{}{}{r}{}\otimes A,
\]
so
\[
(\widetilde{G}_1\otimes \id)(\Delta_{13}(a)(1\otimes b\otimes c)) = \delta_{rs} \Delta_{13}(a)(1\otimes b\otimes c) = \Delta_{13}(a)(1\otimes E(b\otimes c)),
\]
proving $\widetilde{G}_1 = G_1$. Hence by assumption $\Ker(\widetilde{\can}_r) = \Ker(T_1)$ is the image of $\id- G_1$. 

Similar computations can be done for $G_2,G_3$ and $G_4$.
\end{proof}

The previous proposition allows us to fully use the results of \cite{VDW13}. In particular, \cite[Proposition 4.9]{VDW13} tells us the following (see also the beginning of the proof of Proposition \ref{PropWeakToPartial}):
\begin{Cor}\label{CorLocUn}
A regular $I$-partial multiplier Hopf algebra $A$ admits \emph{local units}: for each finite subset $\{a_1,\ldots,a_n\}\subseteq A$ there exists $e\in A$ with $a_i e = ea_i = a_i$ for all $i$. 
\end{Cor}

The following gives an alternative characterisation of regular $I$-partial multiplier Hopf algebras. We will now employ the Sweedler notation 
\[
\Delta(a) = a_{(1)}\otimes a_{(2)} \in \widetilde{M}(A\otimes^I A),
\]
making sure that all expressions using this notation are sufficiently `covered', see \cite{VDae08}.

\begin{Prop}\label{PropExAntip}
A regular $I$-partial multiplier bialgebra $(A,\Delta)$ is a regular $I$-partial multiplier Hopf algebra if and only if there exists a bijective anti-homomorphism $S: A \rightarrow A$ such that
\begin{equation}\label{EqImUnit}
S(\Unit_r) = \Unit^r,\qquad S(\Unit^r) = \Unit_r
\end{equation}
and
\begin{equation}\label{EqLeftInv}
S(a_{(1)})a_{(2)}b = \varepsilon(a)\Unit_sb,\qquad  \forall a\in \Gr{A}{}{}{}{s}, b\in A. 
\end{equation}
\begin{equation}\label{EqRightInv}
ba_{(1)}S(a_{(2)}) = \varepsilon(a)b\Unit^r,\qquad \forall a\in \Gr{A}{r}{}{}{},b\in A.
\end{equation}
Moreover, such a map $S$ is unique and satisfies
\begin{equation}\label{EqAntipodComult}
S(a)_{(1)}\otimes S(a)_{(2)} = S(a_{(2)})\otimes S(a_{(1)}),\qquad a\in A, \end{equation}
\begin{equation}\label{EqAntipodCounit}
\varepsilon \circ S = \varepsilon. 
\end{equation}
\end{Prop}
Note that there is no problem to interpret \eqref{EqAntipodComult} as $S$ is assumed to be an anti-isomorphism. We call $S$ the \emph{antipode} of the regular $I$-partial multiplier Hopf algebra.
\begin{proof}
Assume that $(A,\Delta)$ is a regular $I$-partial multiplier bialgebra with $S: A \rightarrow A$ a bijective anti-homomorphism satisfying  \eqref{EqImUnit}, \eqref{EqLeftInv} and \eqref{EqRightInv}. Consider the maps 
\begin{equation}\label{Eqcanrinvgen} 
\widetilde{\can}_r^-: A\otimes A \rightarrow  A\otimes A, \qquad a\otimes b \mapsto a_{(1)}\otimes S(a_{(2)})b,
\end{equation}
\begin{equation}\label{Eqcanlinvgen}
\widetilde{\can}_l^-: A\otimes A \rightarrow  A\otimes A ,\qquad a\otimes b \mapsto a S(b_{(1)})\otimes b_{(2)},
\end{equation}
\begin{equation}\label{Eqcanrinvopgen} 
\widetilde{\can}_r^{c,-}:   A\otimes A \rightarrow  A\otimes A , \qquad a\otimes b \mapsto S^{-1}(b_{(1)})a\otimes b_{(2)},
\end{equation}
\begin{equation}\label{Eqcanlinvopgen}
\widetilde{\can}_l^{c,-}:  A\otimes A \rightarrow  A\otimes A,\qquad a\otimes b \mapsto a_{(1)} \otimes bS^{-1}(a_{(2)}).
\end{equation}
As we are assuming $S$ to be bijective and antimultiplicative, it is easily seen that these maps indeed have ranges in $A \otimes A$, and \eqref{EqImUnit} easily implies that they restrict/corestrict to maps 
\begin{equation}\label{Eqcaninv} 
\can_r^-:  \underset{r}{\oplus}(\Gr{A}{}{}{r}{} \otimes \Gr{A}{r}{}{}{}) \rightarrow \underset{r}{\oplus} (\Gr{A}{}{}{}{r} \otimes \Gr{A}{}{}{r}{}),\qquad \can_l^-: \underset{r}{\oplus} (\Gr{A}{}{}{}{r} \otimes \Gr{A}{}{r}{}{}) \rightarrow \underset{r}{\oplus}(\Gr{A}{}{r}{}{} \otimes \Gr{A}{r}{}{}{}),
\end{equation}
\begin{equation}\label{Eqcaninvop} 
\can_r^{c,-}:  \underset{r}{\oplus}(\Gr{A}{r}{}{}{} \otimes \Gr{A}{}{}{r}{}) \rightarrow  \underset{r}{\oplus} (\Gr{A}{}{r}{}{} \otimes \Gr{A}{r}{}{}{}),\qquad 
\can_l^{c,-}:   \underset{r}{\oplus} (\Gr{A}{}{r}{}{} \otimes \Gr{A}{}{}{}{r}) \rightarrow \underset{r}{\oplus}(\Gr{A}{}{}{}{r} \otimes \Gr{A}{}{}{r}{}),
\end{equation}
Further, we have for $a,b\in A$ that
\[
(\widetilde{\can}_r \circ \widetilde{\can}_r^-)(a\otimes b) = a_{(1)}\otimes a_{(2)}S(a_{(3)})b = \sum_r a_{(1)}\otimes \varepsilon(\Unit^ra_{(2)})\Unit^rb = \sum_r \Unit_r a_{(1)}\otimes \varepsilon(a_{(2)})\Unit^r b = E(a\otimes b),
\]
and a similar computation for $\widetilde{\can}_r^-\circ \widetilde{\can}_r$ shows that $\can_r$ is the inverse of $\can_r^-$. These same computations can then be repeated for the other maps, showing that $(A,\Delta)$ is a regular $I$-partial multiplier Hopf algebra. By \cite[Proposition 2.4 and Proposition 2.7]{VDW13}, its (unique) antipode must coincide with the map $S$, which is hence  an anti-cohomomorphism by \cite[Proposition 4.4]{VDW13}. Then \eqref{EqAntipodCounit} follows from the uniqueness of the counit $\varepsilon$. 

Assume now conversely that $(A,\Delta)$ is a  regular $I$-partial multiplier Hopf algebra, and hence a regular weak multiplier Hopf algebra. By \cite[Proposition 2.4, Proposition 3.5 and Proposition 4.3]{VDW13} there exists a unique bijective anti-homomorphism $S: A \rightarrow A$ such that the maps \eqref{Eqcanrinvgen}, \eqref{Eqcanlinvgen}, \eqref{Eqcanrinvopgen}, \eqref{Eqcanlinvopgen} satisfy
\begin{equation}\label{EqRange1}
(\widetilde{\can}_r \circ \widetilde{\can}_r^-)(a\otimes b) = E(a\otimes b),\qquad (\widetilde{\can}_l\circ \widetilde{\can}_l^-)(a\otimes b) = (a\otimes b)E
\end{equation}
\begin{equation}\label{EqRange2}
(\widetilde{\can}_r^{c} \circ \widetilde{\can}_r^{c,-})(a\otimes b) = E(a\otimes b),\qquad (\widetilde{\can}_l^{c}\circ \widetilde{\can}_l^{c,-})(a\otimes b) = (a\otimes b)E,
\end{equation}
Note that ranges of the maps $\widetilde{\can}_r^-,\ldots$ can automatically be interpreted as elements in $A\otimes A$ as we can already use that $S$ is anti-multiplicative and bijective (in the more general setting of \cite[Proposition 2.4]{VDW13} some more care is needed as these properties of $S$ can not be invoked). By \cite[Proposition 2.6]{VDW13}, the map $S$ then moreover satisfies 
\begin{equation}\label{EqWeakInv}
a_{(1)}S(a_{(2)})a_{(3)} = a,\qquad a_{(3)}S^{-1}(a_{(2)})a_{(1)} = a,\qquad \forall a\in A.
\end{equation}

Now \cite[Lemma 3.2]{VDW13} tells us that for each $a\in A$, the element $S(a_{(1)})a_{(2)} \in M(A)$ is a linear combination of elements of the form $(\id\otimes \omega(b \emdash c))E$ for $b,c\in A$, so we obtain in particular a linear map 
\[
A \rightarrow \msB_A^{\mbs},\qquad a \mapsto S(a_{(1)})a_{(2)}.
\]
In particular, there exists a linear map $F:A  \rightarrow \msB_A$ such that 
\[
S(a_{(1)})a_{(2)} = \sum_s F(a)(s) \Unit_s. 
\]
Put $f: A \rightarrow \C$ the linear map with $f(a) = \sum_{s\in I} F(a)(s)$.  Since $F(a\Unit_s) = F(a) \Unit_s$ by the definition of $F$, we see that $F(a)(s) = f(a\Unit_s)$. 

Since now $\widetilde{\can}_l \circ \widetilde{\can}_l^-$ is the same as right multiplication with $E$, we find that 
\[
\sum_s a\Unit_s \otimes b\Unit^s= aS(b_{(1)})b_{(2)}\otimes b_{(3)}  = \sum_s f(b_{(1)}\Unit_s)a\Unit_s \otimes b_{(2)} = \sum_s f(b_{(1)})a\Unit_s \otimes b_{(2)}\Unit^s.
\]
It follows that $(f\otimes \id)\Delta(b)= b$ for all $b\in A$, and the fullness of $\Delta$ implies that $f =\varepsilon$. Similarly one can show \eqref{EqRightInv}. 

Note now that \eqref{EqLeftInv} and \eqref{EqRightInv} already imply \eqref{EqImUnit}: given $r,s\in I$, choose $a\in A_s$ with $\varepsilon(a)=1$. Then for $b\in \Gr{A}{}{}{s}{}$ we have 
\begin{multline*}
S(\Unit^r)b = S(\Unit^r)\varepsilon(a)b = S(\Unit^r)S(a_{(1)})a_{(2)}b = S(a_{(1)}\Unit^r)a_{(2)}b = S((a\Unit^r)_{(1)})(a\Unit^r)_{(2)}b  \\
= \varepsilon(a\Unit^r)b = \delta_{r,s} \varepsilon(a)b = \delta_{r,s}b,
\end{multline*}
where in the beginning of the second line we use that $\varepsilon$ has support on $\underset{t,u}{\oplus}\Gr{A}{t}{u}{t}{u}$. Hence $S(\Unit^r) = \Unit_r$, and similarly one shows $S(\Unit_r) =\Unit^r$. 
\end{proof}

\begin{Rem}\label{RemSourceTarget}
Recall from \cite[Proposition 2.3]{VDW20} that the maps 
\begin{equation}\label{EqSourceTarget}
\varepsilon_{\mbs}: A \rightarrow M(A),\quad a\mapsto S(a_{(1)})a_{(2)},\qquad \varepsilon_{\mbt}: A \rightarrow M(A),\quad a \mapsto a_{(1)}S(a_{(2)}),
\end{equation}
are called respectively the \emph{source} and \emph{target map}. From Proposition \ref{PropExAntip}, we see that they are given by
\begin{equation}\label{EqAntDef1}
\varepsilon_{\mbs}(a)  = S(a_{(1)})a_{(2)} = \sum_s \varepsilon(a\Unit^s)\Unit_s,\qquad \varepsilon_{\mbt}(a) = a_{(1)}S(a_{(2)})= \sum_t \varepsilon(\Unit_ta)\Unit^t. 
\end{equation}
It is then clear that the images of these maps coincide respectively with the source algebra $\msB_A^{\mbs}$ and target algebra $\msB_A^{\mbt}$ defined before. Note that this terminology is slightly different from the one in \cite[Definition 2.1]{VDW20}: the source and target algebra defined there are rather given by the \emph{multiplier} algebras of $\msB_A^{\mbs}$ and $\msB_A^{\mbt}$, which are copies of the algebra $\Fun(I)$ of all complex valued functions on $I$. 
\end{Rem}

Note that \eqref{EqAntDef1} can also be rewritten as 
\begin{equation}\label{EqAntDef2}
a_{(2)}S^{-1}(a_{(1)}) = \sum_t \varepsilon(\Unit^ta)\Unit_t,\qquad S^{-1}(a_{(2)})a_{(1)} = \sum_s \varepsilon(a\Unit_s)\Unit^s.
\end{equation}
We then see that the antipode of $(A,\Delta^{\opp})$ is given by $S^{-1}$ (see also the discussion following \cite[Proposition 4.3]{VDW13}). Similarly, $(A^{\opp},\Delta)$ has the antipode $a^{\opp}\mapsto S^{-1}(a)^{\opp}$. 

We now state a converse to Proposition \ref{PropPartialToWeak}.

\begin{Prop}\label{PropWeakToPartial}
Let $A$ be an $I^2$-partial algebra, and assume the total algebra $A$ is non-degenerate and idempotent, with $\Unit_k \neq 0 \neq \Unit^k$ for all $k\in I$. Assume further that $(A,\Delta)$ defines a regular weak multiplier Hopf algebra, with antipode $S$, such that 
\begin{equation}\label{EqIm1}
\Delta(1) = \sum_s \Unit_s \otimes \Unit^s.
\end{equation}
Then $A$ is non-degenerate and idempotent as an $I^2$-algebra, and $(A,\Delta)$ defines a regular $I$-partial multiplier Hopf algebra.
\end{Prop}
\begin{proof}
We first prove that $A$ is necessarily a non-degenerate and idempotent $I^2$-algebra. In fact, we will show that 
\begin{equation}\label{EqNondegId}
a\in \Gr{A}{}{r}{}{s}\qquad \Rightarrow \qquad a\in a\Gr{A}{}{r}{}{s},
\end{equation}
which implies at once the non-degeneracy and idempotency condition on $\Gr{A}{}{r}{}{s}$. The non-degeneracy and idempotency condition on $\Gr{A}{r}{}{s}{}$ is then proven similarly.  

To prove \eqref{EqNondegId}, we can simply adapt the proof of \cite[Proposition 3.9]{VDW13}. Indeed, it is sufficient to have that $\omega(a)=0$ for all functionals $\omega: A\rightarrow \C$ with $\omega(a\Gr{A}{}{r}{}{s})=0$. Pick such an $\omega$. Then for all $b\in A,c\in \Gr{A}{}{r}{}{s}$ we have that 
\[
\omega(aS(b_{(1)})c) b_{(2)} =0.
\]
Hence $\omega(aS(b_{(1)})-\UnitC{r}{s})\otimes b_{(2)} = 0$, and applying the evaluation map we obtain that $\omega(aS(b_{(1)})b_{(2)}\UnitC{r}{s})=0$. As $\Unit^s \neq 0$ by assumption, non-degeneracy of $A$ implies $\Gr{A}{}{s}{}{}\cdot\Gr{A}{s}{}{}{} \neq 0$, hence by \cite[Lemma 3.2]{VDW13} we obtain $\omega(a\Unit_s\UnitC{r}{s}) = \omega(a)=0$. 

Having proven that $A$ is a non-degenerate and idempotent $I^2$-algebra, we now show that $\Delta$ is a comultiplication for $A$. Since $\Delta(A) = \Delta(1)\Delta(A)\Delta(1)$, it follows immediately from Definition \ref{DefCorrWHOPF}, \eqref{EqIm1} and Lemma \ref{LemEqQuotMap} that $\Delta: A \rightarrow M(A\otimes A)$ corestricts to a homomorphism
\[
\Delta: A \rightarrow \widetilde{M}(A\otimes^I A),
\]
and that 
\begin{equation}\label{EqRangeGal}
\Delta(A)(1\otimes A)  = \Delta(A)(A\otimes 1)= \Delta(1)(A\otimes A),\qquad (A\otimes 1)\Delta(A) = (1\otimes A)\Delta(A) = (A\otimes A)\Delta(1).
\end{equation}
It then follows from \eqref{EqCondEBI} that
\[
\Delta(\Unit^r) = (\Unit^r\otimes 1)\Delta(1),\qquad \Delta(\Unit_r) = (1\otimes \Unit_r)\Delta(1),
\]
so that $\Delta$ is indeed a morphism from $A$ to $A\otimes^I A$, with coassociativity simply inherited. 

Let $\varepsilon: A \rightarrow \C$ be the counit of $(A,\Delta)$ as a regular weak multiplier Hopf algebra. Then \eqref{EqCounit4} in Definition \ref{DefPMB} holds by definition. On the other hand, we also have from \eqref{EqRangeGal} that
\[
(A\otimes 1)\Delta(\Gr{A}{}{r}{}{s}) = (A\otimes A)\Delta(\UnitC{r}{s}).  
\]
Applying $\varepsilon$ to the second leg and using idempotency, we find
\[
\Gr{A}{}{r}{}{s} = \sum_k \varepsilon(\Gr{A}{}{k}{}{s}) \Gr{A}{}{r}{}{k}. 
\]
Hence if $k\neq s$, we must have $\varepsilon(\Gr{A}{}{k}{}{s})=0$ since otherwise $\UnitC{r}{k}= 0$ for all $r$, in contradiction with $\Unit_k\neq 0$. Similarly, we find $\varepsilon(\Gr{A}{k}{}{s}{})=0$ for all $k\neq s$, so  \eqref{EqCounit1} in Definition \ref{DefPMB} is satisfied. The multiplicativity condition in \eqref{EqCounit3} is then verified as follows: let $a\in \Gr{A}{}{r}{}{s}$ and $b\in \Gr{A}{r}{}{s}{}$, and choose $c\in \Gr{A}{}{}{}{r}$ and $d\in \Gr{A}{}{}{r}{}$ with $cd \neq 0$.  Write 
\[
c\otimes a = \sum_i (p_i\otimes 1)\Delta(q_i),\qquad d\otimes b = \sum_j \Delta(z_j)(w_j\otimes 1).
\]
Then
\begin{multline*}
\varepsilon(ab)cd = \sum_{i,j} (\id\otimes \varepsilon)((p_i\otimes 1)\Delta(q_iz_j)(w_j\otimes 1)) = \sum_{i,j} p_iq_iz_jw_j \\ =  \left(\sum_{i} (\id\otimes \varepsilon)((p_i\otimes 1)\Delta(q_i))\right)\left(\sum_j (\id\otimes \varepsilon)(\Delta(z_j)(w_j\otimes 1))\right) = \varepsilon(a)\varepsilon(b)cd,
\end{multline*}
hence $\varepsilon(ab) = \varepsilon(a)\varepsilon(b)$. 

We have shown up to now that $(A,\Delta)$ is an $I$-partial multiplier bialgebra. It now follows from the computation of the maps $G_1,G_2,\ldots$ in Proposition \ref{PropPartialToWeak} that the requirements in Definition \ref{DefCorrWHOPF} imply the isomorphisms in Definition \ref{DefCanMaps}, hence $(A,\Delta)$ is an $I$-partial multiplier Hopf algebra.
\end{proof}

To end, we comment on the connection with the theory of \emph{regular multiplier Hopf algebroids} \cite{TVD16}. In fact, it follows from \cite[Theorem 3.5]{TVD15} (see also \cite[Example 2.2]{Tim16}) that any regular weak multiplier Hopf algebra gives rise to a regular multiplier Hopf algebroid (although we note that the reader should be careful that some conventions are switched when passing between \cite{TVD15} and \cite{TVD16}!). In the theorem below, we can then simply apply this correspondence with respect to our $I$-partial multiplier Hopf algebras by first using Proposition \ref{PropPartialToWeak}. 

We do not recall here the general definition of a regular multiplier Hopf algebroid, but provide in the next theorem the basic elements of the resulting object (specializing \cite[Theorem 3.5]{TVD15} to our setting), so that the original references can be easily checked when referring to a result in the original works. For consistency of later use, we will refer to the definition of regular multiplier Hopf algebroid as presented in \cite[Definition 2.4]{Tim16}.

\begin{Theorem}\label{TheoEqMAlgbroid}
Let $(A,\Delta)$ be an $I$-partial multiplier Hopf algebra $A$, and set 
\[
B =  C = \msB_A,\qquad s_B = \mbs = t_C,\qquad t_B = \mbt = s_C. 
\]
Then in the notation of \cite[Definition 2.1]{Tim16}, we can identify 
\[
{}_BA \otimes A^B \cong \underset{s}{\oplus}  \Gr{A}{}{}{s}{} \otimes \Gr{A}{s}{}{}{},
\]
and we obtain a left multiplier bialgebroid 
\[
\Delta_B: A \rightarrow \End(\underset{s}{\oplus} \Gr{A}{}{}{s}{} \otimes \Gr{A}{s}{}{}{})
\]
which maps $a\in A$ to left multiplication with $\Delta(a)$. 

Similarly, in the notation of \cite[Definition 2.2]{Tim16}, we can identify 
\[
{}^CA \otimes A_C \cong \underset{s}{\oplus} \Gr{A}{}{}{}{s} \otimes \Gr{A}{}{s}{}{},
\]
and we obtain a right multiplier bialgebroid 
\[
\Delta_C: A \rightarrow \End(\underset{s}{\oplus}  \Gr{A}{}{}{}{s} \otimes \Gr{A}{}{s}{}{})
\]
which maps $a\in A$ to right multiplication with $\Delta(a)$. 

The above left and right multiplier bialgebroids then combine into a regular multiplier Hopf algebroid as defined in \cite[Definition 2.4]{Tim16}. The antipode of this regular multiplier Hopf algebroid coincides with the antipode $S$ of $(A,\Delta)$ as a regular $I$-partial multiplier Hopf algebra.  
\end{Theorem}

\subsection{The hyperobject set}

%Fix $(A,\Delta)$ a regular $I$-partial multiplier Hopf
%We also record here the following result. 

\begin{Prop}
Let $(A,\Delta)$ be a regular $I$-partial multiplier Hopf algebra. Define $r\sim s$ if and only if $\UnitC{r}{s}\neq 0$. Then $\sim$ is an equivalence relation. 
\end{Prop}
\begin{proof}
The relation $\sim$ is reflexive by the first and second condition in Definition \ref{DefPMB}, and it is symmetric by \eqref{EqImUnit}. The transitivity follows from the fact that 
\[
\Delta(\UnitC{r}{s}) = \sum_t \UnitC{r}{t}\otimes \UnitC{t}{s}. 
\]
\end{proof}

\begin{Def}\label{DefHypObj}
We call $\msI = I/\sim$ the \emph{hyperobject set} associated to $(A,\Delta)$. We write 
\[
\msH\msB_A = \Fun_f(\msI) = \{f: \msI \rightarrow \C\mid f \textrm{ has finite support}\},
\]
and refer to $\msH\msB_A$ as the \emph{hyperbase algebra}.
\end{Def}

\begin{Lem}
There is a natural identification 
\[
M(\msH\msB_A) \cong M(\msB_A^{\mbs}) \cap M(\msB_A^{\mbt}) \subseteq M(A),
\]
given by the identification 
\[
f \mapsto \sum_{x\in \msI} \sum_{s\in x} f(x)\Unit_s = \sum_{x\in \msI}\sum_{r\in x} f(x)\Unit^r. 
\]
\end{Lem}
\begin{proof}
If $f \in M(\msH\msB_A)$, we have by the definition of $\sim$ that 
\[
\sum_{x\in \msI} \sum_{s\in x} f(x)\Unit_s =  \sum_{x\in \msI} \sum_{r,s\in x} f(x)\UnitC{r}{s} = \sum_{x\in \msI} \sum_{r\in x} f(x)\Unit^r,
\]
so that we clearly get an embedding of $M(\msH\msB_A)$ into $M(\msB_A^{\mbs}) \cap M(\msB_A^{\mbt})$. On the other hand, if $\mbs(f) = \mbt(g)$ is an element in the latter intersection, we see that 
\[
\delta_{r\sim s} f(s) = \mbs(f)\UnitC{r}{s} = \mbt(g)\UnitC{r}{s} = \delta_{r\sim s} g(r),  
\]
so that for $r\sim s$ and $x = [r]=[s]$ we get $f(s) = g(r) = h(x)$ for $h$ some function on the hyperobject set. Clearly we then have that 
\[
\mbs(f) = \sum_{x\in \msI}\sum_{s\in x} h(x)\Unit_s = \sum_{x\in \msI}\sum_{r\in x} h(x)\Unit^r = \mbt(g),
\]
so our embedding is surjective onto the intersection.
\end{proof}

In the following we will then simply identify $M(\msH\msB_A) \subseteq M(A)$, which is unambiguously defined.

\section{Partial algebraic quantum groups}

\subsection{Invariant functionals}

\begin{Def}\label{DefInvInt}
Let $(A,\Delta)$ be a regular $I$-partial multiplier bialgebra. 

A functional $\varphi: A \rightarrow \C$ is called \emph{left invariant} if it has support on $\underset{r,s}{\oplus}\Gr{A}{r}{r}{s}{s}$ and
\begin{equation}\label{EqLeftInv2}
(\id\otimes \varphi)\Delta(a) = \sum_r \varphi(\Unit^ra)\Unit^r = \sum_r \varphi(a\Unit^r)\Unit^r,\qquad a \in A. 
\end{equation}
A functional $\psi: A \rightarrow \C$ is called \emph{right invariant} if it has support on $\underset{r,s}{\oplus}\Gr{A}{r}{r}{s}{s}$ and
\begin{equation}\label{EqRightInv2}
(\psi\otimes \id)\Delta(a) = \sum_s \psi(\Unit_s a)\Unit_s = \sum_r\psi(a\Unit_s) \Unit_s ,\qquad a\in A. 
\end{equation}
\end{Def}

\begin{Lem}\label{LemSameNot}
Let $(A,\Delta)$ be a regular $I$-partial multiplier Hopf algebra $(A,\Delta)$. Then a functional $\varphi: A \rightarrow \C$ is left invariant if and only if it has support on $\underset{r,s}{\oplus}\Gr{A}{r}{r}{s}{s}$ and satisfies 
\begin{equation}\label{EqLeftInvVDW}
(\id\otimes \varphi)\Delta(a) \in \msB_A^{\mbt},\qquad \forall a\in A.
\end{equation}
Similarly, a functional $\psi: A \rightarrow \C$ is right invariant if and only if it has support on $\underset{r,s}{\oplus}\Gr{A}{r}{r}{s}{s}$ and satisfies 
\begin{equation}\label{EqRightInvVDW}
(\psi\otimes \id)\Delta(a) \in \msB_A^{\mbs},\qquad \forall a\in A.
\end{equation}
\end{Lem}
Note that \eqref{EqLeftInvVDW} (resp. \eqref{EqRightInvVDW}) is the definition of left invariance of $\varphi$ (resp.\ right invariance of $\psi$) on $(A,\Delta)$ as a regular weak multiplier Hopf algebra, as defined in \cite[Definition 1.1]{VDW17}.
\begin{proof} The `only if' direction is clear by Remark \ref{RemSourceTarget}. Conversely, if $\varphi: A \rightarrow \C$ is left invariant in the sense of \cite[Definition 1.1]{VDW17}, then by \cite[Proposition 1.4]{VDW17} and \eqref{EqImUnit} the formula \eqref{EqLeftInv2} holds.
\end{proof}
%CHECK
\begin{Rem}
It follows from \cite[Proposition 5.1]{Kha21} that any left invariant functional on $(A,\Delta)$ as in \cite[Definition 1.1]{VDW17} automatically has support on $\oplus_s\Gr{A}{s}{s}{}{}$, and from the remark following \cite[Proposition 1.8]{VDW17}, it follows that there exists an isomorphism $\theta: I \rightarrow I$ such that $\varphi$ has support on $\oplus_{r,s} \Gr{A}{s}{s}{r}{\theta(r)}$ (see also the discussion in \cite[Lemma 3.3 and Theorem 6.2]{Tim16}). However, it is not clear if automatically $\theta= \id$, so we do not know if the support condition in Definition \ref{DefInvInt} can be dropped. 
\end{Rem}

For $\omega: A \rightarrow \C$ a linear functional, we define the following subspaces, called respectively the \emph{left} and \emph{right} kernel:
\[
\Ker_l(\omega) = \{a\in A \mid \forall b\in A: \omega(ab) = 0\},\qquad \Ker_r(\omega) = \{a\in A\mid \forall b\in A: \omega(ba)=0\}.
\]

The following proposition depends on the particular nature of the base algebra $\msB_A$, but its proof is completely standard, and is in principle contained in the proof of \cite[Lemma 1.16]{VDW17}.  

\begin{Prop}
Let $(A,\Delta)$ be a regular $I$-partial multiplier Hopf algebra, and let $\varphi$ be a left invariant functional. Then there exists $I_{\varphi} \subseteq I$ such that 
\[
\Ker_r(\varphi) = \underset{s\in I \setminus I_{\varphi}}{\oplus} \Gr{A}{}{}{}{s},\qquad \Ker_l(\varphi) = \underset{s\in I \setminus I_{\varphi}}{\oplus} \Gr{A}{}{}{s}{}
\]
Similarly, if $\psi$ is a right invariant functional there exists $I_{\psi}\subseteq I$ such that 
\[
\Ker_r(\psi) = \underset{s\in I \setminus I_{\psi}}{\oplus} \Gr{A}{}{s}{}{},\qquad \Ker_l(\psi) = \underset{s\in I \setminus I_{\psi}}{\oplus} \Gr{A}{s}{}{}{}.
\]
\end{Prop}
\begin{proof}
Assume $a\in \Ker_r(\varphi)$, so $\varphi(Aa)= 0$. Then by the surjectivity of $\can_l$ also 
\[
(\id\otimes \varphi)((A\otimes A)\Delta(a)) = (\id\otimes \varphi)((A\otimes 1)\Delta(Aa)) = 0
\]
and, with $\Delta^{(2)} = (\Delta\otimes \id)\Delta= (\id\otimes \Delta)\Delta$,
\[
(\id\otimes \id\otimes  \varphi)((A\otimes A\otimes A)\Delta^{(2)}(a)) =  (\id\otimes \id\otimes  \varphi)((((A\otimes 1)\Delta(A))\otimes A)\Delta^{(2)}(a)) =0. 
\]
Hence, by non-degeneracy of the multiplication,
\[
(\id\otimes \id\otimes \varphi)((1\otimes 1\otimes A)(\id\otimes \Delta)((A\otimes 1)\Delta(a)))=0.
\]
Taking $d\in A$, $\omega$ a linear functional on $A$ and putting $p = (\omega\otimes \id)((d\otimes 1)\Delta(a))$, we find by applying $\omega\otimes \id$ to the above identity that
\[
(\id\otimes \varphi)((1\otimes A)\Delta(p))=0.
\]
Taking $q\in A$ arbitrary and writing $\Delta(p)(q\otimes 1) = \sum_i a_i \otimes b_i$ with $a_i$ linearly independent, we find by multiplying with $q$ to the right that 
\[
\sum_i \varphi(Ab_i)a_i = 0,
\]
hence $\varphi(Ab_i) = 0$ for all $i$. So also 
\[
0 = \sum_i \varphi(AS(a_i)b_i) = \varphi(AS(q)S(p_{(1)})p_{(2)}) = \sum_s \varphi(AS(q)\Unit_s)\varepsilon(p\Unit_s). 
\]
But $\varepsilon(p\Unit_s) = \varepsilon(p\Unit^s) = \omega(da\Unit_s)$. Since $AS(A) = A^2= A$ and $\omega,d$ were arbitrary, we find
\[
\sum_{s} \varphi(A\Unit_s) a\Unit_s=0. 
\]
Hence if $I_{\varphi} := \{s\in I\mid \varphi(A_s) \neq 0\}$, then 
\[
a\Unit_s \neq 0 \qquad \Rightarrow \qquad s\notin I_{\varphi}. 
\]
Conversely, if $s\in I\setminus I_{\varphi}$, then clearly $\varphi(AA_s)= \varphi(A_s) = 0$. Hence 
\[
\Ker_r(\varphi) = \underset{s\in I\setminus I_{\varphi}}{\oplus} \Gr{A}{}{}{}{s}.
\]

It now follows similarly that there exists $I_{\varphi}' \subseteq I$ such that 
\[
\Ker_l(\varphi) =\underset{s\in I \setminus I_{\varphi}'}{\oplus} \Gr{A}{}{}{s}{}.
\]
But since $\varphi$ has support on $\underset{r,s}{\oplus} \Gr{A}{r}{r}{s}{s}$, we have
\[
\varphi(\Gr{A}{}{}{}{s}) = 0\qquad \iff\qquad \varphi(\Gr{A}{}{}{s}{s})=0\qquad \iff\qquad \varphi(\Gr{A}{}{}{s}{})=0,
\]
hence $I_{\varphi} = I_{\varphi}'$. 

The proof for $\psi$ follows by symmetry, as $\psi$ is a left invariant functional for $A^{\cop}$. 
\end{proof}

In particular, we see that $\varphi$ is faithful if and only if $\varphi(A_s) = \varphi(\Gr{A}{}{}{s}{})\neq 0$ for every $s\in I$, which is just a special case of \cite[Proposition 1.17]{VDW17}. Also note that if $\varphi$ is left invariant, then $\varphi(-\Unit_s)$ is left invariant for any $s\in I$. It follows that we do not need to differentiate between there being a faithful left invariant functional and a faithful \emph{set} of left invariant functionals (cf.~ the discussion in \cite[Section 2]{TVDW22}), since in the latter case we can take $\varphi =\sum_s \varphi_{i_s}(\emdash\Unit_s)$ for some choice of $\varphi_{i_s}$ with $\varphi_{i_s}$ non-zero on $A_s$. 

\begin{Def}\label{DefpartAQG}
We say that a regular $I$-partial multiplier Hopf algebra $(A,\Delta)$ is an $I$-partial algebraic quantum group if there exists a left invariant functional $\varphi$ and a right invariant functional $\psi$ with $\varphi(A_s) \neq 0$ and $\psi(A^s)\neq 0$ for all $s\in I$.
\end{Def}
By the previous observations, we can equivalently ask that $\varphi$ and $\psi$ are faithful functionals.

Combining the foregoing with \cite[Definition 5.2]{Tim16} (and some further elementary observations), we obtain the following.

\begin{Theorem}\label{TheoEquiQgrIPart}
Let $(A,\Delta)$ be an $I$-partial algebraic quantum group. 

Then $(A,\Delta)$ is an \emph{algebraic quantum groupoid} in the sense of \cite[Definition 1.14]{VDW17}, and hence a \emph{measured regular multiplier Hopf algebroid} as in \cite[Definition 5.2]{Tim16}. For the latter, the \emph{base weights} are given (using again notation as in Definition \ref{TheoEqMAlgbroid}) by 
\[
\mu_B(\Unit_r) = 1 = \mu_C(\Unit_r),\qquad \forall r\in I,
\]
and the associated partial left and right integrals are given by 
\[
{}_C\varphi_C = (\id\otimes \varphi)\Delta,\qquad {}_B\psi_B = (\psi \otimes \id)\Delta.
\]
We further have, in the notation of \cite{Tim16}, that
\[
\varphi_B(a) = {}_B\varphi(a) = \sum_r \varphi(a\Unit^r)\Unit_r,\qquad \psi_C(a) = {}_C\psi(a) = \sum_r \psi(\Unit_ra)\Unit^r,\qquad a\in A. 
\]
\end{Theorem}

Following the notation \eqref{EqSourceTarget} of Remark \ref{RemSourceTarget}, we will also write 
\[
\varphi_{\mbs}(a)= \sum_r \varphi(a\Unit^r)\Unit_r,\qquad \psi_{\mbr}(a) = \sum_r \psi(\Unit_r a)\Unit^r,\qquad a\in A. 
\]

\subsection{Modular structure}

In what follows, we let $(A,\Delta)$ be an $I$-partial algebraic quantum group with a fixed faithful left invariant functional $\varphi$ and faithful right invariant functional $\psi$. We also use the notation 
\[
\varphi_S = \varphi \circ S,\qquad \varphi_{S^{-1}} = \varphi \circ S^{-1}. 
\]
Clearly $\varphi_S$ and $\varphi_{S^{-1}}$ are faithful right invariant functionals, as $S$ is anti-comultiplicative.

\begin{Lem}\label{LemInvFunctLtoR}
Let $\varphi'$ be a left invariant functional on $A$, and let $a\in A$. Then there exist $b,c\in A$ and $f\in M(\msB_A)= \Fun(I)$ such that 
\begin{equation}\label{EqIdInvFunctt}
\varphi'(\emdash a) = \psi(\emdash b) = \psi(c\emdash),
\end{equation}
\begin{equation}\label{EqIdInvFunctt2}
 \varphi' =  \varphi(\emdash \mbs(f)) = \varphi(\mbs(f) \emdash).
\end{equation}
Similarly, if $\psi'$ is a right invariant functional on $A$, and $a\in A$, there exist $b,c\in A$ and $g\in M(\msB_A)=\Fun(I)$ with 
\begin{equation}
\psi'(\emdash a) = \varphi(\emdash b)= \varphi(c\emdash),\qquad \psi' = \psi(\emdash \mbt(g)) = \psi(\mbt(g)\emdash). 
\end{equation}
\end{Lem} 
\begin{proof}
This follows from \cite[Proposition 1.6 and Proposition 1.8]{VDW17}. More precisely, \cite[Proposition 1.6]{VDW17} tells us that if we assume 
\begin{equation}\label{EqSpecForm}
a = (\psi\otimes \id)((S^{-1}\otimes \id)(\Delta(q))(p\otimes 1)),
\end{equation}
for some $p,q\in A$, then $\varphi'(\emdash a) = \psi(\emdash b)$ for 
\begin{equation}
b = (\id\otimes \varphi')((\id\otimes S)(\Delta(p))(1\otimes q)), 
\end{equation}
as follows by a direct computation. Now 
\[
\mathrm{span}\{(\psi\otimes \id)((S^{-1}\otimes \id)(\Delta(q))(p\otimes 1))\mid p,q\in A\} = \mathrm{span}\{(\psi\circ S^{-1}\otimes \id)((p\otimes 1)\Delta(q))\mid p,q\in A\}.
\]
The latter set equals $\mathrm{span}\{(\psi\circ S^{-1}\otimes \id)((p \otimes q)\Delta(1))\mid p,q\in A\}$, by the assumptions on a regular $I$-partial multiplier Hopf algebra. By faithfulness of $\psi\circ S^{-1}$, it then easily follows that \emph{any} element of $A$ is a linear combination of elements of the form  \eqref{EqSpecForm}. 

This gives the first equality in \eqref{EqIdInvFunctt}. The second equality in \eqref{EqIdInvFunctt} follow similarly. 

The identity \eqref{EqIdInvFunctt2} is simply a restatement of \cite[Proposition 1.8]{VDW17}. Note that $\mbs(f)$ can be switched from right to left by the support condition on $\varphi$. 

The results for $\psi'$ and $\varphi$ follow by symmetry.
\end{proof}

%TODO: extra formula for \Delta \circ S^2. 
\begin{Prop}\label{PropModAutStruct}
There exist unique algebra isomorphisms $\sigma^{\varphi}$ and $\sigma^{\psi}$ of $A$, called the \emph{modular automorphisms}, such that
\[
\varphi(ab) = \varphi(b\sigma^{\varphi}(a)),\qquad \psi(ab) = \psi(b\sigma^{\psi}(a)),\qquad \forall a,b\in A. 
\]
Moreover, we have the following identities:
\begin{equation}\label{EqModAutUnit}
\sigma^{\varphi}(\Unit_s) = \sigma^{\psi}(\Unit_s) = \Unit_s,\qquad \sigma^{\varphi}(\Unit^s) = \sigma^{\psi}(\Unit^s)= \Unit^s.
\end{equation}
and
\begin{equation}\label{EqFormCoprodMod}
 \Delta\circ \sigma^{\varphi} = (S^2\otimes \sigma^{\varphi})\circ \Delta,\qquad \Delta \circ \sigma^{\psi} = (\sigma^{\psi}\otimes S^{-2})\circ \Delta
\end{equation}
\end{Prop} 
\begin{proof}
The existence and uniqueness of $\sigma^{\varphi},\sigma^{\psi}$ follows from Lemma \ref{LemInvFunctLtoR} and faithfulness of $\varphi,\psi$ (see also \cite[Proposition 1.7]{VDW17}). The behaviour of $\sigma^{\varphi},\sigma^{\psi}$ on the base algebras follows immediately from the condition on the support of $\varphi,\psi$. 

The identity \eqref{EqFormCoprodMod} follows from \cite[Proposition 5.2]{Kha21} (see also \cite[Theorem 6.2]{Tim16} upon using Theorem \ref{TheoEquiQgrIPart}). In fact, all that needs to be observed are the following \emph{strong invariance condition}, as stated in \cite[Proposition 5.1]{VDW17}: 
\begin{equation}
S((\id\otimes \varphi)(\Delta(a)(1\otimes b))) = (\id\otimes \varphi)((1\otimes a)\Delta(b)),\qquad a,b\in A,
\end{equation}
together with its companion, obtained by considering $(A^{\opp},\Delta)$:
\begin{equation}
S^{-1}((\id\otimes \varphi)((1\otimes b)\Delta(c))) = (\id\otimes \varphi)(\Delta(b)(1\otimes c)),\qquad b,c\in A.
\end{equation}
Indeed, taking then $c = \sigma^{\varphi}(a)$ and using the defining property of $\sigma^{\varphi}$ immediately leads to the first identity in \eqref{EqFormCoprodMod}. The other identity follows by symmetry.
\end{proof}

\begin{Prop}\label{PropModElStruct}
There exists a unique invertible element $\delta_{\varphi} \in M(A)$, called the \emph{modular element}, such that 
\begin{equation}\label{EqDefModElDir}
\varphi_S(a) = \varphi(a\delta_{\varphi}),\qquad \varphi_{S^{-1}}(a) = \varphi(\delta_{\varphi} a),\qquad \forall a\in A. 
\end{equation}
We have
\begin{equation}\label{EqModElUnit}
\delta_{\varphi} \Unit_s = \Unit_s\delta_{\varphi},\qquad \delta_{\varphi}\Unit^s = \Unit^s\delta_{\varphi},
\end{equation}
\begin{equation}\label{EqDeldel}
\Delta(\delta_{\varphi}) =(\delta_{\varphi}\otimes \delta_{\varphi})\Delta(1)=  \Delta(1)(\delta_{\varphi}\otimes \delta_{\varphi}),\qquad S(\delta_{\varphi}) = \delta_{\varphi}^{-1},\qquad \varepsilon(\delta_{\varphi} a) = \varepsilon(a\delta_{\varphi}) = \varepsilon(a),\qquad \forall a\in A,
\end{equation}
and moreover
\begin{equation}\label{EqModEl}
(\varphi\otimes \id)\Delta(a) = \sum_s\varphi(a\Unit_s)\delta_{\varphi} \Unit_s,
\end{equation}
\end{Prop} 
\begin{proof}
By Lemma \ref{LemInvFunctLtoR} and faithfulness of $\varphi,\psi$, there exist linear maps $\lambda,\rho: A \rightarrow A$ such that 
\[
\varphi_S(\emdash b) = \varphi(\emdash \rho(b)),\qquad \varphi_S((\sigma^{\varphi})^{-1}(a)\emdash) = \varphi(\emdash \lambda(a)),\qquad a,b\in A. 
\]
Again, faithfulness of $\varphi$ implies $\rho$ is a right multiplier and $\lambda$ is a left multiplier. But then also for any $c\in A$ we have 
\[
\varphi(cb\lambda(a)) = \varphi_S((\sigma^{\varphi})^{-1}(a)cb) = \varphi((\sigma^{\varphi})^{-1}(a)c\rho(b)) = \varphi(c\rho(b)a),
\]
so $\delta_{\varphi} = (\lambda,\rho)$ is a multiplier of $A$. Applying once more the faithfulness of $\varphi$, we see that $\delta_{\varphi}$ must be invertible. The identities \eqref{EqModElUnit} then follow from the defining first equality in \eqref{EqDefModElDir} for $\delta_{\varphi}$, together with the support condition on $\varphi$. 

We now derive \eqref{EqModEl}: compute for $b\in A$ that 
\begin{eqnarray*}
(\varphi\otimes \varphi)((1\otimes b)\Delta(a)) 
&=& (\varphi\otimes \varphi)((1\otimes b)\Delta(1)\Delta(a))\\ 
&=& \varphi(S(b_{(1)})b_{(2)}a_{(1)}) \varphi(b_{(3)}a_{(2)}) \\
&=& \sum_s\varphi(S(b_{(1)})\Unit^s)\varphi(\Unit^s b_{(2)}a)\\
&=& \varphi(S(b_{(1)}))\varphi(b_{(2)}a) \\
&=& \varphi_S(b_{(1)})\varphi(b_{(2)}a) \\
&=& \sum_s \varphi_S(\Unit_s b) \varphi(\Unit_sa)\\
&=& \sum_s \varphi(b \delta_{\varphi} \Unit_s) \varphi(a\Unit_s).
\end{eqnarray*}
Hence \eqref{EqModEl} follows by non-degeneracy and idempotency of $A$ and faithfulness of $\varphi$. 

Applying the foregoing reasoning to $(A^{\opp},\Delta)$, we see that there must also exist an invertible multiplier $\delta_{\varphi'}$ such that 
\[
\varphi_{S^{-1}}(a) = \varphi(\delta_{\varphi}'a),\qquad \forall a\in A. 
\]
But as \eqref{EqModEl} still needs to hold, we deduce that 
\[
\sum_s \varphi(a\Unit_s)\delta_{\varphi}\Unit_s = (\varphi\otimes \id)\Delta(a) = \sum_s \varphi(\Unit_sa)\Unit_s \delta_{\varphi}' = \sum_s \varphi(a\Unit_s) \delta_{\varphi}' \Unit_s,\qquad \forall a\in A.
\]
Faithfulness of $\varphi$ now shows that $\delta_{\varphi} = \delta_{\varphi}'$. 

To see that $S(\delta_{\varphi})= \delta_{\varphi}^{-1}$, we compute for $a\in A$ that  
\[
\varphi(a) = \varphi(S^{-1}(a)\delta_{\varphi}) = \varphi(S^{-1}(S(\delta_{\varphi})a)) = \varphi_{S^{-1}}(S(\delta_{\varphi})a) = \varphi(\delta_{\varphi}S(\delta_{\varphi})a),
\]
so indeed $S(\delta_{\varphi})$ must be inverse to $\delta_{\varphi}$. 

Finally, choosing $t$ fixed  and $a\in \Gr{A}{}{}{}{t}$ such that $\varphi(a)=1$, we get from \eqref{EqModEl} that $\delta_{\varphi}\Unit_t = (\varphi\otimes \id)\Delta(a)$. Applying $\Delta$ and using coassociativity, another application of \eqref{EqModEl} then leads straightforwardly to 
\[
\Delta(\delta_{\varphi}\Unit_t) = (\delta_{\varphi}\otimes \delta_{\varphi})\Delta(1)(1\otimes \Unit_t), 
\]
from which the first identity in \eqref{EqDeldel} follows. The counit identity in \eqref{EqDeldel} follows similarly: choosing $q,t$ fixed and $b \in \Gr{A}{q}{}{t}{}$, we clearly have that both $\varepsilon(\delta_{\varphi}b)$ and $\varepsilon(b)$ are zero when $q\neq t$, while for $q=t$ we compute, with $a \in  \Gr{A}{}{}{t}{t}$ such that $\varphi(a) =1$, 
\[
\varepsilon(\delta_{\varphi}b) = (\varphi \otimes \varepsilon)(\Delta(a)(1\otimes b)) = (\varphi(\Unit_t\emdash) \otimes \varepsilon)(\Delta(a))\varepsilon(b) = \varepsilon(b). 
\]
\end{proof}
\begin{Rem}
Alternatively, one can make use of the results of \cite{Tim16}. Indeed, from \cite[Theorem 6.4]{Tim16}, we have the existence of $\delta_{\varphi}^+$ and $\delta_{\varphi}^-$ such that 
\[
\varphi_S(a) = \varphi(a\delta_{\varphi}^+),\qquad \varphi_{S^{-1}}(a) = \varphi(\delta_{\varphi}^- a),\qquad \forall a\in A. 
\]
Since $S^2,\sigma^{\varphi}$ and $\sigma^{\psi}$ are trivial on the base algebras, it follows easily from the second part of \cite[Theorem 6.4]{Tim16} that $\delta_{\varphi}^+$ and $\delta_{\varphi}^-$ commute elementwise with $\msB_A^{\mbs}$ and $\msB_A^{\mbt}$. Since further ${}_B\varphi = \varphi_B$ in our case (see Theorem \ref{TheoEquiQgrIPart}), it follows from the first part of \cite[Theorem 6.4]{Tim16} that in fact $\delta_{\varphi}^+ = \delta_{\varphi}^-$, proving the existence of the element $\delta_{\varphi}$. The formulas in \eqref{EqDeldel} now follow also from \cite[Theorem 6.4]{Tim16}. 
\end{Rem}

Note that by Lemma \ref{LemInvFunctLtoR}, there also exists an invertible element $\delta_{\varphi,\psi} \in M(A)$ and an invertible element $\mu_{\varphi,\psi} \in M(\msB_A)$ such that 
\begin{equation}\label{EqGenModEl}
\psi(a) = \varphi(a\delta_{\varphi,\psi}),\qquad \delta_{\varphi,\psi} = \mbt(\mu_{\varphi,\psi})\delta_{\varphi},\qquad \forall a\in A. 
\end{equation}
It then follows immediately that 
\[
\sigma^{\psi}(a) = \delta_{\varphi,\psi}\sigma^{\varphi}(a)\delta_{\varphi,\psi}^{-1},\qquad \forall a\in A. 
\]
In particular,
\begin{equation}\label{EqSigmaPhiS}
\sigma^{\varphi_S}(a) = \delta_{\varphi} \sigma^{\varphi}(a)\delta_{\varphi}^{-1},\qquad \forall a\in A.
\end{equation}

\begin{Prop}\label{PropSModAut}
The automorphisms $S^2,\sigma^{\varphi}$ and $\sigma^{\psi}$ all commute pairwise. 
\end{Prop}
\begin{proof}
This follows from \cite[Proposition 5.11]{Kha21} (the quasi-invariance condition required there is clearly fulfilled in our situation). Let us briefly give the argument. 

By \eqref{EqFormCoprodMod}, we have $\Delta \circ \sigma^{\varphi}\circ S^{-2} = (\id\otimes \sigma^{\varphi}\circ S^{-2})\circ \Delta$, so applying $\psi\otimes \id$ to this and using right invariance of $\psi$, together with the fact that $\sigma^{\varphi}$ and $S^2$ are the identity on $\msB_A^{\mbs}$, leads to $\psi \circ \sigma^{\varphi}\circ S^{-2} = \psi$. This easily implies that $\sigma^{\varphi}\circ S^{-2}$ commutes with $\sigma^{\psi}$. Similarly, $\sigma^{\psi} \circ S^2$ commutes with $\sigma^{\varphi}$. 

Applying now \eqref{EqFormCoprodMod} to the commutation between $\sigma^{\psi}$ and $\sigma^{\varphi}\circ S^{-2}$ shows that $(\id\otimes \sigma^{\varphi}\circ S^{-2})$ and $(\id\otimes S^{-2})$ commute on the range of $\Delta$. By fullness, it follows that $\sigma^{\varphi}\circ S^{-2}$ and $S^{-2}$ commute on $A$. Hence $S^2$ and $\sigma^{\varphi}$ commute, and then $S^2$ and $\sigma^{\psi}$ commute as well. By the commutation in the first paragraph, we then find the commutation between $\sigma^{\varphi}$ and $\sigma^{\psi}$.
\end{proof}

%The result is then also true for general $\psi$, since $\sigma^{\psi}$ is obtained from $\sigma^{\varphi_S}$ by conjugating with an invertible element in the multiplier algebra of the target algebra.
%Note that the results in \cite[Section 5.3]{Kha21} 

%%commute up to the scaling element as follows: 
%\begin{equation}\label{EqCommSigmaS}
%\sigma^{\varphi}(S^2(a)) = \mbs(\nu) S^{2}(\sigma^{\varphi}(a))\mbs(\nu^{-1}),\qquad %\forall a\in A.
%\end{equation}
%\end{Prop}

%\begin{Prop}\label{CorSModAut}
%The automorphisms $S^2$ and $\sigma^{\varphi}$ commute up to the scaling element as follows: 
%\begin{equation}\label{EqCommSigmaS}
%\sigma^{\varphi}(S^2(a)) = \mbs(\nu) %S^{2}(\sigma^{\varphi}(a))\mbs(\nu^{-1}),\qquad %\forall a\in A.
%\end{equation}
%\end{Prop}
%\begin{proof}
%This follows from the following straightforward computation for $a,b\in A$:
%\[
%\varphi(b \sigma^{\varphi}(S^2(a))) = \varphi(S^2(a)b) = \varphi(aS^{-2}(b)\mbs(\nu)) = \varphi(S^{-2}(b)\mbs(\nu)\sigma^{\varphi}(a))) = \varphi(b \mbs(\nu) S^{2}(\sigma^{\varphi}(a))\mbs(\nu^{-1})). 
%\]
%\end{proof}

Recall now the hyperobject algebra $\msH\msB_A$ introduced in Definition \ref{DefHypObj}.

\begin{Lem}
There exists a unique invertible element $\nu \in M(\msH\msB_A) = \Fun(\msI)$ such that $\nu$ is central in $M(A)$ and such that
\begin{equation}\label{EqDefIdScalEl}
\varphi\circ S^2 = \varphi(\emdash\nu),\qquad \psi \circ S^2 = \psi(\emdash\nu). 
\end{equation}
Moreover, we have 
\begin{equation}\label{EqPropdelt}
 \sigma^{\varphi}(\delta_{\varphi})  = \delta_{\varphi}\nu^{-1}.
\end{equation}
\end{Lem}
\begin{proof}
Since $\varphi \circ S^2$ is clearly still a faithful left invariant functional, there exists by Lemma \ref{LemInvFunctLtoR} a unique invertible $\nu \in \Fun(I)$ such that $\varphi\circ S^2 = \varphi(\emdash\mbs(\nu))$. 

On the other hand, with respect to the right invariant functional $\psi = \varphi_S$ we have
\[
\varphi_S \circ S^2 = \varphi \circ S^3 =  \varphi(S(\emdash)\mbs(\nu)) = \varphi_S(\mbt(\nu)\emdash) = \varphi_S(\emdash\mbt(\nu)),
\]
where in the last two equalities we used \eqref{EqImUnit} and \eqref{EqModAutUnit}. This identity then still holds for general $\psi$ by applying again Lemma \ref{LemInvFunctLtoR}.

Now we compute for $a,b\in A$ that
\[
\varphi(b \sigma^{\varphi}(S^2(a))) = \varphi(S^2(a)b) = \varphi(aS^{-2}(b)\mbs(\nu)) = \varphi(S^{-2}(b)\mbs(\nu)\sigma^{\varphi}(a))) = \varphi(b \mbs(\nu) S^{2}(\sigma^{\varphi}(a))\mbs(\nu^{-1})). 
\]
Since $\sigma^{\varphi}$ and $S^2$ commute, it follows that necessarily $\mbs(\nu)$ is a central element in $M(A)$.

%The identity in \eqref{EqPropdelt} follows from the computation
Now on the one hand, we have
\[
\varphi(S^2(a)) = \varphi(S(a)\delta_{\varphi}) = \varphi(S(\delta_{\varphi}^{-1}a)) = \varphi(\delta_{\varphi}^{-1}a\delta_{\varphi}),\qquad a\in A,
\]
so that $\sigma^{\varphi}(\delta_{\varphi}) = \delta_{\varphi} \mbs(\nu)^{-1}$. On the other hand, we have for $a\in A$ that 
\[
\varphi(S(a)) = \varphi(S^2(S^{-1}(a))) = \varphi(S^{-1}(a) \mbs(\nu)) = \varphi(S^{-1}(\mbt(\nu)a)) = \varphi(\delta_{\varphi}\mbt(\nu)a),
\]
so that also $\sigma^{\varphi}(\delta_{\varphi}) = \delta_{\varphi} \mbt(\nu)^{-1}$. It follows that $\mbt(\nu)= \mbs(\nu)$, so $\nu$ descends to a function on the hyperobject set and then $\nu = \mbs(\nu)=\mbt(\nu)$ under our identification $\msH\msB_A \subseteq M(A)$.

%and the defining identity \eqref{EqDefIdScalEl}.
\end{proof}

\begin{Def}\label{DefScalingEl}
The element $\nu \in M(\msH\msB_A)$ is called the \emph{scaling element}.
\end{Def}

We end with the following formula, which is proven in \cite[Proposition 5.14]{Kha21}. Although the proof is stated in the context of weak multiplier Hopf $*$-algebras with invariant integrals, the proof holds ad verbatim in the context of our partial multiplier Hopf algebras as we have the formula $S(\delta_{\varphi})= \delta_{\varphi}^{-1}$ at our disposal.

\begin{Prop}
The automorphisms $\sigma^{\varphi},\sigma^{\varphi_S}$ and $S^2$ interact in the following way with the coproduct: 
\begin{equation}\label{EqInterDelSigS}
\Delta \circ S^2 = (\sigma^{\varphi} \otimes (\sigma^{\varphi_S})^{-1})\circ \Delta.
\end{equation}

\end{Prop}

\subsection{Duality}

In the following, we fix an $I$-partial algebraic quantum group $(A,\Delta)$ with faithful left invariant functional $\varphi$ and a faithful right invariant functional $\psi$.

\begin{Def}
We define the vector space $\check{A}$ by 
\[
\check{A} = \{\varphi(a\emdash)\mid a\in A\}. 
\]
\end{Def}

By Proposition \ref{PropModAutStruct} and \eqref{EqGenModEl}, one actually also has
\[
\check{A}= \{\varphi(\emdash a)\mid a\in A\} = \{\psi(a\emdash)\mid a\in A\} = \{\psi(\emdash a)\mid a\in A\}. 
\]
so all choices lead to the same object. 

The space $\check{A}$ is again a regular weak multiplier Hopf algebra with a faithful invariant functional, see \cite[Theorem 2.21]{VDW17}. We endow $\check{A}$ however with the opposite comultiplication to the one loc.\ cit. Concretely, the product and coproduct of $\check{A}$ are defined as follows. By definition of $\check{A}$, we have $(\omega\otimes \id)(\Delta(a))\in A$ and $(\id\otimes \omega)(\Delta(a))\in A$ for $\omega \in \check{A}$ and $a\in A$, and the product on $\check{A}$ is uniquely defined by 
\begin{equation}\label{EqMultDual}
(\omega\cdot \chi)(a) =\chi((\omega\otimes \id)\Delta(a)) = \omega((\id\otimes \chi)\Delta(a)),\qquad a\in A,\omega,\chi\in \check{A}. 
\end{equation}
The coproduct on $\check{A}$ is uniquely determined by 
\begin{equation}\label{EqComultDual1}
(\check{\Delta}(\omega)(\chi\otimes 1),a\otimes b) = (\omega\otimes \chi)((b\otimes 1)\Delta(a)) = \omega(b(\id\otimes \chi)(\Delta(a))),\qquad a,b\in A,\omega,\chi \in \check{A}, 
\end{equation}
\begin{equation}\label{EqComultDual2}
((1\otimes \chi)\check{\Delta}(\omega),a\otimes b) = (\chi\otimes \omega)(\Delta(b)(1\otimes a)) = \omega((\chi \otimes \id)(\Delta(b))a),\qquad a,b\in A,\omega,\chi \in \check{A}. 
\end{equation}
see \cite[Proposition 2.7]{VDW17}. 

We note the following convenient characterisation of the multiplier algebra of $\check{A}$: by \cite[Proposition 2.6]{VDW17}, we have the concrete model
\begin{equation}\label{EqMultiDual}
M(\check{A}) = \{\omega \in \Hom(A,\C)\mid \forall a\in A: (\omega \otimes \id)\Delta(a)\in A \textrm{ and }(\id\otimes \omega)\Delta(a)\in A\},
\end{equation}
where we define left and right multiplication with elements of $\check{A}$ by the formulas in \eqref{EqMultDual}.

\begin{Lem}\label{LemEvExtMult}
The formulas in \eqref{EqComultDual2} still hold if $\omega,\chi \in M(\check{A})$. 
\end{Lem}
\begin{proof}
Note that we can also consider the pairing of $M(\check{A}\otimes \check{A})$ with $A \otimes A$, for example by observing that $A\otimes A$ is also a regular weak multiplier Hopf algebra with faithful invariant functionals for the obvious tensor product comultiplication. It is immediate that $(A\otimes A)^{\vee} = \check{A}\otimes \check{A}$. 

Endow now $A$ with the strict topology, where a net $m_{\alpha}$ converges to $m$ if and only if $m_{\alpha} a = ma$ and $am_{\alpha} = am$ eventually for each $a \in A$. Since $A$ has local units (Corollary \ref{CorLocUn}), each element of $M(A)$ can be approximated in this strict topology by elements in $A$. Similarly, $A\otimes A$ and $\check{A}$ can be endowed with the strict topology, and the same remark holds.  

By the surjectivity of the canonical maps in Definition \ref{DefCanMaps}, it is easily seen that the homomorphism 
\[
\Delta: M(A) \rightarrow M(A\otimes^I A) = \Delta(1) M(A\otimes A)\Delta(1) \subseteq M(A\otimes A)
\]
is strictly continuous, and hence the maps of the form 
\[
M(A) \times M(A) \rightarrow M(A\otimes A),\qquad (a,b) \mapsto \Delta(a)(1\otimes b)
\]
etc.\ are jointly continuous. Applying this remark to $\check{A}$, the conclusion will then follow once we can show that the evaluation maps 
\[
\ev_a: M(\check{A}) \rightarrow \C,\qquad \omega \mapsto \omega(a)
\]
are strictly continuous (for fixed $a\in A$). But this will be the case once we know that 
\[
A = \mathrm{span}\{(\omega \otimes \id)\Delta(b)\mid \omega \in \check{A},b\in A\} = \mathrm{span}\{(\id \otimes \omega)\Delta(b)\mid \omega \in \check{A},b\in A\}. 
\]
By the bijectivity of the canonical maps, this follows since for example
\[
(\id\otimes \psi)((1\otimes A)\Delta(A)) = (\id\otimes \psi)((A\otimes A)\Delta(1)) = A, 
\]
where in the last step we used the faithfulness of $\psi$ (cf.\ Definition \ref{DefpartAQG} and the remark following it). 
\end{proof}

\begin{Theorem}
Endow $\check{A}$ with the $I^4$-grading 
\[
\Gr{\check{A}}{r}{s}{t}{u} = \{\omega\left(\UnitC{t}{u}\emdash\UnitC{r}{s}\right)\mid \omega \in \check{A}, r,s,t,u \in I\}. 
\]
Then $(\check{A},\check{\Delta})$ is an $I$-partial algebraic quantum group.
\end{Theorem}
Note that the ordering of the grading can be remembered as follows: when simply dualizing the algebra and co-algebra structure, the horizontal and vertical grading interchange by a flip around the anti-diagonal, so 
\[
\Gr{\Hom(A,\C)}{t}{u}{r}{s} \cong \Hom(\Gr{A}{t}{r}{u}{s},\C). 
\]
Then, because we are considering the opposite comultiplication, we have one extra flip around the horizontal axis. 

\begin{proof}
Let us first verify that $\check{A}$ is indeed an $I^2$-partial algebra, possibly degenerate and/or non-idempotent. By \cite[Proposition 2.1]{VDW17}, we have for $a\in A$ that
\begin{equation}\label{EqFormVDW}
\omega \cdot \varphi(a\emdash) = \varphi(b\emdash),\qquad b = (\omega \circ S \otimes \id)\Delta(a).
\end{equation}
Now $\varphi(a\emdash) \in \Gr{\check{A}}{m}{n}{p}{q}$ if and only if  $a\in \Gr{A}{m}{p}{n}{q}$. Hence if $\omega \in \Gr{\check{A}}{t}{u}{r}{s}$, we get $b=0$ unless $p=s$ and $m=u$, in which case $b\in \Gr{A}{t}{r}{n}{q}$ and hence $\varphi(b\emdash) \in \Gr{\check{A}}{t}{n}{r}{q}$. So 
\[
\Gr{\check{A}}{t}{u}{r}{s}\cdot  \Gr{\widehat{A}}{m}{n}{p}{q} \subseteq \delta_{p,s}\delta_{m,u} \Gr{\check{A}}{t}{n}{r}{q},
\]
and $\check{A}$ is an $I^2$-partial algebra.

By \cite[Proposition 2.2]{VDW17}, the total algebra $\check{A}$ is non-degenerate and idempotent. By Proposition \ref{PropWeakToPartial} and Lemma \ref{LemSameNot}, it is now sufficient to show that 
\begin{equation}\label{EqImHatDel}
\check{\Delta}(\check{1}) = \sum_s \check{\Unit}_s\otimes \check{\Unit}^s,
\end{equation}
where $\check{1} = \varepsilon$ is the unit of $M(\check{A})$ under the identification \eqref{EqMultiDual}, and where the $\check{\Unit}_s$ and $\check{\Unit}^s$ are the corresponding units for the grading on $\check{A}$. 

Now it is easily checked from Lemma \ref{LemEvExtMult} that $\check{A}$ has 
\begin{equation}\label{EqLocUnitsDual}
\check{\mbunit}_s = \varepsilon(\mbunit_s\emdash),\qquad \check{\mbunit}^s = \varepsilon(\emdash\mbunit_s),\qquad s\in I, 
\end{equation}
where we use again the identification \eqref{EqMultiDual}. Then \eqref{EqImHatDel} follows from applying once more Lemma \ref{LemEvExtMult} upon observing that 
\[
(\check{\Delta}(\check{1}),a\otimes b) = \varepsilon(ba) = \sum_{s} \varepsilon(b \Unit_s a) = \sum_s \varepsilon(b\Unit_s)\varepsilon(\Unit_sa) = (\sum_s \check{\Unit}_s \otimes \check{\Unit}^s,a\otimes b).  
\]
\end{proof}

By \cite[Proposition 2.11]{VDW17}, the antipode of $\check{A}$ is given by 
\begin{equation}\label{EqDualAntipod}
\check{S}(\omega) = \omega \circ S^{-1}. 
\end{equation}
We also note the following formula for further use, cf.\ also the proof of \cite[Proposition 2.7]{VDW17}. 

\begin{Lem} 
Let $\omega = \varphi(\emdash a)$ and $\chi = \varphi(\emdash b)$. Then 
\begin{equation}\label{EqIdCanMapDual}
\check{\Delta}(\omega)(\chi\otimes 1) = \varphi(\emdash b_{(2)})\otimes \varphi(\emdash S^{-1}(b_{(1)})a).
\end{equation}
\end{Lem}
\begin{proof}
By definition, we have for $c,d\in A$ that
\[
(\check{\Delta}(\omega)(\chi\otimes 1))(c\otimes d) = \varphi(dc_{(1)}a)\varphi(c_{(2)}b) = \sum_r \varphi(dc_{(1)}\mbunit_ra)\varphi(c_{(2)}\mbunit_rb).
\]
By the defining property of the antipode, we can write 
\[
\varphi(dc_{(1)}\mbunit_ra)\varphi(c_{(2)}\mbunit^rb) = \varphi(dc_{(1)}(\mbunit^rb)_{(2)}S^{-1}((\mbunit^rb)_{(1)})\mbunit_r a)\varphi(c_{(2)}(\mbunit^rb)_{(3)}),
\]
which by the defining property of a left invariant functional simplifies to 
\begin{eqnarray*}
\sum_s \varphi(d\mbunit^s S^{-1}((\mbunit^rb)_{(1)})\mbunit_ra)\varphi(c(\mbunit^rb)_{(2)}\mbunit^s)& =&  \varphi(d S^{-1}((\mbunit^rb)_{(1)})\mbunit_ra)\varphi(c(\mbunit^rb)_{(2)}) \\
&=&  \varphi(d S^{-1}(b_{(1)})\mbunit_ra)\varphi(c b_{(2)}).
\end{eqnarray*}
Summing over $r$, we then find \eqref{EqIdCanMapDual}.
\end{proof}

Note that, by \cite[Proposition 2.20]{VDW17} and the discussion following it, we have that $\check{A}$ has respectively the faithful left and right invariant functionals $\check{\varphi}$ and $\check{\psi}$ where
\begin{equation}\label{EqDualInvFuncRight}
\check{\varphi}(\omega) = \varepsilon(a),\qquad \omega = \varphi(\emdash a),\qquad a\in A,
\end{equation}
\begin{equation}\label{EqDualInvFuncLeft}
\check{\psi}(\omega) = \varepsilon(a),\qquad \omega = \psi(a\emdash),\qquad a\in A.
\end{equation}

Fixing these invariant functionals, we have the following. We endow also the modular structure of the dual with $\vee$.

\begin{Prop}\label{PropDualModStr}
The modular structure of $\check{A}$ is given by 
\[
\check{\sigma}^{\check{\varphi}}(\varphi(\emdash a)) = \varphi(-S^2(a)\delta_{\varphi}^{-1}),\qquad \check{\delta}_{\check{\varphi}}(a) = \varepsilon(\sigma^{\varphi}(a)),\qquad \check{\nu} = \nu^{-1}.  
\]
\end{Prop}
Note that we view here $\check{\delta}_{\check{\varphi}} \in A^*$ via \eqref{EqMultiDual}.
\begin{proof}
By a similar computation as for \eqref{EqFormVDW} (or by using the formulas in Proposition \ref{PropModAutStruct}), we see that for $a\in A$ and $\omega \in \check{A}$
\[
\omega \cdot \varphi(\emdash a) =  \varphi(\emdash b),\qquad b = (\omega \circ S^{-1}\otimes \id)\Delta(a). 
\]
Hence for $a,b\in A$, we find 
\[
\check{\varphi}(\varphi(\emdash a)\varphi(\emdash b)) = \varphi(S^{-1}(b)a). 
\]
Hence if we write $\check{\sigma}^{\check{\varphi}}(\varphi(\emdash a)) = \varphi(\emdash a')$, we find that 
\[
\varphi(S^{-1}(b)a) = \varphi(S^{-1}(a')b),\qquad \forall b\in A. \]
Using Proposition \ref{PropModElStruct}, we arrive at the expression for $\check{\sigma}^{\check{\varphi}}$ in the statement of the proposition.

Now again by a similar calculation as before, we can also check that for $a \in A$ and $\omega \in M(\check{A})$ one has 
\[
\varphi(\emdash a)\omega = \varphi(\emdash b),\qquad b = (\id\otimes \omega \circ S)(\Delta(a)(1\otimes \delta_{\varphi}^{-1})),
\]
while
\[
\check{S}(\varphi(\emdash a)) = \varphi(\emdash \sigma^{\varphi}(\delta_{\varphi} S(a))). 
\]
Taking $\omega = \check{\delta}_{\check{\varphi}}$ and applying $\check{\varphi}$ to both expressions, we conclude $\check{\delta}_{\check{\varphi}} = \varepsilon \circ \sigma^{\varphi}$. 

The identity $\check{\nu} = \nu^{-1}$ follows from a straightforward computation, using the defining identity in Definition \ref{DefScalingEl}. 
\end{proof}

As is to be expected, one also has the following biduality statement, which is a special case of \cite[Proposition 2.19]{VDW17}. 

\begin{Theorem}\label{TheoDualPQG}
Let $(A,\Delta)$ be an $I$-partial algebraic quantum group, with faithful left and right invariant functionals $\varphi,\psi$. Then we have an isomorphism 
\[
(A,\Delta) \cong (A^{\vee\vee},\Delta^{\vee\vee}),\qquad a \mapsto \check{\varphi}(\emdash \varphi(\emdash a)), 
\]
and under this isomorphism 
\[
\varphi = \varphi^{\vee\vee},\qquad \psi = \psi^{\vee\vee}.
\]
\end{Theorem}

\subsection{Partial $*$-algebraic quantum groups}

\begin{Def}
An $I$-partial $*$-algebra is an $I$-partial algebra $A$ so that the total algebra $A$ is a $*$-algebra and
\[
(\Gr{A}{}{}{r}{s})^* = \Gr{A}{}{}{s}{r}. 
\]

An $I$-partial multiplier Hopf $*$-algebra is a regular $I$-partial multiplier Hopf algebra $(A,\Delta)$ such that $A$ is an $I^2$-partial $*$-algebra and 
\[
\Delta: A \rightarrow \widetilde{M}(A\otimes^I A)\subseteq M(A\otimes A)
\]
is a $*$-algebra homomorphism. 

An $I$-partial $*$-algebraic quantum group is an $I$-partial multiplier Hopf $*$-algebra possessing faithful left invariant functional $\varphi$ which is \emph{positive}: 
\begin{equation}\label{EqPos}
\varphi(a^*a) \geq 0, \forall a\in A. 
\end{equation}
\end{Def}

We will show later that $(A,\Delta)$ then automatically also possesses a right invariant functional $\psi$ which is \emph{positive}, and that we may in fact take $\psi = \varphi_S$.

It is straightforward to extend Proposition \ref{PropPartialToWeak} and Theorem \ref{TheoEqMAlgbroid} to this setting, so an $I$-partial multiplier Hopf $*$-algebra is automatically a weak multiplier Hopf $*$-algebra (cf.\ the remark following \cite[Definition 1.14]{VDW13}).

\begin{Lem}\label{LemSAStat}
Any positive functional $\varphi: A \rightarrow \C$ as above is automaticaly self-adjoint:
\[
\varphi(a^*) = \overline{\varphi(a)},\qquad a\in A. 
\]
\end{Lem} 
\begin{proof}
By polarisation, we have $\varphi(b^*a) = \overline{\varphi(a^*b)}$ for all $a,b\in A$. The statement then follows from the fact that $A$ is idempotent.
\end{proof}

The following lemma follows from \cite[Proposition 4.11]{VDW13}.

\begin{Prop}
Let $(A,\Delta)$ be a partial multiplier Hopf $*$-algebra. Then the antipode $S$ satisfies \begin{equation}\label{EqSStarStar}
    S(S(a)^*)^* = a,\qquad \forall a\in A,
\end{equation} 
\end{Prop}

\begin{Prop}\label{PropdelSA}
Let $(A,\Delta)$ be a partial $*$-algebraic quantum group. Then the modular element $\delta_{\varphi}$ is self-adjoint,
\[
\delta_{\varphi}^* =\delta_{\varphi}.
\]
\end{Prop}
\begin{proof}
This can be derived straightforwardly using directly Lemma \ref{LemSAStat} and \eqref{EqSStarStar}: we have for $a\in A$ that 
\[
\varphi(\delta_{\varphi}^*a) = \overline{\varphi(a^*\delta_{\varphi})} = \overline{\varphi(S(a))} = \varphi(S^{-1}(a^*)) = \varphi(\delta_{\varphi}a^*).
\]
\end{proof}

\begin{Prop}\label{PropDualPos}
Let $(A,\Delta)$ be an $I$-partial $*$-algebraic quantum group. Then also $(\check{A},\check{\Delta})$ is an $I$-partial $*$-algebraic quantum group for the $*$-structure
\begin{equation}\label{EqStarDual}
\omega^*(a) = \overline{\omega(S(a)^*)},\qquad \omega \in \check{A},a\in A. 
\end{equation}
\end{Prop} 
\begin{proof}
By \eqref{EqSStarStar}, the $*$-structure on $\check{A}$ is a well-defined antilinear involution. As $\Delta$ is a $*$-homomorphism and $S$ is anti-comultiplicative, we immediately obtain that $*$ is anti-multiplicative,  and as $S$ is anti-multiplicative we also see immediately that $\check{\Delta}$ is then a $*$-homomorphism. The fact that $\left(\Gr{\check{A}}{r}{s}{t}{u}\right)^* = \Gr{\check{A}}{s}{r}{u}{t}$ follows from the definition of the grading. 

We are left to show that $(\check{A},\check{\Delta})$ has a positive faithful left invariant functional. But consider the faithful left invariant functional $\check{\varphi}$ given by \eqref{EqDualInvFuncRight}. Borrowing the computation in the discussion following \cite[Proposition 2.20]{VDW17} we see that
\[
\check{\varphi}(\omega^*\omega)  = \omega^*(S^{-1}(a)) = \varphi(a^*a),\qquad \omega = \varphi(\emdash a),\qquad a\in A,
\]
so $\check{\varphi}$ is indeed positive. 
\end{proof}

We can now actually show that the scaling element of an $I$-partial $*$-algebraic quantum group is trivial, and that in fact all structure automorphisms are diagonalizable. We follow the arguments of \cite[Section 3]{DCVD10}. We fix an $I$-partial $*$-algebraic quantum group $(A,\Delta)$ with a faithful left invariant positive functional $\varphi$.

\begin{Lem}
If $b \in A$ is non-zero and $n$ is an even integer, then $b^*(\sigma^{\varphi_S})^n(S^{2n}(b))\neq 0$.
\end{Lem}
\begin{proof}
We can follow \cite[Lemma 3.1]{DCVD10} verbatim: by centrality of $\nu$, we get from \eqref{EqSigmaPhiS} and \eqref{EqPropdelt} that $b^*(\sigma^{\varphi_S})^n(S^{2n}(b))=0$ implies
\begin{equation}\label{EqTempEqMod}
b^* \delta_{\varphi}^{n}(((\sigma^{\varphi})^n)(S^{2n}(b)))=0.
\end{equation}
Using that $S^2(a)^* = S^{-2}(a^*)$ and similarly $\sigma^{\varphi}(a)^* = (\sigma^{\varphi})^{-1}(a^*)$ (as follows easily from the selfadjointness of $\varphi$), we obtain by applying $(\sigma^{\varphi})^{-n/2}S^{-n}$ to \eqref{EqTempEqMod} that 
\[
((\sigma^{\varphi})^{n/2} S^{n}(b))^* \delta_{\varphi}^n ((\sigma^{\varphi})^{n/2} S^{n}(b)) = 0.
\]
As we know that $\delta_{\varphi}$ is self-adjoint by Proposition \ref{PropdelSA}, we can apply $\varphi$ to the previous identity and use positivity and faithfulness of $\varphi$ to arrive at 
\[
\delta_{\varphi}^{n/2} ((\sigma^{\varphi})^{n/2} S^n(b)) =0, 
\]
so $b =0$.
\end{proof}

Put 
\[
\kappa = (\sigma^{\varphi})^{-1} \circ S^2,\qquad \rho = \sigma^{\varphi_S} \circ S^2.
\]
\begin{Lem}\label{LemSpanFin}
If $a \in A$, then the linear span of the $\{\kappa^n(a)\mid n \in \Z\}$ is finite-dimensional.
\end{Lem}
\begin{proof}
Pick $a\in \Gr{A}{}{}{r}{}$. Choose a non-zero  $b\in \Gr{A}{r}{}{}{}{}$, and write 
\[
a\otimes b = \sum_{i=1}^n \Delta(p_i)(1\otimes q_i),
\]
which is possible by surjectivity of $\can_r$. We obtain by \eqref{EqInterDelSigS}  that 
\[
\kappa^n(a)\otimes \rho^{n}(b) =   \sum_{i=1}^n \Delta(p_i)(1\otimes \rho^n(q_i))
\]
Multiplying to the left with $1\otimes b^*$, we get 
\[
\kappa^n(a)\otimes b^*\rho^{n}(b) =   \sum_{i=1}^n ((1\otimes b^*) \Delta(p_i))(1\otimes \rho^n(q_i)).
\]
Choose now $a_{ij},b_{ij}\in A$ with 
\[
(1\otimes b^*)\Delta(p_i) = \sum_{j=1}^{m_i} a_{ij}\otimes b_{ij},
\]
and let $L$ be the finite-dimensional space spanned by the $a_{ij}$. We then see that 
\[
\kappa^n(a) \otimes b^*\rho^n(b) \in L \otimes A
\]
for all $n \in \Z$. By the previous lemma, this implies $\kappa^{2n}(a) \in L$ for all integers $n$. But then clearly also the span of the $\kappa^n(a)$ for $n$ ranging over all integers is finite-dimensional.
\end{proof}

\begin{Theorem}
Let $(A,\Delta)$ be an $I$-partial $*$-algebraic quantum group with positive, faithful left invariant functional $\varphi$. Then there exists a joint eigenbasis for the automorphisms $\sigma^{\varphi},S^2$ and left and right multiplication with $\delta_{\varphi}$, with all eigenvalues strictly positive real numbers. Moreover, the scaling element is trivial, i.e.\ $\nu=1$, and $(A,\Delta)$ possesses the positive faithful right invariant functional 
\[
\psi = \varphi_S = \varphi_{S^{-1}}.
\]
\end{Theorem}
\begin{proof}
We begin by showing that left and right multiplications with $\delta_{\varphi}$ are diagonalizable. By Theorem \ref{TheoDualPQG} and Proposition \ref{PropDualPos}, it is sufficient to show this with respect to $(\check{A},\check{\Delta})$ and $\check{\varphi}$. But by Proposition \ref{PropDualModStr}, we have that $\check{\delta}_{\check{\varphi}} = \varepsilon \circ \kappa^{-1}$, and so, using that $\Delta \circ \kappa = (\id\otimes \kappa)\circ \Delta$, we find that
\[
(\varphi(-a)\check{\delta}_{\check{\varphi}})(b) = \varphi(\kappa^{-1}(b)a) = (\varphi\circ S^2)(b\kappa(a))) = \varphi(b\kappa(a)\nu^{-1}),\qquad a,b\in A. 
\]
By Lemma \ref{LemSpanFin}, it follows that for fixed $\omega$ the span of the $\omega \check{\delta}_{\check{\varphi}}^n$ over all integers $n$ form a finite-dimensional space, and hence, by applying the $*$-operation, the same holds for the span of the $\check{\delta}_{\check{\varphi}}^n\omega$. 

Let us now return to $(A,\Delta)$, where by the above we may assume that the span of $\delta_{\varphi}^na$ over all integers $n$ is finite-dimensional for any $a\in A$. Since multiplication with $\delta$ is self-adjoint with respect to the inner product $\langle a,b\rangle = \varphi(a^*b)$, it follows that $A$ has an eigenbasis for left multiplication with $\delta$. Moreover, since $\delta$ commutes with the basis algebra, we may assume that this eigenbasis lives in homogeneous components.

Choose now a non-zero homogeneous eigenvector $b \in \Gr{A}{}{}{r}{}$ for left multiplication with $\delta$, say $\delta b = \lambda b$ for some real number $\lambda \neq 0$. Then by self-adjointness of $\delta_{\varphi}$ and the fact that $\sigma^{\varphi}(\delta_{\varphi}) = \delta_{\varphi}\nu^{-1}$, we find that 
\[
\lambda \varphi(bb^*) = \varphi(bb^*\delta_{\varphi}) = \varphi (\nu \delta_{\varphi}bb^*) = \nu_r \lambda \varphi(bb^*). 
\]
By faithfulness of $\varphi$, we deduce $\nu_r=1$. As this must hold for all $r$, we find $\nu = 1$. 

Note now that it follows from the triviality of $\nu$ that $\kappa,\rho$ and left and right multiplication with $\delta$ all commute. On the other hand, they leave finite-dimensional spaces fixed whose span is $A$, and they are either selfadjoint or anti-selfadjoint operators with respect to the above inner product on $A$. It follows that these operations admit a joint eigenbasis of $A$, and hence the same is true for $S^2,\sigma^{\varphi}$ and left and right multiplication with $\delta_{\varphi}$. 

Finally, we must show that these operators have positive eigenvalues. This follows since by \eqref{EqModEl} we have for $a,b\in A$ that 
\[
0\leq (\varphi\otimes \varphi)((1\otimes b^*)\Delta(a^*a)(1\otimes b)) = \sum_s \varphi(a\Unit_s)\varphi(b^* \delta_{\varphi}\Unit_sb),
\]
so choosing $a$ and $b$ non-zero homogeneous elements with $b$ an eigenvector for left multiplication with $\delta_{\varphi}$, we see that the corresponding eigenvalue must necessarily be positive. The same is then of course true for right multiplication with $\delta_{\varphi}$. It follows in particular that we can define uniquely multipliers $\delta_{\varphi}^z \in A$ for $z\in \C$ such that $\delta_{\varphi}^z b  = \lambda^zb$ for each $\lambda$-eigenvector $b$ for left multiplication with $\delta_{\varphi}$. In particular, $(\delta_{\varphi}^z)^* = \delta_{\varphi}^{\overline{z}}$ and $\sigma^{\varphi}(\delta_{\varphi}^z) = \delta_{\varphi}^z$ (as we already know the scaling element to be trivial). We conclude that 
\[
\varphi_S(a^*a) = \varphi(a^*a\delta_{\varphi}) = \varphi(\delta_{\varphi}^{1/2}a^*a\delta_{\varphi}^{1/2}) \geq 0,
\]
so $\varphi_S$ is a positive functional. 

By the first part of the proof and duality, we then know that also $\kappa$ must have positive eigenvalues, and by considering $(A,\Delta^{\opp})$ with the positive faithful left invariant functional $\varphi_S$, we know that also $\rho$ must be diagonizable with positive eigenvalues. Since $S^2$ and $\sigma_{\varphi}$ can be obtained as combinations of $\kappa,\rho$ and left and right multiplications with $\delta_{\varphi}$, it follows that also the eigenvalues of $S^2,\sigma^{\varphi}$ must be positive.  
\end{proof}

Note that in particular that it follows immediately from this that an $I$-partial $*$-algebraic quantum group is indeed automatically a measured Hopf $*$-algebroid as in \cite[Definition 5.2]{Tim16}.

%TODO: positivity scaling element, positivity $\delta$. 

\subsection{Representation theory}

Partial $*$-algebraic quantum groups have a good  $*$-representation theory. 

\begin{Def}\label{DefModAndRep}
Let $(A,\Delta)$ be an $I$-partial multiplier Hopf algebra. 

By a (left) \emph{module} for $(A,\Delta)$ we mean a unital $A$-module $V$, i.e.~ $AV =V$. By Corollary \ref{CorLocUn} we can extend such a module to a unital left $M(A)$-module, and in particular we can endow $V$ with the natural (vertical) $I^2$-grading $\Gr{V}{r}{}{s}{} := \UnitC{r}{s}V$. 

If $(A,\Delta)$ is a partial multiplier Hopf $^*$-algebra, we mean by a (left) \emph{unitary module} a module $V$ for the underlying partial multiplier Hopf algebra, together with a pre-Hilbert structure on $V$ such that for all $a\in A$ and $v,w\in V$ we have
\[
\langle v,aw\rangle = \langle a^*v,w\rangle,\qquad \forall a\in A,v,w\in V.
\]
We mean by a (left, bounded) \emph{$*$-representation} a Hilbert space $\Hsp$ together with a $*$-homomorphism
\[
A \rightarrow B(\Hsp)
\]
into the bounded operators such that $A\Hsp$ is dense in $\Hsp$. 
\end{Def}

Note that, for $V$ a unitary module, the $\Gr{V}{r}{}{s}{}$ are automatically all mutually orthogonal to each other. 

The unitary modules of partial $*$-algebraic quantum groups are automatically well-behaved.
\begin{Prop}
If $V$ is a  unitary module of a partial $*$-algebraic quantum group $A$, then $A$ automatically acts by bounded operators on $V$.
\end{Prop}
\begin{proof}
We follow the proof as for $*$-algebraic quantum groups \cite{KVD97}. Note that $A$ becomes a pre-Hilbert space by the inner product 
\[
\langle a,b\rangle = \varphi(a^*b). 
\]
Consider the linear map 
\[
T: V\otimes A \rightarrow V\otimes A,\qquad v\otimes a \mapsto a_{(1)}v\otimes a_{(2)}.
\]
This is well-defined by unitality of $V$ as a module, by non-degeneracy of $A$ as an algebra and as $\Delta$ has values in $\widetilde{M}(A\otimes^I A) \subseteq \widetilde{M}(A\otimes A)$. Moreover, $T$ is bounded since, with $\omega_{w,v}(x) = \langle w,xv\rangle$, 
\[
\langle T(dw\otimes b),T(cv\otimes a)\rangle = \omega_{w,v}((\id\otimes \varphi)((d^* \otimes 1)\Delta(b^*a)(c\otimes 1)))  = \sum_r  \langle dw,\mbunit^rcv\rangle  \varphi(\mbunit^rb^*a),
\]
so 
\[
\|T(\sum_i w_i\otimes a_i)\|^2 = \sum_r \|\sum_i \mbunit^r w_i\otimes a_i\mbunit^r\|^2 \leq \|\sum_i w_i\otimes a_i\|^2
\]
as the right multiplications with $\mbunit^r$ on $A$ are self-adjoint projections.  

This shows that the elements 
\[
(\id\otimes \omega_{b,a})(T) = \varphi(b^*a_{(2)})a_{(1)}
\]
act boundedly on $V$. By surjectivity of the map $\can_l^c$, we see that the $\varphi(c)d$ act boundedly for $r\in I$, $c\in \Gr{A}{}{r}{}{}$ and $d\in \Gr{A}{}{}{}{r}$. As $\varphi$ is faithful, $\varphi(A^r) = \C$ for all $r$, hence it follows that $A$ acts boundedly on $V$.  
\end{proof}

\begin{Cor}
Any unitary module of a partial $*$-algebraic quantum group completes in a unique way to a $*$-representation.
\end{Cor}

Conversely, if $(\Hsp,\pi)$ is a $*$-representation of $(A,\Delta)$, we have that $V_{\pi} = A\Hsp$ is a unitary $A$-module by the idempotency of $A$, and clearly $(\Hsp,\pi)$ is the completion of $V_{\pi}$. We call $V_{\pi}$ the \emph{canonical unitary module} associated to the $*$-representation $(\Hsp,\pi)$. 

%the existence of local units in $A$ \cite[Proposition 4.9]{VDW20}, 

\section{Partial $*$-algebraic quantum groups of compact type and of discrete type, and the Drinfeld double construction}

\subsection{Partial $*$-algebraic quantum groups of compact type and discrete type}

\begin{Def}
We say that a regular partial multiplier Hopf ($*$-)algebra $(A,\Delta)$ is a regular \emph{partial Hopf ($*$-)algebra} if $\UnitC{r}{s}\in A$ for all $r,s$.

We say that a partial ($*$-)algebraic quantum group $(A,\Delta)$ is of \emph{compact type} if the underlying regular partial multiplier Hopf ($*$-)algebra is a regular partial Hopf ($*$-)algebra.

We say that a partial ($*$-)algebraic quantum group $(A,\Delta)$ is of \emph{discrete type} if it arises as the dual of a partial ($*$-)algebraic quantum group of compact type. 
\end{Def}

It is easily verified that the above definition of a regular partial Hopf algebra coincides with the one in \cite[Definition 1.11]{DCT15}. We also have the following. 

\begin{Prop}
A partial $*$-algebraic quantum group of compact type coincides with the notion of a partial compact quantum group as introduced in \cite[Definition 1.18]{DCT15}. 
\end{Prop}
\begin{proof}
The only thing to verify is that we can choose our left invariant functional $\varphi$ such that $\varphi$ is also right invariant and $\varphi(\UnitC{r}{r}) = 1$ for all $r\in I$. 

However, choose $\varphi$ a general faithful positive left invariant functional. If $r\sim s$ (cf.\ Definition \ref{DefHypObj}), it follows from the left invariance of $\varphi$, applied to $\Delta(\UnitC{r}{s})$, that
\[
\varphi(\UnitC{t}{s})\UnitC{r}{t} = \varphi(\UnitC{r}{s})\UnitC{r}{t},\qquad \forall t\in I. 
\]
In particular, $\varphi(\UnitC{t}{s}) =\varphi( \UnitC{r}{s})$ whenever $r \sim t$. Furthermore, by faithfulness of $\varphi$ we can choose for each $s$ at least one $r$ (necessarily with $r\sim s$) such that
\[
0 \neq \varphi(\Gr{A}{}{r}{}{s}) = \varphi(\Gr{A}{r}{r}{s}{s}). 
\]
From the positivity of $\varphi$, the Cauchy-Schwarz inequality then forces $\varphi(\UnitC{r}{s}) \neq 0$. 

It follows from the above that for each $s\in I$, there exists $0<\gamma_s$ such that $\varphi(\UnitC{r}{s}) = \gamma_s$ for all $r \sim s$. Rescaling with $\gamma_s$, it follows that we may assume $\gamma_s =1$ for all $s$. 

It remains to show that $\varphi$ is also right invariant. But if $\psi$ is any positive right invariant functional, we can similarly rescale $\psi$ so that $\psi(\UnitC{r}{s}) = 1$ whenever $r\sim s$. We then find that for $0\neq a\in \Gr{A}{}{r}{}{s}$
\[
\varphi(a) = \varphi(a)\psi(\UnitC{r}{s}) = (\psi \otimes \varphi)((\Unit_s \otimes 1)\Delta(a)) = (\psi \otimes \varphi)((1\otimes \Unit^s)\Delta(a)) = \psi(a) \varphi(\UnitC{s}{s}) = \psi(a), 
\]
so $\varphi= \psi$ and $\varphi$ is right invariant. 
\end{proof}

\begin{Rem}
In correspondence with the above result, we will then also refer to a partial $*$-algebraic quantum group of discrete type as a \emph{partial discrete quantum group}.
\end{Rem}

Whenever we work with a partial compact quantum group as above, we will assume that our positive left invariant functional $\varphi$ is \emph{normalized} as above, so $\varphi(\UnitC{r}{s}) = 1$ whenever $r\sim s$. Then $\varphi$ is also right invariant, so we refer to $\varphi$ as an \emph{invariant functional} or \emph{invariant integral}. As an immediate corollary, we obtain:

\begin{Cor}\label{CorModScalingEl}
Let $(A,\Delta)$ be a partial $*$-algebraic quantum group of compact type. If $\varphi$ is a normalized faithful positive left invariant functional, then $\varphi$ is also a normalized positive right invariant functional, and 
\[
\varphi = \varphi \circ S.
\]
In particular, we have that the modular element is trivial,  $\delta_{\varphi} = 1$.
\end{Cor}

We will show that if $(A,\Delta)$ is a partial algebraic quantum group of compact type, the module theory of $(\check{A},\check{\Delta})$ can be linked to a comodule theory for $(A,\Delta)$. The latter was defined in \cite{DCT15}, but we will take a slightly different approach. We will comment on this in Remark \ref{RemCompCorep}. Note that a general theory of comodules for weak multiplier Hopf algebras was developed in \cite{Boh14}, but we will not need this generality for our purposes. 

If $V$ is a horizontally $I^2$-graded vector space $V = \underset{r,s}{\oplus} \Gr{V}{}{}{r}{s}$, we view $V$ as a $\Fun(I)$-bimodule via 
\[
\lambda(f)v = f \cdot v   = f(r)v,\qquad \rho(g)v= v \cdot g =  g(s)v,\qquad v\in \Gr{V}{}{}{r}{s}. 
\]
If $V$ is a pre-Hilbert space, we assume by default that the $\Gr{V}{}{}{r}{s}$ are mutually orthogonal to each other. 

It will also be convenient to use the following notation in the case of partial Hopf algebras:
\[
\Gr{\Delta}{}{}{k}{l}(a) = (\mbunit_k\otimes \mbunit^k)\Delta(a)(\mbunit_l\otimes \mbunit^l) =(\UnitC{r}{k}\otimes \UnitC{k}{t})\Delta(a)(\UnitC{s}{l}\otimes \UnitC{l}{u}) \in \Gr{A}{}{}{k}{l}\otimes \Gr{A}{k}{l}{}{},\qquad a\in \Gr{A}{r}{s}{t}{u} \subseteq A.
\]

\begin{Def}\label{DefComodPar}
Let $(A,\Delta)$ be a partial Hopf algebra.

A (right) \emph{comodule} for $(A,\Delta)$ consists of an $I^2$-graded vector space $V = \underset{r,s}{\oplus}\Gr{V}{}{}{r}{s}$ and a collection of linear maps $\delta = \{\Gr{\delta}{}{}{r}{s}\}$ with 
\[
\Gr{\delta}{}{}{r}{s}: V \rightarrow V\otimes A,\qquad v\mapsto v_{(0;rs)}\otimes v_{(1;rs)}
\]
 such that
\begin{enumerate}
\item\label{EqRangeDel} $\Gr{\delta}{}{}{r}{s}(\Gr{V}{}{}{k}{l})\subseteq \Gr{V}{}{}{r}{s}\otimes \Gr{A}{r}{s}{k}{l}$ for all $k,l,r,s$,
\item\label{EqComod} we have for all $k,l,r,s\in I$ that
\[
(\id\otimes \Gr{\Delta}{}{}{k}{l})\Gr{\delta}{}{}{r}{s}= (\Gr{\delta}{}{}{r}{s}\otimes \id)\Gr{\delta}{}{}{k}{l},\qquad (\id\otimes \varepsilon)\Gr{\delta}{}{}{r}{s} = \lambda(\mbunit_r)\rho(\mbunit_s).
\]
\item\label{EqFinSum1} for fixed $v$ and $s$ there are only finitely many $r$ with $\Gr{\delta}{}{}{r}{s}(v)$ non-zero, leading to a map  $\Gr{\delta}{}{}{}{s}: V \rightarrow V\otimes A$ with $\Gr{\delta}{}{}{}{s}(v) = \sum_r \Gr{\delta}{}{}{r}{s}(v)$,
\item\label{EqFinSum2} for fixed $v$ and $r$ there are only finitely many $s$ with $\Gr{\delta}{}{}{r}{s}(v)$ non-zero, leading to a map $\Gr{\delta}{}{}{r}{}: V \rightarrow V\otimes A$ with $\Gr{\delta}{}{}{r}{}(v) = \sum_s \Gr{\delta}{}{}{r}{s}(v)$.
\end{enumerate}

We then call the linear map 
\begin{equation}\label{EqDefCorep}
X=X_{\delta}: V \otimes A \rightarrow V\otimes A,\qquad v\otimes a \mapsto \sum_{rs}\Gr{\delta}{}{}{r}{s}(v)(1\otimes a)
\end{equation}
the associated \emph{corepresentation}. 

If $(A,\Delta)$ is a partial Hopf $*$-algebra, we call (right) \emph{unitary comodule} a comodule $(V,\delta)$ for the underlying partial Hopf algebra such that $V$ is a pre-Hilbert space and the associated corepresentation $X$ is \emph{isometric}, in the sense that with respect to the $A$-valued inner product on $V\otimes A$ given by
\[
\langle w\otimes b,v \otimes a\rangle_A = \langle w,v\rangle b^*a,\qquad v,w\in V,a,b\in A,
\] 
it satisfies 
\begin{equation}\label{EqUnitaryX}
\langle X(w\otimes b), X(v\otimes a)\rangle_A = \sum_t \langle v\mbunit_t,w\mbunit_t\rangle b^*\mbunit_ta,\qquad v,w\in V,a,b\in A.
\end{equation}
\end{Def}
Note that the linear map in \eqref{EqDefCorep} is indeed well-defined by condition \eqref{EqFinSum1} in the above definition, as we can write 
\[
X(v\otimes a) = \delta_s(v)(1\otimes a),\qquad v\in V,a\in \Gr{A}{}{}{s}{}. 
\]
Then $X$ vanishes on the $\Gr{V}{}{}{}{s}\otimes \Gr{A}{}{}{s'}{}$ with $s\neq s'$ and has range in $\oplus_s \Gr{V}{}{}{s}{}\otimes \Gr{A}{s}{}{}{}$. By condition \eqref{EqFinSum2}, we can also define a map 
\begin{equation}\label{EqDefAdjX}
X^-: V \otimes A \rightarrow V\otimes A,\qquad v\otimes a \mapsto \sum_{rs} v_{(0;rs)}\otimes S(v_{(1;rs)})a. 
\end{equation}
The defining relations for the antipode and comodule structure easily show that then
\begin{equation}\label{EqAdjCorepAlt}
X^-X(v\otimes a) = \sum_s v\Unit_s\otimes \Unit_s a,\qquad XX^-(v\otimes a) = \sum_s \Unit_sv\otimes \Unit^sa,\qquad v\in V,a\in A,
\end{equation}
so $X^-$ is a pseudo-inverse to $X$.

If $(A,\Delta)$ is a partial Hopf $*$-algebra, we can use the alternative notation $(w\otimes b)^*(v \otimes a) =  \langle w\otimes b,v\otimes a\rangle_A$ to rewrite the isometry condition for $X$ as
\begin{equation}\label{EqIdUnitComod}
\sum_s \delta_s(v)^*\delta_s(w) = \sum_t \langle v\mbunit_t,w\mbunit_t\rangle \mbunit_t,\qquad v,w\in V,
\end{equation}
where the left hand side, being a sum of elements in distinct $\Gr{A}{s}{s}{}{}$, converges in $M(A)$ for the strict topology. The following lemma reinterprets this identity in the setting of corepresentations.

\begin{Lem}
Let $(A,\Delta)$ be a  partial Hopf $*$-algebra, and let $(V,\delta)$ be a comodule for the underlying partial Hopf algebra. Assume $V$ is a pre-Hilbert space with the $\Gr{V}{}{}{r}{s}$ mutually orthogonal. Then $(V,\delta)$ is a unitary comodule if and only if there exists an adjoint $X^*: V\otimes A \rightarrow V\otimes A$ for the associated corepresentation $X$ with respect to the $A$-valued inner product on $V\otimes A$, and 
\begin{equation}\label{EqAdjCorep}
X^*X = \sum_s \rho(\mbunit_s) \otimes \mbunit_s,\qquad XX^* = \sum_r \lambda(\mbunit_r)\otimes \mbunit^r.  
\end{equation}
\end{Lem} 
\begin{proof}
If \eqref{EqAdjCorep} holds, then by definition of the adjoint we get 
\[
\langle X(w\otimes b), X(v\otimes a)\rangle_A = \langle w\otimes b,X^*X(v\otimes a)\rangle_A = \sum_t \langle v\mbunit_t,w\mbunit_t\rangle b^*\mbunit_ta,\qquad v,w\in V,a,b\in A,
\]
hence $\delta$ defines a unitary comodule. 

Conversely, if $\delta$ is a unitary comodule, it easily follows from \eqref{EqUnitaryX} and the kernel/image conditions of $X$ that $X^-$ is an adjoint for $X$, and \eqref{EqAdjCorep} follows from \eqref{EqAdjCorepAlt}. 
\end{proof}

%\textcolor{red}{TO DO: check papers of Wang for similar kind of duality.}
We can now state the following duality. We will need to switch between horizontal and vertical $I^2$-gradings. If $V$ is horizontally graded, then the associated compatible vertical grading is given by 
\[
\Gr{V}{r}{}{s}{} = \Gr{V}{}{}{s}{r} 
\]
and vice versa.

\begin{Prop}\label{Prop1to1RepCorep}
Let $(A,\Delta)$ be a partial algebraic quantum group of compact type, and let $V$ be an $I^2$-graded vector space.

Then there is a one-to-one correspondence between left $\check{A}$-module structures on $V$ (compatible with the vertical grading) and right $A$-comodule structures $\delta$ on $V$ (compatible with the horizontal grading), determined by 
\begin{equation}\label{EqModFromComod}
\omega \cdot v = \sum_{r,s} (\id\otimes \omega)(\Gr{\delta}{}{}{r}{s}(v)),\qquad v\in V,\omega \in \check{A}. 
\end{equation}
If $V$ is an $I^2$-graded pre-Hilbert space, this correspondence respects the unitarity of the modules/comodules. 
\end{Prop}
\begin{proof}
Clearly $V$ becomes a well-defined $\check{A}$-module by \eqref{EqModFromComod} for the total algebra $\check{A}$. To see that this module is non-degenerate, choose $v\in \Gr{V}{}{}{k}{}$ and choose $a\in \Gr{A}{k}{}{}{}{}$ with $\psi(a)=1$. Since the corepresentation $X$ associated to $\delta$ has range $\oplus_s {}_sV \otimes \Gr{A}{s}{}{}{}$ by \eqref{EqAdjCorepAlt}, it follows that we can write $v\otimes a = \sum_i X(v_i \otimes a_i) = \sum_i \sum_{rs} \delta_{rs}(v_i)(1\otimes a_i)$, and hence
\[
v = (\id\otimes \psi)(v\otimes a) = \sum_i \psi(\emdash a_i)\cdot v_i,
\]
proving non-degeneracy. It is then easily checked that $\Gr{\check{A}}{r}{s}{t}{u}$ maps $\Gr{V}{s}{}{u}{}$ into $\Gr{V}{r}{}{t}{}$ and is zero elsewhere. Hence the given grading on $V$ coincides with its grading as a module for the partial algebraic quantum group $\check{A}$. If moreover $A$ is a $*$-algebraic quantum group and $\delta$ is unitary, it follows from the fact that $X^-$ given by \eqref{EqDefAdjX} is the adjoint of $X$ that the module is unitary, taking into account the $*$-structure \eqref{EqStarDual} on $\check{A}$. 

Assume now conversely that $V$ is a (non-degenerate) $\check{A}$-module. Then by non-degeneracy of the module $V$ and the property of the comultiplication of a weak multiplier Hopf algebra, we have a well-defined map 
\[
\msX: V \otimes \check{A}\rightarrow V\otimes \check{A},\qquad v\otimes \omega \mapsto \omega_{(1)}v \otimes \omega_{(2)}.
\]
Let $\msF: A\rightarrow \check{A}$ be the bijective linear map $a\mapsto \varphi(-a)$, and define
\[
X = (\id\otimes \msF^{-1})\msX(\id\otimes \msF): V\otimes A \rightarrow V\otimes A.
\]
We claim that $X$ arises as the corepresentation associated to a (unique) comodule $\delta$ implementing the original $\check{A}$-module. In fact, note that if 
\[
\msX^-: V \otimes \check{A}\rightarrow V\otimes \check{A},\qquad v\otimes \omega \mapsto \check{S}^{-1}(\omega_{(1)})v \otimes \omega_{(2)},
\] 
then the properties of the antipode give
\begin{equation}\label{EqCoinvmsX}
\msX^- \msX(v\otimes \omega) = \sum_s \check{\mbunit}^sv\otimes \omega\check{\mbunit}^s,\qquad \msX\msX^-(v\otimes \omega) =\sum_s \check{\mbunit}_sv\otimes \check{\mbunit}^s\omega.
\end{equation}
It follows in particular that, using the horizontal grading on $V$, $X$ is zero on the $\Gr{V}{}{}{}{r}\otimes \Gr{A}{}{}{s}{}$ with $r\neq s$ and $X$ defines a linear bijection between $\oplus_r \Gr{V}{}{}{}{r}\otimes \Gr{A}{}{}{r}{}$ and $\oplus_s \Gr{V}{}{}{s}{} \otimes \Gr{A}{s}{}{}{}$. 

Now by invoking \eqref{EqIdCanMapDual}, it is not hard to see that $X$ satisfies
\[
X(v\otimes ab) = (X(v\otimes a))(1\otimes b),\qquad v\in V,a,b\in A. 
\]
If we hence put
\[
\Gr{\delta}{}{}{}{s}(v) := X(v\otimes \UnitC{s}{t}),\qquad v\in \Gr{V}{}{}{}{t},
\]
it follows that we can write 
\begin{equation}\label{EqComodCorep}
X(v\otimes a) = \sum_s \delta_s(v)(1\otimes a),\qquad v\in V,a\in A. 
\end{equation}
Putting $\Gr{\delta}{}{}{r}{s}(v) = (1\otimes \mbunit_r)\delta_s(v)$, it then follows easily that the $\Gr{\delta}{}{}{r}{s}$ satisfy \eqref{EqRangeDel} in Definition \ref{DefComodPar}, and by \eqref{EqComodCorep} and the definition of $X$ that the original module structure is given by the formula \eqref{EqModFromComod}. The module property for the latter then immediately gives that  the $\Gr{\delta}{}{}{r}{s}$ satisfy the comodule property \eqref{EqComod}. Property \eqref{EqFinSum1} holds by construction. 

Finally, note that by \eqref{EqCoinvmsX} the map $X$ has a coinverse $X^-$. It is not hard, using the properties of the antipode, that we must then have 
\[
\sum_{rs} v_{(0;rs)}\otimes bS(v_{(1;rs)})a = (1\otimes b)(X^-(v\otimes a)),\qquad v\in V,a,b\in A. 
\]
It follows that, fixing $v\in V$ and $r\in I$, there can only be finitely many $s\in I$ with $\Gr{\delta}{}{}{r}{s}(v)\neq 0$, proving that also Property \eqref{EqFinSum2} holds. 

If the original module was a unitary module, it is not hard to see from \eqref{EqModFromComod} that $X^{-}$ is the adjoint of $X$, proving that then also $\delta$ is a unitary comodule. 
\end{proof}

%\begin{Rem}
%Note that to define the dual $\check{A}$-module structure on $V$ with the proper grading properties, property \eqref{EqFinSum2} in Definition \ref{DefComodPar} was not needed, and can hence, by the above reconstruction result, be dropped from the axioms for comodules of partial algebraic quantum groups of compact type.
%\end{Rem} 

\begin{Rem}\label{RemCompCorep} 
Recall from \cite[Definition 2.1]{DCT15} that an $I$-bigraded vector space $V = \oplus_{r,s}\Gr{V}{}{}{r}{s}$ is called \emph{rcfd} (row-and column finite-dimensional) if all $\oplus_r \Gr{V}{}{}{r}{s}$ and $\oplus_s \Gr{V}{}{}{r}{s}$ are finite-dimensional. A \emph{corepresentation $\msX$} of an $I$-partial Hopf algebra $A$ on $V$ was then defined in \cite[Definition 2.2]{DCT15} as a collection of elements 
\[
\Gr{X}{k}{l}{m}{n} \in \Hom_{\C}(\Gr{V}{}{}{m}{n},\Gr{V}{}{}{k}{l}) \otimes \Gr{A}{k}{l}{m}{n} = \Hom_{\C}(\Gr{V}{}{}{m}{n},\Gr{V}{}{}{k}{l} \otimes \Gr{A}{k}{l}{m}{n})
\]
with 
\[
(\id\otimes\Gr{ \Delta}{}{}{p}{q})(\Gr{X}{k}{l}{m}{n}) = \left(\Gr{X}{k}{l}{p}{q}\right)_{12}\left(\Gr{X}{p}{q}{m}{n}\right)_{13},\qquad (\id\otimes \varepsilon)(\Gr{X}{k}{l}{m}{n}) = \delta_{k,m}\delta_{l,n}\id_{\Gr{V}{}{}{k}{l}}
\]
(we have switched here from left to right corepresentation for compatibility with the conventions in this paper). Then putting 
\[
\Gr{\delta}{}{}{r}{s}(v) = \sum_{mn} \Gr{X}{r}{s}{m}{n} (\Unit_m v\Unit_n) \in \oplus_{mn} (\Gr{V}{}{}{r}{s} \otimes \Gr{A}{r}{s}{m}{n}),
\]
it is easily verified that $\delta$ satisfies Definition \ref{DefComodPar}, where \eqref{EqFinSum1} and \eqref{EqFinSum2} are automatic as we have assumed $V$ to be rcfd.  Clearly the associated corepresentation $X= X_{\delta}$ then coincides with the total corepresentation for $\msX$ as defined in \cite[Definition 2.4]{DCT15}. 

Under this correspondence, also the notion of unitarity of a corepresentation, \cite[Definition 2.26]{DCT15}, coincides. 
\end{Rem}

\subsection{The Drinfeld double}

The Drinfeld double for finite-dimensional weak Hopf algebras was defined in \cite{Nen02}. We generalize this construction in the setting of partial algebraic quantum groups of compact type. A full generalization to the setting of arbitrary partial algebraic quantum groups or even arbitrary algebraic quantum groupoids \cites{Tim16,VDW17,WZ17,TVDW22} seems feasible, but we prefer to restrict to this simpler setting which will be sufficient for the applications we have in mind (these will be dealt with in a separate paper). Note that for (pairs of) multiplier Hopf algebras, the Drinfeld double was constructed in \cite{DrVD01}, see also \cite{DeVD04}. 

Fix in what follows an algebraic quantum group of compact type $(A,\Delta)$, with fixed faithful left invariant functional $\varphi$ and faithful right invariant functional $\psi = \varphi \circ S$. 
% normalized so that $\varphi\left(\UnitC{r}{s}\right)=1$ whenever $\UnitC{r}{s}\neq 0$. Put $\psi = \varphi \circ S$. 

\begin{Def}
We define $\msU$ as the space of linear functionals $A\rightarrow \C$ such that, if $r,s$ are fixed, there are only finitely many $t,u$ with $\omega(\UnitC{r}{t} - \UnitC{s}{u})$ or $\omega(\UnitC{t}{r}-\UnitC{u}{s})$ non-zero. 
\end{Def}

It is easily seen that the conditions on the functionals in $\msU$ allow to define an associative algebra structure on $\msU$ by
\[
(\omega \cdot \chi)(a) = (\omega \otimes \chi)\Delta(a),\qquad \omega,\chi \in \msU,a\in A,
\]
and in fact we have $\msU \subseteq M(\check{A})$ by \eqref{EqMultiDual}. Moreover, $\msU$ contains the unit $\varepsilon$. 
We denote as before $\UnitCC{r}{s} = \varepsilon(\Unit_s-\Unit_r)$ for the partial units. Then
\begin{equation}\label{EqPartUnitDual}
\UnitCC{r}{s}\omega = \omega(\Unit^s\emdash\Unit^r),\qquad \omega\UnitCC{r}{s} = \omega(\Unit_s\emdash\Unit_r).
\end{equation}

\begin{Def}
We define $\msD(A)$ to be the universal unital algebra generated by $A$ and $\msU$ with relations $\UnitC{r}{s} = \UnitCC{r}{s}$ and interchange relations
\begin{equation}\label{EqInterRel}
\omega g =  g_{(2)}\omega(g_{(3)} \emdash S^{-1}(g_{(1)})). 
\end{equation}
\end{Def}

Note that since $\msU$ is unital, there is no ambiguity in defining the free product of $A$ and $\msU$ as the usual unital free product of the unitalisation $A_1$ and $\msU$: in the latter the added unit of $A_1$ gets identified with the already existing unit of $\msU$, and this free product contains copies of $A$ and $\msU$. Moreover, the interchange relation \eqref{EqInterRel} is meaningful: if $g\in \Gr{A}{r}{s}{t}{u}$, then by the condition on $\omega$ there exist finitely many $x,y\in I$ such that $\omega(\Unit_t \emdash \Unit_r) = \sum_{xy} \omega(\UnitC{x}{t}\emdash \UnitC{y}{r})$. But for each of these $x,y$, we have
\[
\UnitC{r}{y}g_{(1)}\otimes \UnitC{y}{x}g_{(2)}\otimes \UnitC{x}{t}g_{(3)} \in A\otimes A\otimes A,
\]
hence we can write 
\[
\omega g =  \sum_{xy} \UnitC{y}{x}g_{(2)}\omega(\UnitC{x}{t}g_{(3)} \emdash S^{-1}(\UnitC{r}{y}g_{(1)}))
\]
with the linear combination on the right finite. Moreover, for each $a\in A$ and $\omega \in \msU$ we have that $\omega(a\emdash)$ and $\omega(\emdash a)$ again lie in $\msU$. 

\begin{Lem}\label{LemMultOppOrd}
Inside $\msD(A)$, the following relation holds: 
\begin{equation}\label{EqInterRelRev}
g\omega = \omega(S^{-1}(g_{(3)})\emdash g_{(1)})g_{(2)}.
\end{equation}
\end{Lem} 
\begin{proof}
A similar reasoning as before shows that the above expression is well-defined. In fact, if $g\in \Gr{A}{r}{s}{t}{u}$, then we can write 
\[
\omega(S^{-1}(g_{(3)})\emdash g_{(1)})g_{(2)} = \sum_{xy} \omega(S^{-1}(g_{(3)}\UnitC{x}{u})\emdash  g_{(1)}\UnitC{s}{y})g_{(2)}\UnitC{y}{x},
\]
where only finitely many terms on the right, only depending on $s,u$ and $\omega$, are non-zero. 

Now we compute that 
\[
\omega(S^{-1}(g_{(3)})\emdash g_{(1)})g_{(2)} =  g_{(3)}\omega(S^{-1}(g_{(5)})g_{(4)}\emdash S^{-1}(g_{(2)})g_{(1)}) = g \omega(\Unit^u \emdash  \Unit^s)
\]
where we must note that the expressions inside are meaningful using again (tedious) local unit arguments. By \eqref{EqPartUnitDual} and the relations of $\msD(A)$, we can write $\omega(\Unit^u\emdash \Unit^s) = \UnitCC{s}{u}\omega = \UnitC{s}{u}\omega$ in $\msD(A)$, hence the above expression simplifies to $g\omega$. 
\end{proof} 

\begin{Def}
We define $\mcD(A) = A\check{A} \subseteq \msD(A)$.
\end{Def}

Clearly $\check{A}\subseteq \msU$, so $\mcD(A)$ has meaning as a vector subspace of $\mcD(A)$. 

\begin{Lem}\label{LemMultBij}
\begin{enumerate}
\item The vector space $\mcD(A)$ is a subalgebra of $\msD(A)$. 
\item The multiplication maps 
\begin{equation}\label{EqDecompDrinf1}
\oplus_{a,b} \left(\Gr{A}{}{a}{}{b}\otimes \Gr{\check{A}}{a}{}{b}{}\right) \rightarrow \mcD(A),\qquad a\otimes \omega \mapsto a\omega,
\end{equation}
\begin{equation}\label{EqDecompDrinf2}
\oplus_{a,b} \left(\Gr{\check{A}}{}{a}{}{b}\otimes \Gr{A}{a}{}{b}{}\right) \rightarrow \mcD(A),\qquad \omega\otimes a \mapsto \omega a
\end{equation}
\end{enumerate}
are isomorphisms.
\end{Lem}
\begin{proof}
Since we can write 
\[
\check{A} = \{\varphi(g \emdash h) \mid g,h\in A\}
\]
by idempotency of $A$ and the existence of the modular automorphism for $\varphi$, it is clear from Lemma \ref{LemMultOppOrd} that both products land in $\mcD(A)$, and that $\mcD(A)$ is a subalgebra. We note in passing that it is also easily seen that $\Gr{A}{r}{a}{t}{b}\cdot \Gr{\check{A}}{a'}{s}{b'}{u}= \Gr{\check{A}}{r}{a}{t}{b}\cdot \Gr{A}{a'}{s}{b'}{u}=0$ if $a\neq a'$ or $b\neq b'$. 

The first multiplication map is then surjective by definition. 
By  Lemma \ref{LemMultOppOrd} and the range properties of the canonical maps (Definition \ref{DefCanMaps}), the second multiplication map contains in its range all $\omega(\UnitC{u}{t} \emdash \UnitC{r}{s})\UnitC{s}{t}a$ for $\omega \in \check{A}$ and $a\in A$, which by the relation $\UnitC{x}{y} = \UnitCC{x}{y}$ and \eqref{EqPartUnitDual} reduces to $\UnitCC{r}{u}\omega \UnitC{s}{t} a$. As $r,u,s,t$ were arbitrary, the surjectivity of the second multiplication map follows. 

Finally, to see that the multiplication maps are injective, note that we can define on $\oplus_{rs} \Gr{A}{}{r}{}{s}\otimes \Gr{\check{A}}{r}{}{s}{}$ a left $\msD(A)$-module structure by 
\[
\pi(a)(b\otimes \chi) = ab\otimes \chi,\qquad \pi(\omega)(b\otimes \chi) =  b_{(2)} \otimes \omega(b_{(3)}\emdash S^{-1}(b_{(1)}))\chi,
\]
where the second formula is meaningful by the same reasons as before and the fact that $\msU \subseteq M(\check{A})$. It is straightforward that this indeed defines a module for $\msD(A)$. If now $\sum_{i} a_i \otimes \omega_i \in \oplus\left(\Gr{A}{r}{a}{t}{b}\otimes \Gr{\check{A}}{a}{s}{b}{u}\right)$ and $\sum_i \pi(a_i)\pi(\omega_i) =0$, then it follows that 
\[
\sum_i  \omega_i(b_{(3)}cS^{-1}(b_{(1)}))a_ib_{(2)}d = 0,\qquad \forall b,c,d\in A. 
\]
By the properties of the canonical maps in Definition \ref{DefCanMaps}, this is the same as having 
\[
\sum_i a_i\UnitC{s}{t} \otimes \omega_i(\UnitC{t}{u}\emdash \UnitC{s}{r}) = 0,\qquad \forall r,s,t,u.
\]
But by the specific form of the subspace of $A\otimes \check{A}$ in which our element lies, this can only happen if $\sum_i a_i \otimes \omega_i =0$. Injectivity of the second multiplication map can be proven similarly.
\end{proof}

\begin{Lem}
The algebra $\mcD(A)$ is a non-degenerate and idempotent $I^2$-algebra by the grading
 \[
\Gr{\mcD(A)}{r}{s}{t}{u} = \UnitC{r}{t}\mcD(A)\UnitC{s}{u}.
\]
\end{Lem}
\begin{proof}
The lemma will follow at once if we show that $\mcD(A) = \sum_{r,s,t,u} \UnitC{r}{t}\mcD(A)\UnitC{s}{u}$. But this is immediate from the fact that (by \eqref{EqDecompDrinf1} and \eqref{EqDecompDrinf2}) $\mcD(A) = A\check{A} = \check{A}A \subseteq \msD(A)$.
\end{proof} 

Note also that from the above, we obtain non-degenerate embeddings $A \rightarrow M(\mcD(A))$ and $\check{A}\rightarrow M(\mcD(A))$.

Let us now start constructing a coproduct on $\mcD(A)$. Note first that by Lemma \ref{LemMultBij}, we may uniquely write elements in $\mcD(A)$ in the form $\sum a_i \omega_i$ for $a_i \otimes \omega_i \in \oplus_{ab} \Gr{A}{}{a}{}{b}\otimes \Gr{\check{A}}{a}{}{b}{}$, or in the form $\sum \omega_i a_i$ for $\omega_i \otimes a_i \in  \oplus_{ab} \Gr{\check{A}}{}{a}{}{b}\otimes \Gr{A}{a}{}{b}{}$. We can hence define a unique multiplier $\Delta(a)$ of $\mcD(A)\otimes^I \mcD(A)$ by putting 
\begin{equation}\label{EqDefDelDrin1}
\Delta(a)(b\chi\otimes c\theta) = (\Delta(a)(b\otimes c))(\chi \otimes \theta),\qquad
(\chi b\otimes \theta c)\Delta(a) = (\chi\otimes\theta)((b\otimes c)\Delta(a)).
\end{equation}
Similarly, we can define $\check{\Delta}(\omega)$ for $\omega \in \msU$ as a multiplier on $\mcD(A)\otimes^I \mcD(A)$ by 
\begin{equation}\label{EqDefDelDrin2}
 \check{\Delta}(\omega)(\chi b \otimes \theta c) = (\check{\Delta}(\omega)(\chi\otimes \theta))(b\otimes c),\qquad (b\chi \otimes c\theta)\check{\Delta}(\omega) = (b\otimes c)((\chi\otimes \theta)\check{\Delta}(\omega).
\end{equation}
Here we use that $\msU \subseteq M(\check{A})$, and that $\check{\Delta}$ can be uniquely extended as a strictly continuous homomorphism to $M(\check{A})$. 

% Next Theorem might need some more details
\begin{Theorem}
Let $(A,\Delta)$ be an $I$-partial algebraic quantum group of compact type, with faithful left invariant functional $\varphi$ and faithful right invariant functional $\psi$. 

Define
\begin{equation}\label{EqDefDelDrin}
\Delta_{\mcD}: \mcD(A) \rightarrow M(\mcD(A) \otimes^I \mcD(A)),\qquad a\omega \mapsto \Delta(a)\check{\Delta}(\omega)
\end{equation}
by means of \eqref{EqDefDelDrin1}. 

Then $(\mcD(A),\Delta_{\mcD})$ defines an $I$-partial algebraic quantum group, with antipode given by 
\begin{equation}\label{EqDefAntDD}
S_{\mcD}(a\omega) = \check{S}(\omega)S(a),\qquad S_{\mcD}(\omega a) = S(a)\check{S}(\omega) \end{equation}
and with faithful left, resp.\ right invariant functionals given (with respect to the decompositions in \eqref{EqDecompDrinf1} and \eqref{EqDecompDrinf2}) by 
\begin{equation}\label{EqInvLeftFuncDrin}
\varphi_{\mcD}: \mcD(A) \rightarrow \C,\qquad a\omega \mapsto \varphi(a)\check{\varphi}(\omega),\qquad  a\in \Gr{A}{}{r}{}{s},\omega \in \Gr{\check{A}}{r}{}{s}{},
\end{equation}
\begin{equation}\label{EqInvRightFuncDrin}
\psi_{\mcD}: \mcD(A) \rightarrow \C,\qquad a\omega \mapsto \psi(a)\check{\psi}(\omega),  \qquad a\in \Gr{A}{}{r}{}{s},\omega \in \Gr{\check{A}}{r}{}{s}{},
\end{equation}
\end{Theorem}
%Note that when specifying the range of $\Delta_{\mcD}$, we are using that for $B$ a non-degenerate algebra, we can view $M(B) \subseteq L(B)$, with $L(B)$ the algebra of left multipliers. 
\begin{proof}
Let us first show that $\Delta_{\mcD}$ is well-defined as a homomorphism $\mcD(A) \rightarrow M(\mcD(A) \otimes \mcD(A))$.

First of all, note that since $\mcD(A) = A\check{A} = \check{A}A$, it is clear that $\mcD(A)\otimes^I \mcD(A)$ is a left unital module for $A$ and for $\check{A}$, where $a\in A$ acts by left multiplication with $\Delta(a)$ and where $\omega \in \check{A}$ acts by left multiplication with $\check{\Delta}(\omega)$. These module structures then extend canonically to module structures for $M(A)$ and $M(\check{A})$, and it is easy to see that then 
\[
\Delta(\UnitC{r}{s}) = \check{\Delta}(\UnitCC{r}{s}) = \sum_{t} \UnitC{r}{t}\otimes \UnitC{t}{s} = \sum_t \UnitCC{r}{t}\otimes \UnitCC{t}{s},
\]
the latter series converging pointwise. 

To show now that $\Delta_{\mcD}$ exists as a homomorphism 
\[
\Delta_{\mcD}: \mcD(A) \rightarrow M(\mcD(A)\otimes^I \mcD(A))\subseteq M(\mcD(A)\otimes \mcD(A)),
\]
it is by Lemma \ref{LemMultBij} sufficient to show that if $a\in \Gr{A}{r}{}{s}{}$ and $\omega \in \Gr{\check{A}}{}{r}{}{s}$ and 
\[
\sum  a' \otimes \omega' = a_{(2)}\otimes \omega(a_{(3)}\emdash S^{-1}(a_{(1)})) \quad \in \oplus_{r',s'} \Gr{A}{}{r'}{}{s'} \otimes \Gr{\check{A}}{r'}{}{s'}{}, \]
then
\[
 \check{\Delta}(\omega) \Delta(a)= \sum \Delta(a')\check{\Delta}(\omega'), 
\]
or hence that for all $\chi \in \check{A}$ and $b\in A$ we have in $\mcD(A)\otimes \mcD(A)$ that
\[
 \check{\Delta}(\omega) \Delta(a)(\chi \otimes b)= \sum \Delta(a')\check{\Delta}(\omega')(\chi \otimes b).
\]
Now applying the commutation rules \eqref{EqInterRel} and \eqref{EqInterRelRev} of $\msD(A)$, and writing $m: \check{A}\otimes A \rightarrow \mcD(A)$ for the multiplication map, it is sufficient to check that
\[
(m\otimes m)X = (m\otimes m)(Y), 
\]
where 
\[
X = \omega_{(1)}\theta(S^{-1}(a_{(3)})\emdash a_{(1)}) \otimes a_{(2)} \otimes \omega_{(2)}\otimes a_{(4)}b,
\]
\[
Y =\sum (\omega_{(1)}'\theta)(S^{-1}(a_{(3)}')\emdash a_{(1)}') \otimes a_{(2)}' \otimes  \omega_{(2)}'(S^{-1}(a_{(6)}')\emdash a_{(4)}')\otimes a_{(5)}'b.
\]
Take now $c,d\in A$. Evaluating the first leg of $X$ in $c\in A$ and the third leg in $d\in A$, we obtain
\begin{equation}\label{EqCompProdDD}
\omega(dc_{(1)}) \theta(S^{-1}(a_{(3)})c_{(2)}a_{(1)}) a_{(2)}\otimes a_{(4)}b.
\end{equation}
On the other hand, evaluating the first leg of $Y$ in $c$ and the third leg in $d$, we obtain
\begin{eqnarray*}
&& \hspace{-2cm} \sum \omega'(S^{-1}(a_{(8)}') d a_{(6)}' S^{-1}(a_{(5)}') c_{(1)}a_{(1)}') \theta(S^{-1}(a_{(4)}')c_{(2)} a_{(2)}') a_{(3)}'\otimes a_{(7)}'b \\
&\underset{\eqref{EqAntDef2}}{=}& 
\sum \sum_t \omega'(S^{-1}(a_{(7)}') d \Unit_t c_{(1)} a_{(1)}') \varepsilon(\Unit^t a_{(5)}') \theta(S^{-1}(a_{(4)}') c_{(2)} a_{(2)}') a_{(3)}'\otimes a_{(6)}'b \\
&=& \sum \sum_t\omega'(S^{-1}(a_{(6)}') d \Unit_t c_{(1)} a_{(1)}')  \theta(S^{-1}(\Unit_ta_{(4)}') c_{(2)} a_{(2)}') a_{(3)}'\otimes a_{(5)}'b \\
&=& \sum \omega'(S^{-1}(a_{(6)}') d  c_{(1)} a_{(1)}')  \theta(S^{-1}(a_{(4)}') c_{(2)} a_{(2)}') a_{(3)}'\otimes a_{(5)}'b \\
&=& \omega(a_{(8)} S^{-1}(a_{(7)})dc_{(1)} a_{(2)}S^{-1}(a_{(1)})) \theta(S^{-1}(a_{(5)}) c_{(2)}a_{(3)}) a_{(4)}\otimes a_{(6)}b \\
&\underset{\eqref{EqAntDef2}}{=}& \sum_{t,u} \omega(\Unit_tdc_{(1)}\Unit_u) \theta(S^{-1}(a_{(3)}) c_{(2)}\Unit^u a_{(1)}) a_{(2)}\otimes \Unit_t a_{(4)}b \\
&=& \omega(\Unit_sdc_{(1)}) \theta(S^{-1}(a_{(3)}) c_{(2)} a_{(1)}) a_{(2)}\otimes \Unit_s a_{(4)}b\\
&=& \omega(dc_{(1)}) \theta(S^{-1}(a_{(3)}) c_{(2)} a_{(1)}) a_{(2)}\otimes a_{(4)}b.
\end{eqnarray*} 
Comparing with \eqref{EqCompProdDD}, we see that indeed $(m\otimes m)(Y) =  (m\otimes m)(X)$. 

Using now again that $\mcD(A) = A\check{A} = \check{A}A$, it is clear that $\Delta_{\mcD}(\mcD(A)) \in \widetilde{M}(\mcD(A)\otimes^I \mcD(A))$, and that $\Delta_{\mcD}$ defines a non-degenerate homomorphism $\mcD(A) \rightarrow M(\mcD(A)\otimes^I \mcD(A))$. Also the coassociativity of $\Delta_{\mcD}$ is immediate, so $\Delta_{\mcD}$ is a regular comultiplication. 

%It is also clear that  with then 
%\[
%\Delta_{\mcD}(1_{\mcD}) = \Delta_A(1_A) = %\Delta_{\check{A}}(\check{1}_{\check{A}}) \in %M(\mcD(A)\otimes \mcD(A)). 
%\]

It is now further routine to check that $(\mcD(A),\Delta_{\mcD})$ becomes a regular $I$-partial multiplier bialgebra with counit given explicitly by 
\[
\varepsilon_{\mcD}(a\omega) = \varepsilon(a)\check{\varepsilon}(\omega),\qquad a\in \Gr{A}{}{r}{}{s},\omega \in \Gr{\check{A}}{r}{}{s}{},
\]
\[
\varepsilon_{\mcD}(\omega a)  = \check{\varepsilon}(\omega)\varepsilon(a), \qquad \omega \in \Gr{\check{A}}{}{r}{}{s},a\in \Gr{A}{r}{}{s}{}.  
\]
In fact, $(\mcD(A),\Delta_{\mcD})$ is a regular $I$-partial multiplier Hopf algebra with antipode map as in \eqref{EqDefAntDD}.

The invariance conditions for $\varphi_{\mcD}$ and $\psi_{\mcD}$ are immediately clear, as is their faithfulness. 
Hence $(\mcD(A),\Delta_{\mcD})$ is an $I$-partial algebraic quantum group. 
\end{proof}
\begin{Def}
We call $(\mcD(A),\Delta_{\mcD})$ the \emph{Drinfeld double} of $(A,\Delta)$. 
\end{Def}

In the next proposition, we determine the modular data associated to $\mcD(A)$. 

\begin{Prop}\label{PropModDatDrinDoub}
Let $(A,\Delta)$ be an $I$-partial algebraic quantum group of compact type, with faithful left invariant functional $\varphi$ and faithful right invariant functional $\psi$. 

Then the modular data with respect to \eqref{EqInvLeftFuncDrin} is given by 
\begin{equation}\label{EqModStructDrinDoub}
\sigma^{\varphi_{\mcD}}(a) = S^{2}(a),\qquad \sigma^{\varphi_{\mcD}}(\omega) = \check{S}^2(\omega),\qquad \delta_{\varphi_{\mcD}} = \check{\delta}_{\check{\varphi}}\delta_{\varphi} \nu\qquad \nu_{\mcD} = 1, 
\end{equation}
where $a\in A,\omega \in \check{A}$. 
\end{Prop}
\begin{proof}
Take $a \in \Gr{A}{r}{s}{p}{q}, b\in \Gr{A}{x}{z}{y}{l}$ and $\omega \in \Gr{\check{A}}{s}{x}{q}{y}$. Then by the definition of the product in $\mcD(A)$, we obtain 
\[
\varphi_{\mcD}(a\omega b) = \varphi(ab_{(2)}) \check{\varphi}(\omega(b_{(3)}\emdash S^{-1}(b_{(1)}))). 
\]
Writing $\omega = \varphi(-d)$ for $d\in \Gr{A}{s}{q}{x}{y}$, we find that 
\[
\varphi_{\mcD}(a\omega b) = \varphi(ab_{(2)}) \varepsilon(S^{-1}(b_{(1)}) d \sigma^{\varphi}(b_{(3)})).
\]
Now the indices of the components $S^{-1}(b_{(1)})$, $d$ and $\sigma^{\varphi}(b_{(3)})$ are such that we can apply the multiplicative property of $\varepsilon$, so we can simplify 
\[
\varphi_{\mcD}(a\omega b) = \varphi(ab_{(2)}) \varepsilon(S^{-1}(b_{(1)})) \varepsilon(d) \varepsilon(\sigma^{\varphi}(b_{(3)})))
\]
(note that this expression is meaningful as the factor $\varphi(ab_{(2)})$ caps of the triple comultiplication of $b$ in all factors). Using now the counit property, we can further simplify 
\[
\varphi_{\mcD}(a\omega b) = \varphi(ab_{(1)}) \varepsilon(\sigma^{\varphi}(b_{(2)})) \check{\varphi}(\omega).
\]
Finally, by \eqref{EqFormCoprodMod} we see that this becomes
\[
\varphi_{\mcD}(a\omega b) = \varphi((\sigma^{\varphi})^{-1}S^{-2}\sigma^{\varphi}(b)a) \check{\varphi}(\omega) = \varphi_{\mcD}(S^{-2}(b)a \omega),
\]
so by Proposition \ref{PropSModAut} we obtain the formula for $\sigma^{\varphi_{\mcD}}(b)$ in \eqref{EqModStructDrinDoub}.

We can do a similar computation to find how the modular structure of the Drinfeld double acts on $\check{A}$. Namely, take $\chi,\omega \in \check{A}$ and $a \in A$. Then using the invariance and the $\delta$-invariance of $\varphi$ from \eqref{EqModEl}, we find 
\begin{eqnarray*}
\varphi_{\mcD}(\chi a \omega) &=& \varphi_{\mcD}(a_{(2)} \chi(a_{(3)}\emdash S^{-1}(a_{(1)}))\omega) \\
&=& \varphi(a_{(2)}) \check{\varphi}(\chi(a_{(3)}\emdash S^{-1}(a_{(1)}))\omega)\\
&=& \sum_r \varphi(a_{(1)}\Unit^r) \check{\varphi}(\chi(a_{(2)}\emdash \Unit_r)\omega)\\
&=& \sum_r \varphi((a\Unit^r)_{(1)}) \check{\varphi}(\chi((a\Unit^r)_{(2)}\emdash \Unit_r)\omega)\\
&=& \sum_{r,s} \varphi(a\UnitC{r}{s})\check{\varphi}(\chi(\delta_{\varphi}\Unit_s \emdash \Unit_r)\omega)\\
&=& \sum_{r,s} \varphi(a\UnitC{r}{s})\check{\varphi}(\chi(\delta_{\varphi}\emdash)\UnitCC{r}{s}\omega) \\
&=& \varphi_{\mcD}(a \omega \sigma^{\check{\varphi}}(\chi(\delta_{\varphi}\emdash))). 
\end{eqnarray*}
Using now Proposition \ref{PropDualModStr}, we find from \eqref{EqPropdelt} that $\sigma^{\varphi_{\mcD}}(\chi) = \check{S}^2(\chi)$. 

Note that from the above formulas for the modular structure, we immediately obtain that also
\[
\varphi_{\mcD}(\omega a) = \varphi(a)\check{\varphi}(\omega),\qquad a\in \Gr{A}{}{r}{}{s},\omega \in \Gr{\check{A}}{r}{}{s}{}. 
\]
It follows that 
\begin{multline*}
\varphi_{S_{\mcD}}(a\omega) = \varphi_{\mcD}(\check{S}(\omega)S(a)) = \varphi(a\delta_{\varphi}) \check{\varphi}(\omega \check{\delta}_{\check{\varphi}}) =\varphi(\delta_{\varphi}\nu a) \check{\varphi}(\omega \check{\delta}_{\check{\varphi}}) 
 =  \varphi_{\mcD}(\delta_{\varphi} \nu a \omega \check{\delta}_{\check{\varphi}}) =  \varphi_{\mcD}(a \omega \check{\delta}_{\check{\varphi}}\delta_{\varphi} \nu)
\end{multline*}

The fact that $\nu_{\mcD} =1$ follows now immediately by applying $\sigma^{\varphi_{\mcD}}$ to $\delta_{\varphi_{\mcD}}$. 
\end{proof}

Assume now that $(A,\Delta)$ defines a partial $*$-algebraic quantum group of compact type, i.e.\ a partial compact quantum group, and let $\varphi$ be a normalized faithful positive invariant functional. 

\begin{Theorem}
Let $(A,\Delta)$ be an $I$-partial $*$-algebraic quantum group of compact type, and let $(\mcD(A),\Delta_{\mcD})$ be its Drinfeld double as an $I$-partial algebraic quantum group. Then $(\mcD(A),\Delta_{\mcD})$ is an $I$-partial $*$-algebraic quantum group for the $*$-structure
\[
(a\omega)^* = \omega^*a^*,\qquad (\omega a)^* = a^*\omega^*. 
\]
\end{Theorem} 
\begin{proof}
It follows from straightforward calculations that the $*$-structure is in fact well-defined on $\msD(A)$, and that it in particular turns $\mcD(A)$ into an $I$-partial multiplier Hopf $*$-algebra. 

To see that the associated left invariant functional is positive, note that by Proposition \ref{PropModDatDrinDoub} (and the fact that $\varphi = \varphi \circ S$ and $\nu= 1$ in our situation) we have 

\begin{equation}\label{EqPDFormD}
\varphi_{\mcD}((\omega a)^*\chi b) = \check{\varphi}(\omega^*\chi) \varphi(S^{-2}(b)a^*) = \check{\varphi}(\omega^*\chi) \varphi(S^{-1}(a)^*S^{-1}(b)).
\end{equation}
Hence the left hand side of \eqref{EqPDFormD} defines a positive definite sesqui-linear form on $\mcD(A)$, so $\varphi_{\mcD}$ is a positive faithful left invariant functional.  

%The proof that $\psi_{\mcD}$ is a positive faithful right invariant functional is completely similar. 
\end{proof}

% need to introduce Sweedler notation for comodules.

\subsection{Yetter-Drinfeld modules} 
% Add extra binding text

In this last section, we show that modules of the Drinfeld double of a partial algebraic quantum group of compact type can be described in terms of Yetter-Drinfeld modules of the underlying partial Hopf algebra, as described by the next definition. 

\begin{Def}
Let $(A,\Delta)$ be a partial Hopf algebra. A (left) \emph{Yetter-Drinfeld module} for $(A,\Delta)$ consists of an $I^2$-graded vector space $V$ with a left $A$-module and right $A$-comodule structure $\delta$ (compatible with the associated horizontal and vertical gradings) such that for all $a\in A,v\in V$
\begin{equation}\label{EqIdYD}
\Gr{\delta}{}{}{r}{s}(av) =\sum_{pq} \UnitC{s}{r}a_{(2)}v_{(0;pq)} \otimes \mbunit^ra_{(3)}v_{(1;pq)}S^{-1}(a_{(1)})\mbunit^s. 
\end{equation}
If $(A,\Delta)$ is a partial Hopf $*$-algebra, we call the Yetter-Drinfeld module \emph{unitary} if $V$ is an $I^2$-graded pre-Hilbert space and the $A$-module and $A$-comodule structure are unitary. 
\end{Def}
Note that \eqref{EqIdYD} is meaningful, since we can write the right hand side as 
\[
\sum_{pq} (\UnitC{s}{r}a_{(2)}\otimes  \UnitC{r}{m} a_{(3)})\Gr{\delta}{}{}{p}{q}(v)(1\otimes S^{-1}(\UnitC{k}{s}a_{(1)})),\qquad a\in \Gr{A}{k}{l}{m}{n},
\]
where all except a finite number of terms in the sum disappear.

% Next proposition a bit quick... Also: add unitary structure
\begin{Prop}
Let $V$ be an $I^2$-graded vector space, and let $(A,\Delta)$ be a partial algebraic quantum group of compact type. Then there is a one-to-one correspondence between Yetter Drinfeld $A$-module structures on $V$ and left $\mcD(A)$-module structures on $V$ by 
\[
(a\omega)v = a(\id\otimes \omega)\delta(v),\qquad a\in A,\omega \in\check{A}.
\]
Furthermore, if $(A,\Delta)$ defines a partial compact quantum group, the above correspondence induces a one-to-one correspondence between unitary Yetter-Drinfeld modules and unitary $\mcD(A)$-modules. 
\end{Prop} 
Here we write 
\begin{equation}\label{EqDefActOm}
(\id\otimes \omega)\delta(v) = \sum_{rs} (\id\otimes \omega)\Gr{\delta}{}{}{r}{s}(v),
\end{equation}
    which is allowed as all but a finite number of terms on the right hand side disappear.
\begin{proof}
By Proposition \ref{Prop1to1RepCorep}, we have a one-to-one correspondence between right $A$-comodules and left $A$-modules. 

If $(V,\delta)$ is a Yetter-Drinfeld module for $(A,\Delta)$, then the associated $\check{A}$-module will satisfy $\UnitC{r}{s} = \UnitCC{r}{s}$ by the compatibility of the gradings. Furthermore, we can extend the left $\check{A}$-module to a unital $\msU$-module by the same formula \eqref{EqDefActOm}, which is still meaningful by the defining condition on $\msU$. The identity \eqref{EqIdYD} then turns into the defining interchange relation \eqref{EqInterRel}, so we obtain a unital $\msD(A)$-module, which necessarily restricts (by its construction) to a non-degenerate $\mcD(A)$-module. 

Conversely, if we have a non-degenerate $\mcD(A)$-module $V$, then through the natural map 
\[
\check{A} \rightarrow \mcD(A),\qquad \Gr{\check{A}}{r}{s}{t}{u}\ni \omega \mapsto \UnitC{r}{t}\omega,
\] 
we obtain that $V$ becomes a non-degenerate left $\check{A}$-module. Let $(V,\delta)$ be the dual $A$-comodule. Then it is again clear that the gradings on $V$ are compatible, and \eqref{EqIdYD} must hold since it holds true when applying any element of $\check{A}$ to the second leg, and these elements separate $A$.

In case $(A,\Delta)$ defines a partial compact quantum group, it is clear that the above correspondence preservers unitarity, again by Proposition \ref{Prop1to1RepCorep}.
\end{proof}

\end{document}